%% file: main.tex
\documentclass{amsart}
\usepackage[english]{babel}
\usepackage{graphicx}
\usepackage[hidelinks]{hyperref} 
\usepackage{amssymb}
\usepackage{amsmath}
\usepackage[foot]{amsaddr}
\usepackage{amsfonts}
\usepackage{xcolor}
\usepackage{tikz}
\usepackage{tikz-3dplot}
\usepackage{siunitx}
\usepackage{pgfplots}
\usetikzlibrary{decorations.pathmorphing} 
\usepackage{mathtools}
\usepackage{algpseudocode}
\usepackage{algorithm}
\usepackage[top=1in, bottom=1.25in, left=1.25in, right=1.25in]{geometry}
\usepackage{subcaption}

\newtheorem*{remark}{Remark}

\graphicspath{{images/}}
\usepackage{float}

\pgfplotsset{compat=1.18}
\begin{document}
\nocite{*}
\title[Detecting bifurcations for multiparametric PDEs]{Bifurcation Curve Detection with Deflation for Multiparametric PDEs}
\author{Nitin Kumar$^1$}
\author{Federico Pichi$^1$}
\author{Gianluigi Rozza$^1$}

\address{$^1$ mathLab, Mathematics Area, SISSA, via Bonomea 265, I-34136 Trieste, Italy}

\input{abstract}

\maketitle
\tableofcontents

\input{sections/intro}
\input{sections/nonLinear_pde}
\input{sections/numerical_approximation}

\input{sections/results}
\input{sections/Conclusion}
\input{sections/acknowledgments}
\input{sections/Appendix}


\bibliographystyle{abbrv}
\bibliography{main}

\end{document}

%% file: abstract.tex
\begin{abstract}
This work presents a comprehensive framework for capturing bifurcating phenomena and detecting bifurcation curves in nonlinear multiparametric partial differential equations, where the system exhibits multiple coexisting solutions for given values of the parameters. Traditional continuation methods for one-dimensional parameterizations employ the previously computed solution as the initial guess for the next parameter value. These are usually very inefficient, since small step sizes increase computational cost, while larger steps could jeopardize the method convergence jumping to a different solution branch or missing the bifurcation point. To address these challenges, we propose a novel framework that combines: (i) arclength continuation, adaptively selecting new parameter values in higher dimension, and (ii) the deflation technique, discovering multiple branches to construct complete bifurcation diagrams without requiring a costly spectral analysis of the system. In particular, the arclength continuation method is designed to handle multiparametric scenarios, where the parameter vector $\lambda \in \mathbb{R}^p$ traces a curve $g(\lambda)$ within a $p$-dimensional parameter space. In addition, we introduce a zigzag path-following strategy to robustly track the bifurcation curves and surfaces, respectively, for two- and three-dimensional parametric spaces. Finally, we demonstrate its performance on three benchmark problems of increasing complexity: from the 1D/2D Bratu and Allen--Cahn equations to the 2D/3D Rayleigh--Benard convection problem.
\end{abstract}

%% file: sections/intro.tex
\section{Introduction and motivation}
Nonlinear parametric partial differential equations (PDEs) play a crucial role in many scientific and industrial applications, as they provide reliable models for capturing complex physical phenomena in many different physically and geometrically parameterized contexts in continuum mechanics, quantum mechanics, and fluid dynamics. 

In contrast to linear parametric PDEs, which typically exhibit a unique solution that evolves continuously with small changes in the parameter $\lambda \in \mathbb{R}$, nonlinear parametric PDEs may admit multiple solutions for the same parameter value. As a result, a slight variation in the parameter can lead to abrupt changes in the solution's behavior and stability. This phenomenon is known as bifurcation \cite{seydel2009practical,caloz1997numerical,kielhofer2012bifurcation}, and the parameter value at which such qualitative changes occur is referred to as the \textit{bifurcation point} $\lambda^*$.

The complex and non-differentiable dependence of the solution on the parameters often results in the non-uniqueness of parametric PDE solutions. Several physical models exhibit this behavior, including the buckling for the Von Kármán plate model \cite{bauer1965nonlinear,berger1967karman,pichi2019reduced} and hyperelastic beams \cite{PichiReducedOrderModels2024,xia2020nonlinear,yu2022analytic}, the Navier–Stokes equations in channel flow problems \cite{pintore2021efficient,pitton2017computational,pichi2022driving,guevel2018numerical,TomadaSparseIdentificationBifurcating2025a,GonnellaStochasticPerturbationApproach2024a}, the Gross–Pitaevskii equation for Bose–Einstein condensates \cite{charalampidis2018computing,pichi2020reduced}, Turing models \cite{choi2025turing,sun2025spatiotemporal,woolley2025bespoke}, post-contact states in Micro-Electromechanical Systems \cite{naudet2024numerical,gobat2022reduced,opreni2021model} and many more ranging from fluid dynamics \cite{pitton2017application,KhamlichModelOrderReduction2022,deng2020low,venturi_stochastic_2010,boulle2022bifurcation} to chemical reactions \cite{olsen1991bifurcation}.

Various bifurcating branches may emerge from the critical value at which the system undergoes a qualitative change in behavior. A common and intuitive way to represent these changes is through a \textit{bifurcation diagram}, which plots a scalar quantity of interest for the solution against the corresponding parameter value. For example, constructing these plots one could observe that, when crossing a critical parameter value $\lambda^*$, a single solution branch gives rise to two additional symmetric branches of admissible solutions that coexist with the original one, while exchanging stability properties. This type of behavior is characteristic of a \textit{pitchfork bifurcation}. For instance, the Allen--Cahn equation \cite{kuehn2024uncertainty,allen1982phase,feng2003numerical}, commonly used to model phase separation in multi-component alloys, exhibits multiple pitchfork bifurcation.

In contrast, if two solution branches merge at a single point in the parameter space, beyond which no solution exists, the model is said to exhibit a \textit{saddle-node bifurcation}. The Bratu equation \cite{syam2006efficient,shahab2024neural,hosseini2023application}, which arises in combustion theory, features a saddle-node bifurcation, where two solution branches merge and annihilate. 

Beyond these primary bifurcations, many nonlinear PDEs exhibit multiple or secondary bifurcations, leading to increasingly intricate bifurcation structures. Indeed, many models display a rich hierarchy of bifurcating branches associated with qualitatively different states. Moreover, secondary bifurcations can emerge from nontrivial branches, giving rise to an intricate cascade of bifurcations from already bifurcated patterns. These examples illustrate how non-linear parametric PDEs can exhibit rich and complex solution structures depending on their parametric characterization.

The numerical computation of bifurcation diagrams involves three distinct tasks. First, locating the bifurcation points \cite{lopez2002detection,allower1990detection,pichi_deflation-based_2025,pichi_artificial_2023} could require an adaptive parameter strategy, monitoring suitable quantities, such as Jacobian determinant or critical eigenvalues, and the detection of qualitative changes in the solution. Second, once a bifurcation point is identified, it is important to find all coexisting admissible solutions and possibly their stability properties, which can be achieved using techniques such as deflation \cite{farrell2015deflation} or careful initialization strategies \cite{pichi2020reduced,venturi_stochastic_2010}. Finally, when a complete knowledge of the system is acquired, one can adopt a continuation strategy \cite{allgower2003introduction,uecker2022continuation} to reconstruct the solution branches, employing either simple or arclength continuation methods to accurately track the bifurcating behaviors, even through turning points.

However, the simple continuation method comes with limitations: using a very small step size in the parameter space leads to a high computational cost \cite{allgower2003introduction}, while taking large steps may result in convergence failure or cause the solver to jump to an unintended solution branch. To address these issues, we employ the arclength continuation method \cite{chan1982arc,chan1984newton,dickson2007condition}, where the step size is adaptively controlled and the parameter is treated as an unknown. This allows the method to take smaller steps near bifurcation points, where the solution changes rapidly, and larger steps in smoother regions, leading to a more robust and efficient construction of the bifurcation diagram.

On the other hand, deflation can be exploited to reconstruct bifurcation diagrams, since it enables the computation of multiple coexisting solutions for a given parameter value without requiring a spectral analysis, which is computationally expensive for large-scale discretizations of nonlinear PDEs. Moreover, an important advantage of deflation is its ability to identify solutions lying on disconnected branches, which are often difficult or impossible to detect using standard continuation or eigenvector-based branch-switching methods, making the deflation particularly attractive for multiparametric PDE problems.


Indeed, in real scenarios, PDEs are characterized by multiple parameters, resulting in considerably more intricate bifurcation diagrams. In a two-parameter setting, this complexity appears as bifurcation curves instead of isolated bifurcation points, and solution surfaces rather than simple one-dimensional solution branches. With three parameters, the structure becomes even more involved, as bifurcation curves extend into bifurcation surfaces that separate regions of uniqueness or non-existence from those with multiple solutions. 


Accurately detecting bifurcation curves or surfaces in multiparametric settings can be more valuable than computing the entire bifurcation diagram, especially in practical applications. Researchers and engineers are often primarily interested in the regions in which the stability of solutions changes. Indeed, constructing the full bifurcation diagram in a multiparametric space can be computationally unaffordable, while focusing solely on the detection of bifurcation curves significantly reduces the computational burden but still allows capturing the most critical features of the system’s behavior. 


Motivated by the advantages of arclength continuation and deflation, we propose a hybrid strategy that combines deflation with arclength continuation to efficiently compute bifurcation diagrams of nonlinear PDEs in the multiparametric setting. This integrated approach enables robust detection of bifurcation points and accurate construction of solution branches, even in regions with complex solution behavior. We extend the arclength continuation method combined with deflation to provide a systematic framework for exploring efficiently and comprehensively bifurcation diagrams in high-dimensional settings with many parameters. Furthermore, we introduce a zigzag path-following strategy tailored to accurately trace critical bifurcating regions in multiparametric scenarios. This strategy is particularly effective for detecting both saddle-node and pitchfork bifurcations by making small, controlled deviations along turning points, and it naturally extends to three-parameter problems, enabling efficient and reliable detection of bifurcation surfaces.


The main novel contributions of this work are summarized as follows:
\begin{itemize}
\item Deflated arclength continuation algorithm for multiparametric problems, a novel method that combines the deflation technique with arclength continuation to robustly and efficiently construct bifurcation diagrams by capturing multiple solution branches and handling turning points.
\item Detection strategy for bifurcation curves and surfaces, a zigzag path-following approach designed to accurately detect bifurcation curves and surfaces, in two- and three-dimensional parametric problems, respectively, applicable for both saddle-node and pitchfork bifurcating phenomena.
\end{itemize}

The proposed algorithms are validated on three multiparametric problems: the Allen--Cahn equation, which exhibits pitchfork bifurcations, the Bratu equation, which features saddle-node bifurcations, and the Rayleigh--Benard convection problem, which exhibits rich bifurcation structure in the Rayleigh--Prandtl parameter space. The first two examples are studied in one- and two-dimensional spatial domains, while the Rayleigh--Bénard problem is considered in both two- and three-dimensional settings. A key strength of the proposed method lies in its robustness, since it performs reliably without requiring any prior knowledge of the bifurcation type inherent in the multiparametric PDE under consideration. Furthermore, the algorithms effectively construct complete bifurcation diagrams even in the presence of multiple bifurcation points.


In these cases, the detected bifurcation curves are linear, thus, to further demonstrate the robustness and effectiveness of our approach, we designed a custom example that exhibits a nonlinear bifurcation curve, highlighting the algorithm's ability to handle more intricate scenarios. To the best of our knowledge, there is currently no general methodology in the literature capable of both constructing bifurcation diagrams and detecting bifurcation curves in multiparametric context.

The structure of the paper is as follows: Section 2 introduces the nonlinear multiparametric PDE model at both the continuous and discrete levels, with particular emphasis on the non-uniqueness of the solution. Section 3 presents the development of the methodologies for constructing bifurcation diagrams, while Section 4 is devoted to detecting bifurcation curves in multiparametric settings, including an extension of the arclength continuation technique to multiple parameters. In Section 5, the proposed strategies are applied to the multiparametric Allen--Cahn, Bratu and Rayleigh--Benard convection problems to demonstrate their effectiveness. Finally, Section 6 concludes the manuscript and provides potential directions for future works.

%% file: sections/nonlinear_pde.tex
\section{Bifurcating phenomena for nonlinear parametric PDEs}
Nonlinear parametric PDEs often exhibit complex bifurcation phenomena, where small variations in parameters can cause abrupt changes in the number of coexisting solutions or in the stability properties. In such settings, multiple solution branches coexist in regions of the parametric domain, and their identification is essential for theoretical insight, modeling features, and numerical approximation. This section provides an overview on the mathematical description of the bifurcation behaviors that arise in nonlinear PDEs, highlighting the influence of parameter variations and the computational challenges involved in capturing such phenomena \cite{seydel2009practical,kielhofer2012bifurcation,kuznetsov1998elements,chow2012methods}.

Let $\Omega \subset \mathbb{R}^d$ be an open, bounded, and regular domain, where $d$ denotes the spatial dimension, and  $\mathcal{P} \subset \mathbb{R}^p$ the parametric space with $p \geq 1$. Consider a nonlinear parametric operator $G: \mathcal{H} \times \mathcal{P} \rightarrow \mathcal{H}^*$, where $\mathcal{H}$ is a suitable Hilbert space and $\mathcal{H}^*$ its dual. The parametric problem consists of: given a parameter $\lambda \in \mathcal{P}$, find a solution $u = u(\lambda) \in \mathcal{H}$ such that
\begin{equation}\label{sec1eq1}
G(u, \lambda) = 0 \quad \text{in } \mathcal{H}^*.
\end{equation}
To investigate such problems, a numerical approximation is often necessary. Thus, after discretizing the parametric PDE above, one obtains the following nonlinear parametric system, which for a given parameter $\lambda\in \mathcal{P}$ reads as
\begin{equation}\label{sec1eq2}
    G_h(u_h,\lambda)=0,
\end{equation}
where $u_h=u_h(\lambda)\in \mathcal{H}^{N_h}$, with $\mathcal{H}^{N_h}\subset \mathcal{H}$ a finite dimensional subspace of dimension $N_h$. Since the problem exhibits nonlinear terms, the system can be solved using Newton's method \cite{ortega2000iterative,kelley2003solving,knoll2004jacobian,quarteroni2006numerical}.

The existence and uniqueness of the solution of the nonlinear parametric PDEs, which in this setting can be guaranteed by the implicit function theorem \cite{ciarlet2025linear,nussbaum1975global}, are very crucial in the construction of the bifurcation diagram.

Under certain regularity conditions, the solution varies smoothly and remains unique as the parameter changes. However, when these conditions are violated, even slight changes in the parameter near a critical point $\lambda^*$ can lead to abrupt transitions in the system’s response, i.e.\ the system is said to undergo a bifurcation, allowing the coexistence of multiple admissible solutions and the change of their stability properties.

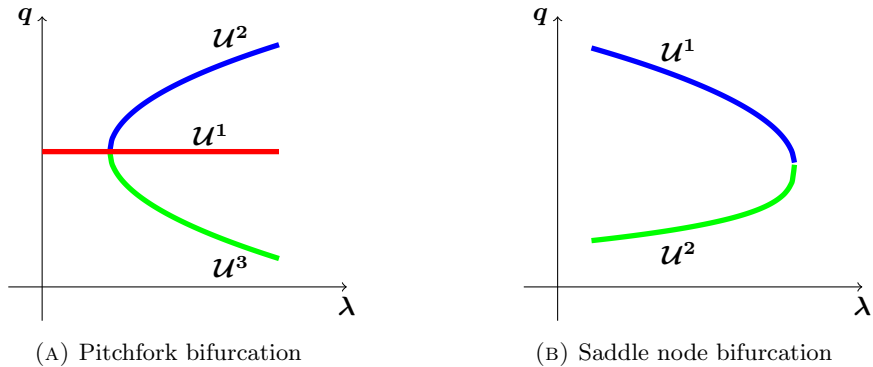
\begin{figure}[htbp]
\centering
\begin{subfigure}[b]{0.28\textwidth}
\centering
\begin{tikzpicture}[scale=0.9]
  \draw[->] (-1.5,-2) -- (3.5,-2) node[below] {$\boldsymbol{\lambda}$};
  \draw[->] (-1,-2.5) -- (-1,2) node[left] {$\boldsymbol{q}$};

  \draw[line width=2pt, blue, domain=0:2.5, samples=100, smooth] plot (\x, {sqrt(\x)});
  \draw[line width=2pt, green, domain=0:2.5, samples=100, smooth] plot (\x, {-sqrt(\x)});
  \draw[line width=2pt, red, domain=-1:2.5, samples=100, smooth] plot (\x, {0});
  
  \node at (1.5, 0.25) {$\boldsymbol{\mathcal{U}^1}$}; 
  \node at (1.8, 1.7) {$\boldsymbol{\mathcal{U}^2}$}; 
  \node at (1.8, -1.7) {$\boldsymbol{\mathcal{U}^3}$};

\end{tikzpicture}
\caption{Pitchfork bifurcation}
\label{fig:pitch}
\end{subfigure}
\quad \hspace{2cm}
\begin{subfigure}[b]{0.28\textwidth}
\centering
\begin{tikzpicture}[scale=0.9]
    \draw[->] (-0.5,0) -- (4.5,0) node[below] {$\boldsymbol{\lambda}$};
  \draw[->] (0,-0.5) -- (0,4) node[left] {$\boldsymbol{q}$};

  \draw[line width=2pt, blue, domain=0.5:3.5, samples=100, smooth] 
    plot (\x, {sqrt(3.5 - \x)+1.8});
  \draw[line width=2pt, green, domain=0.5:3.5, samples=100, smooth] 
    plot (\x, {-pow(3.5 - \x, 0.25) + 2});

\node at (1.8, 3.5) {$\boldsymbol{\mathcal{U}^1}$}; 
  \node at (1.8, 0.5) {$\boldsymbol{\mathcal{U}^2}$};
\end{tikzpicture}
\caption{Saddle node bifurcation}
\label{fig:saddle}
\end{subfigure}
\caption{A schematic representation of bifurcation diagrams.}
\label{fig:Bifurcations}
\end{figure}

We define as solution branches the distinct physical states or behaviors, originating from the bifurcation points, that the system exhibits while varying the parameter values $\lambda \in \mathcal{P}$. Suppose that $\mathcal{K}$ such branches of qualitatively different configurations exist, then we denote by ${\mathcal{U}}^i= \{u_h^i(\lambda): G(u_h^i(\lambda), \lambda) = 0,\lambda\in \mathcal{P}\}$ the $i$-th branch, which contains the set of solutions $u_h^i(\lambda)$. The union of all branches of solutions can be defined as $\mathcal{U}=\cup_{i=1}^{\mathcal{K}} \mathcal{U}^i$.
A schematic representation of ${\mathcal{U}}$ is shown in Figure \ref{fig:Bifurcations} depicting the case of pitchfork and saddle node bifurcations, Figures \ref{fig:pitch} and \ref{fig:saddle}, respectively with $q$ denoting the quantity of interest used to generate the bifurcation diagram.

Capturing bifurcation phenomena presents significant computational challenges, primarily due to the loss of uniqueness and the sensitivity of solutions near bifurcation points. A straightforward application of Newton's method, commonly used for solving nonlinear systems, may perform poorly and diverge in these regions due to near-singular jacobians \cite{hueso2009modified}. 

To address this, advanced strategies such as deflation and continuation methods are employed to explore the bifurcation diagram more effectively. Deflation allows for the computation of multiple distinct solutions at a fixed parameter value, while continuation techniques trace the evolution of the solutions branches as the parameter varies, exploiting previous information as a warm start for the nonlinear solver at next parameter value.

%% file: sections/numerical_approximation.tex
\section{Numerical approximation of bifurcating phenomena}

Several computational methodologies have been developed in the literature for the construction of bifurcation diagrams in nonlinear systems, ranging from classical continuation frameworks to modern machine learning-based approaches. Traditional numerical continuation techniques, as implemented in well-established software such as MATCONT \cite{dhooge2003matcont}, AUTO \cite{doedel1997auto97} and COCO \cite{dankowicz2013recipes} for ordinary differential equations and pde2path \cite{uecker2014pde2path} for nonlinear elliptic PDEs, provide robust tools for detecting bifurcation points and tracking solution branches. These packages rely on predictor–corrector schemes, typically based on simple continuation, and have been extensively applied to study a wide variety of nonlinear dynamical systems. Building upon these foundational frameworks, several recent works have focused on enhancing the robustness and efficiency of continuation methods. For instance, an adaptive homotopy approach was proposed in \cite{hao2020adaptive} to compute bifurcations in nonlinear parametric systems efficiently, while the work \cite{farrell2015deflation} presented a systematic strategy for computing disconnected bifurcation branches. 

More recently, data-driven and neural network-based frameworks have emerged as powerful alternatives to learning and reconstructing bifurcation diagrams directly from the governing equations. The studies in \cite{hao2022learn} and \cite{shahab2024neural} introduced neural network-based approaches that integrate physical constraints or linear stability information to detect and track bifurcation points in parametric PDEs and nonlinear systems. 

Understanding and accurately capturing bifurcation phenomena in nonlinear parametric PDEs requires robust numerical techniques capable of tracing multiple solution branches and detecting qualitative changes in system behavior. In this section, we discuss two key computational methods for our proposed strategy, the arclength continuation and the deflation technique. The former provides a stable path-following mechanism for tracking solution branches, while the latter enables the systematic discovery of multiple distinct solutions. 

Thus, building upon these existing approaches, we propose a unified methodology for the complete reconstruction of complex bifurcation diagrams and the detection of the bifurcation curves.
Furthermore, we extend the methodology to operate in a multiparametric setting by generalizing the arclength continuation method to handle multiple parameters simultaneously. On this foundation, we develop a \textit{zigzag} path-following strategy designed to efficiently trace bifurcation curves and surfaces in multiparametric spaces. A detailed discussion of the proposed framework and its implementation is presented in the following section.

\subsection{Arclength continuation method}\label{sec31}
In this section, we discuss the arclength strategy as the continuation technique exploited for the numerical approximation of the branches \cite{dickson2007condition,de1999determination,mittelmann1986pseudo,dahlke2024practical}.
Let us first assume, to simplify the discussion, that the parametric system in Equation \eqref{sec1eq1} admits a solution $u$ for any value of the parameter $\lambda\in \mathbb{R}$. 
Later, we will develop an arclength continuation method in a more general and multiparametric setting where $\lambda \in \mathbb{R}^p$ with $p>1$. It is important to note that, in practice, we employ the pseudo-arclength continuation method, which provides a numerically efficient approximation of the true arclength formulation. Pseudo-arclength method replaces the geometric arclength constraint by a tangent-based normalization condition, providing a first-order accurate approximation of true arclength parametrization. The resulting augmented system remains well-conditioned and requires only first-order derivative information, offering a robust and computationally efficient alternative \cite{keller1977numerical,chan1982arc}. 

Arclength continuation overcomes the issue related to choosing the step size for the parameter exploration. The main idea behind this technique is to consider the next value of the parameter $\lambda_{i+1}$ as an unknown, and an alternative parameterization of the branch is taken into account. Hence, instead of parameterizing the solution $u$ by $\lambda$, the branches are parameterized using the \textit{arclength} parameter $s$, indicating the distance along the current solution branch we want to travel. 

%


To be more specific, here we consider the arclength $s$ as the parameter, and treat $u$ and $\lambda$ as the function of $s$, i.e., we define $x(s)=(u(s),\lambda(s))$ which depends smoothly on $s$. Now, differentiating the parametric system in Equation \eqref{sec1eq1} with respect to $s$ one obtains
\begin{equation}\label{sec3.1eq1}
    \frac{dG(u(s),\lambda(s))}{ds}=G_u \dot{u}+G_{\lambda} \dot{\lambda}=0.
\end{equation}
Denoting by $\dot{x}$ the derivative with respect to the arclength parameter $s$, it corresponds to the unit tangent vector $(\dot{u},\dot{\lambda})$ at the point $(u,\lambda)$ on the solution branch with 
\begin{equation}\label{sec3.1eq2}
    \|\dot{x}\|^2=\|\dot{u}\|^2+|\dot{\lambda}|^2=1.
\end{equation}
It should be noted that by introducing the parameter 
$s$, the variable $\lambda$ becomes an additional unknown quantity to be determined. As a result, the system becomes undetermined and, to obtain a unique solution, it has to be coupled with the additional constraint 
\begin{equation}
\mathcal{N}(u,\lambda,s)=\dot{u_0}^T(u(s)-u(s_0))+\dot{{\lambda}}_0(\lambda(s)-\lambda(s_0))-(s-s_0)=0,
\label{sec3.1eq3b}
\end{equation}
where $u(s_0)$ denotes the previously computed solution on the solution branch at $s=s_0$, extending the original system as follows
\begin{equation}
\mathcal{G}(u,\lambda,s) \doteq \begin{pmatrix}
G(u,\lambda) \\
\mathcal{N}(u,\lambda,s)
\end{pmatrix} = 0.
\label{sec3.1eq3}
\end{equation}
Equation \eqref{sec3.1eq3b} says that the new point $(u(s),\lambda(s))$ lies on a hyperplane orthogonal to the tangent vector at the current point $(u(s_0),\lambda(s_0))$, and the intersection of the hyperplane with the tangent vector is at a distance $ds=s-s_0$.

\begin{figure}[htbp]
    \centering
\tdplotsetmaincoords{70}{100}
\begin{tikzpicture}[tdplot_main_coords, scale=2]

  \draw[thick,->] (0,0,0) -- (2.,0,0) node[anchor=north east]{$\boldsymbol{s}$};
  \draw[thick,->] (0,0,0) -- (0,2.5,0) node[anchor=north west]{$\boldsymbol{\lambda}$};
  \draw[thick,->] (0,0,0) -- (0,0,2.) node[anchor=east]{$\boldsymbol{u}$};

  \draw[domain=0:1.25, smooth, variable=\t, thick, blue] 
    plot ({0}, {0.5 + \t}, {0.5 + \t*\t});

  \def\xT{0}
  \def\yT{1.0}
  \def\zT{0.75}
  \def\vscale{0.8} 
  \node at (0,1.45,1.05) {\scriptsize$ds$};
  \node at (0,2.6,1.15) {\footnotesize$\mathcal{N}(u,\lambda,s)$};
  
  \node at (0,2.2,1.6) [xshift=10pt] {\scriptsize\textcolor{red}{$(u_{\text{in}},\lambda_{\text{in}})$}};
  
  \node at (0,1.1,1.7) [xshift=-10pt] {\scriptsize\textcolor{blue}{$(u(s),\lambda(s))$}};
  
  \node[blue] at (0,0.5,0.4) [yshift=-8pt] {\footnotesize$G(u,\lambda)$};
  \draw[thick, green!70!black, ->] (\xT,\yT,\zT) -- ++(0,1.1,1.1);
  \filldraw[black] (0,1,0.75) circle (1pt) node[above left] {\tiny$(u(s_0),\lambda(s_0))$};
  
  \pgfmathsetmacro\yEnd{\yT + \vscale}
  \pgfmathsetmacro\zEnd{\zT + \vscale}

  \filldraw[red] (\xT,\yEnd,\zEnd) circle (1pt) node[above right] {};

  \def\rsize{0.35}
  \draw[fill=gray!40, opacity=0.4]
    (\xT - \rsize, \yEnd - \rsize, \zEnd + \rsize) --
    (\xT + \rsize, \yEnd - \rsize, \zEnd + \rsize) --
    (\xT + \rsize, \yEnd + \rsize, \zEnd - \rsize) --
    (\xT - \rsize, \yEnd + \rsize, \zEnd - \rsize) -- cycle;

  \def\tintersect{1.112}
  \pgfmathsetmacro\xI{0}
  \pgfmathsetmacro\yI{0.5 + \tintersect}
  \pgfmathsetmacro\zI{0.5 + \tintersect*\tintersect}

  \filldraw[blue] (\xI,\yI,\zI) circle (1pt) node[below left] {};

\end{tikzpicture}
\caption{Visualization of the arclength continuation method.}
\label{fig:arc_leng}
\end{figure}
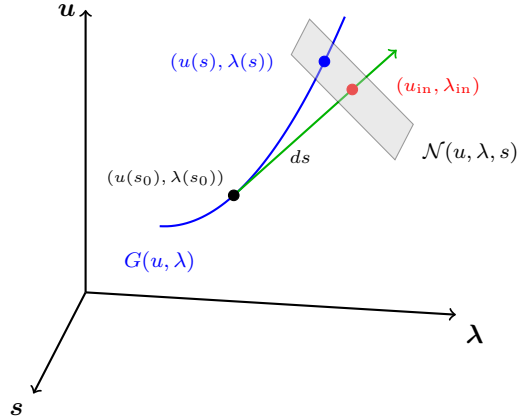

In Figure \ref{fig:arc_leng} we illustrate the arclength continuation method. The black dot indicates the previous solution $(u(s_0), \lambda(s_0))$, at which the tangent is computed. The red dot, $(u_{\text{in}},\lambda_{\text{in}})$, represents the initial guess used to compute the solution $(u(s), \lambda(s))$, shown as the blue dot.  It is evident that using the state corresponding to the red dot as the initial guess leads to faster convergence compared to the original previous solution.

The solution of the extended system in Equation \eqref{sec3.1eq3} with the reparametrization given by $s$ requires the computation of $(\dot{u},\dot{\lambda})$. From Equation \eqref{sec3.1eq1}, $\dot{u}$ can be obtained as
\begin{equation*}
    \dot{u}=-G_u^{-1}G_{\lambda}\dot{\lambda}.
\end{equation*}
Moreover, since $\dot{\lambda}$ is a scalar quantity, it can be obtained from Equation \eqref{sec3.1eq2} as follows
\begin{equation*}
\|G_u^{-1}G_{\lambda}\dot{\lambda}\|^2+|\dot{\lambda}|^2=1 \quad \Rightarrow \quad |\dot{\lambda}|=(1+\|G_u^{-1}G_{\lambda}\|^2)^{-\frac{1}{2}}.
\end{equation*}
The final step is to determine the sign of $\dot{\lambda}$, which is crucial to correctly follow the solution path. If the arclength continuation proceeds forward in the parameter space, the sign should be taken as positive; otherwise, it should be negative.

Starting from a previously computed solution $(u_0,\lambda_0)$ on the branch at the arclength parameter $s=s_0$, the arclength continuation method finds the next parametric sample by increasing the arclength step by $ds=s-s_0$, and solves the extended system in Equation \eqref{sec3.1eq3} by Newton's method with the initial guess (representing the ``corrector" step in the aforementioned analogy). In this way, we enable faster convergence of Newton's method, considering the tangent at point $(u_0,\lambda_0)$, and by providing a better initial guess given by $(u_{\text{in}},\lambda_{\text{in}})=(u_0,\lambda_0)+ds(\dot{u}_0,\dot{\lambda}_0)$. In Algorithm \ref{alg:arclegth} we report a simple implementation of the arclength continuation strategy.

It is important to emphasize that the formulation presented in this subsection is expressed in a continuous setting to facilitate a clear conceptual understanding of the arclength continuation framework. However, in the actual numerical implementation, all quantities, including the solution 
$u(s)$, its derivatives, i.e., \ the associated tangent direction, are obtained from the discretized version of the governing system \eqref{sec1eq2}.


\begin{algorithm}
\caption{Arclength Continuation}\label{alg:arclegth}
\begin{algorithmic}[1]
\Statex{\textbf{Input:} $\text{Initial solution}~(u_i,\lambda_i), \text{nonlinear operator}~ G(u,\lambda),\text{arclength}~ds$}
\Statex{\textbf{Output:} \text{Solution $(u_{i+1},\lambda_{i+1})$}}
\Procedure{Arclength}{$u_i,\lambda_i,G,ds$}
        \State $\text{Compute}~|\dot{\lambda}_i|=(1+\|G_u^{-1}G_{\lambda}\|^2)^{-\frac{1}{2}}$
        \State $\text{Compute}~\dot{u}_i=-G_u^{-1}G_{\lambda}\dot{\lambda}_i$
        \State \text{Set} $(u_{\text{in}},\lambda_{\text{in}})=(u_i,\lambda_i)+ds(\dot{u}_i,\dot{\lambda}_i)$\Comment{Predictor step}
        \State $\text{Solve}~\mathcal{G}(u_{i+1},\lambda_{i+1},s) = 0, ~ \text{with the initial guess} ~(u_{\text{in}},\lambda_{\text{in}})$\Comment{Corrector step}
\EndProcedure
\end{algorithmic}
\end{algorithm}

\subsection{Deflation method}
Understanding and accurately describing all coexisting solution branches that may occur in nonlinear phenomena is a fundamental challenge in bifurcation analysis. One common approach to address this issue is to repeatedly apply Newton's method with a variety of initial guesses, hoping that each guess falls within a different basin of attraction, thereby leading to distinct solutions. Of course, this strategy is not efficient, as it requires many computations to find the different solutions, nor reliable, since there is no guaranty to discover the whole bifurcating scenario. As a potential strategy to overcome this issue, in \cite{farrell2015deflation} the authors proposed a deflation technique that systematically modifies the PDE residual to penalize the convergence to previously identified solutions, thus guiding the solver towards discovering new ones without the need to explore different initial guesses.

Mathematically speaking, assuming that a solution $u_i^1=u^1(\lambda_i)$ to Equation \eqref{sec1eq1} is known (or it has been obtained numerically) for a given value of the parameter $\lambda_i\in \mathcal{P}$, the strategy, reported in Algorithm \ref{alg:def}, deflates the nonlinear system to guide Newton's method towards the convergence to a different solution $u^j(\lambda_i)$ with $j\neq 1$, and the process can be repeated until all admissible solutions for the parameter $\lambda_i$ are obtained, and Newton's method finally diverges.

The procedure works as follows: given $\lambda_i\in \mathcal{P}$ and a corresponding solution,
the original undeflated residual is modified as 
\begin{equation}\label{sec3.2eq1}
    F(u,\lambda_i)=M(u,u_i^1)G(u,\lambda_i),
\end{equation}
where $M(u,u_i^*)$ is the so-called deflation operator defined as
$$M(u,u_i^*)=\left(\frac{1}{\|u-u^*_i\|^p}+\alpha\right)\mathcal{I},$$
where $\mathcal{I}$ is the identity operator, $\alpha$ is the shift parameter, and $p$ is the power parameter. From the deflated system \eqref{sec3.2eq1}, it is clear that when Newton's iterations are converging to the already known solution $u^1_i$, the denominator vanishes $\|u-u^1_i\|\rightarrow 0$, and the solver deviates from its original path, potentially discovering a new solution $u_i^2=u^2(\lambda_i)$. Thus, the deflation process can be performed iteratively composing different deflation operators until the solver diverges for the system 
\begin{equation*}
    F(u,\lambda_i)=M(u,u_i^1)M(u,u_i^2)\ldots M(u,u_i^n)G(u,\lambda_i),
\end{equation*}
possibly having discovered $n$ distinct solutions.

The extra term $\alpha>0$ in the deflation operator $M(u,u_i^*)$, ensures that the norm of the deflated residual does not artificially go to zero as $\|u-u^*_i\|\rightarrow \infty$. The deflation scales the residual by $\alpha$ far away from deflated roots, thus, if $\alpha$ is small then $\|u-u^*_i\|^{-p}$ dominates and Newton's iterates move away from the solution $u_i^*$, while when $\alpha$ is large it is the dominating term and Newton's method remains near the solution. Here, we briefly outline the computation of the deflated solution based on the previously obtained solution starting from the original system \eqref{sec1eq1}. For a more detailed discussion, the reader is referred to \cite{farrell2015deflation}.

Consider the Newton correction step for the discretized form of Equation \eqref{sec1eq1}, which is given by
\begin{equation}\label{sec3.2eq21}
J_G(u)\Delta u_G = -G(u),
\end{equation}
where $\Delta u_G$ denotes the Newton correction corresponding to the original system. It is important to note that, for notational simplicity, we write $G(u)$ instead of $G(u, \lambda_i)$, and $J_G$ represents the Jacobian of $G(u)$ with respect to $u$. Analogously to the original system, the Newton correction step for the deflated system in Equation \eqref{sec3.2eq1} can be expressed as
\begin{equation}\label{sec3.2eq22}
J_F(u)\Delta u_F = -F(u).
\end{equation}
Since the deflated residual is defined as $F(u) = M(u)G(u)$, 
letting $E = M'(u)$, the Newton step becomes
\begin{align*}
\Delta u_F &= -J_F^{-1}F, \\
&= -{(MJ_G + GE^{T})}^{-1}(MG).
\end{align*}
By applying the Sherman–Morrison–Woodbury formula, the correction step simplifies to
\begin{equation*}
\Delta u_F = \tau \Delta u_G,
\end{equation*}
where the scalar $\tau$ is given by
\begin{equation*}
\tau = 1 + \frac{M^{-1}E^{T}\Delta u_G}{1 - M^{-1}E^{T}\Delta u_G}.
\end{equation*}
This shows that the Newton correction step for the deflated system can be efficiently computed by first obtaining the original Newton step $\Delta u_G$ and then evaluating the scalar $\tau$, which only requires a dot product between $E$ and $\Delta u_G$. 

It is important to note that the introduction of the deflation operator modifies the nonlinear residual and, consequently, influences the convergence behavior of the Newton iterations. However, according to the Deflated Rall–Rheinboldt theorem \cite[Theorem 4.5]{farrell2015deflation}, the deflated system still converges to different coexisting solutions. The deflation operator effectively rescales the residual to suppress the attraction toward known solutions, without altering the local quadratic convergence near new roots. In practice, a suitable selection of the shift parameter $\alpha$ and the power parameter $p$ is essential to ensure both numerical stability and robust convergence.

\begin{algorithm}
\caption{Deflation method}\label{alg:def}
\begin{algorithmic}[1]
 \Statex{\textbf{Input:} $\text{Initial guess}~u^1_{i-1}, \text{discovered solutions}~(u^1_i,u^2_i,\ldots u^n_i),\text{nonlinear operator}~ G(u,\lambda)$}
 \Statex{\textbf{Output:} \text{Solution $u^{n+1}_i$}}
 \Procedure{deflation}{$u_{i-1}^1,\{u_i^j\}_{j=1}^n,\lambda_i,G$}
        \State $F(u,\lambda_i)=M(u,u_i^1)M(u,u_i^2)\ldots M(u,u_i^n)G(u,\lambda_i)$\Comment{Deflated system}
        \State $\text{Solve } F(u,\lambda_i) \text{ with initial guess } u^1_{i-1}$
        \If{convergence}
            \State $u^{n+1}_i = u$ \Comment{Deflated solution}
        \Else
            \State \textbf{break}
        \EndIf
\EndProcedure
\end{algorithmic}
\end{algorithm}

\section{A comprehensive strategy to reconstruct multiparameter bifurcation diagrams}\label{sec3a}
In this section, we present the proposed strategy for capturing bifurcating phenomena in multiparametric PDEs, combining the arclength continuation with the deflation technique. We start by introducing the deflated arclength continuation method to efficiently construct the bifurcation diagram in the single-parameter setting. This approach is then extended to handle multiparameteric scenarios by enforcing the arclength continuation to follow pre-defined paths. Then, we build on this to propose an efficient zig-zag path following strategy to detect and reconstruct critical curves and surfaces giving rise to the bifurcating phenomena in the multiparametric setting.

\subsection{Deflated arclength continuation algorithm}\label{sec3.3}

Let us start by denoting with $\mathcal{P}_h$ a discretization of the parametric space $\mathcal{P} \subset \mathbb{R}$. For the first parameter $\lambda_1 \in \mathcal{P}_h$, if no prior information about the system is available, we solve the nonlinear system in Equation \eqref{sec1eq1} using a trivial initial guess $u^1_0$. Assuming the convergence of Newton's method to a solution $u^1_1$, the original system is deflated via Algorithm \ref{alg:def} to seek coexisting solutions for the same value of the parameter $\lambda_1 \in \mathcal{P}_h$. When deflation fails to produce new solutions (i.e., the solver diverges), we utilize the arclength continuation algorithm with the initial solution $(u^1_1,\lambda_1)$ to obtain the next pair $(u^1_2,\lambda_2)$ and repeat the deflation process.

Assume that at $(i-1)$-th step of the arclentgh continuation, we obtain $u_i^1$ for $\lambda_i$, and deflation subsequently yields additional solutions $u_i^j$ for $j = 2, 3, \ldots, n$ (e.g.\ $n=3$ for a pitchfork bifurcation, and $n=2$ for a saddle-node bifurcation). 

At this point, the idea is to reconstruct the $j$-th branch, $\mathcal{U}^j$, using arclength continuation, and proceed with the remaining branches by deflating the system for all newly discovered solutions. Specifically, let us suppose that, at the $i$-th step of the arclength strategy, we have obtained the solution $u_{i+1}^j$ on branch $j$ and the corresponding parameter value $\lambda_{i+1}$. We then deflate $u_{i+1}^j$ and use $u_i^k$ as the initial guess to obtain $u_{i+1}^k$ on some branch $k \neq j$. This step is expected to converge rapidly since the initial guess lies on the same branch. Next, we deflate both $u_{i+1}^j$ and $u_{i+1}^k$, and use $u_{i}^r$ as the initial guess to obtain $u_{i+1}^r$ on some other branch $r\neq j,k$ and so on until all solutions corresponding to $\lambda_{i+1}$ on each branch are recovered.

Concisely speaking, marching along the $j$-th branch using arclength continuation, we generate the sequence ${u_i^j}$ for $i = 1, 2, \ldots$ and the associated parameter values ${\lambda_i}$. By deflating each $u_i^j$, we obtain the remaining solutions $u_i^k$ on each corresponding branch $k$.

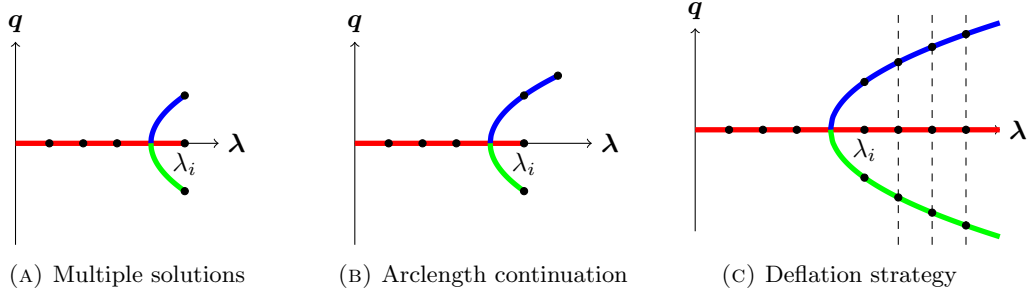
\begin{figure}[htbp]
\centering

\begin{subfigure}[b]{0.28\textwidth}
\centering
\begin{tikzpicture}[scale=0.9]
  \draw[->] (-2,0) -- (1,0) node[right] {$\boldsymbol{\lambda}$};
  \draw[->] (-2,-1.5) -- (-2,1.5) node[above] {$\boldsymbol{q}$};

  \draw[line width=2pt, red, domain=-2:0.5, samples=100, smooth] plot (\x, {0});
  \draw[line width=2pt, blue, domain=0:0.5, samples=100, smooth] plot (\x, {sqrt(\x)});
  \draw[line width=2pt, green, domain=0:0.5, samples=100, smooth] plot (\x, {-sqrt(\x)});

  \node[below] at (0.5,0) {$\lambda_i$};

  \foreach \x in {-1.5, -1, -0.5, 0.5} {
    \filldraw[black] (\x,0) circle (1.5pt);
  }
  \foreach \x in {0.5} {
    \filldraw[black] (\x,{sqrt(\x)}) circle (1.5pt);
    \filldraw[black] (\x,{-sqrt(\x)}) circle (1.5pt);
  }
\end{tikzpicture}
\caption{Multiple solutions}
\label{fig:pitchfork1}
\end{subfigure}
\quad
\begin{subfigure}[b]{0.28\textwidth}
\centering
\begin{tikzpicture}[scale=0.9]
  \draw[->] (-2,0) -- (1.5,0) node[right] {$\boldsymbol{\lambda}$};
  \draw[->] (-2,-1.5) -- (-2,1.5) node[above] {$\boldsymbol{q}$};

  \draw[line width=2pt, red, domain=-2:0.5, samples=100, smooth] plot (\x, {0});
  \draw[line width=2pt, blue, domain=0:1, samples=100, smooth] plot (\x, {sqrt(\x)});
  \draw[line width=2pt, green, domain=0:0.5, samples=100, smooth] plot (\x, {-sqrt(\x)});

  \node[below] at (0.5,0) {$\lambda_i$};

  \foreach \x in {-1.5, -1, -0.5, 0.5} {
    \filldraw[black] (\x,0) circle (1.5pt);
  }
  \foreach \x in {0.5, 1} {
    \filldraw[black] (\x,{sqrt(\x)}) circle (1.5pt);
    }
     \foreach \x in {0.5} {
    \filldraw[black] (\x,{-sqrt(\x)}) circle (1.5pt);
  }
\end{tikzpicture}
\caption{Arclength continuation}
\label{fig:pitchfork2}
\end{subfigure}
\quad
\begin{subfigure}[b]{0.28\textwidth}
\centering
\begin{tikzpicture}[scale=0.9]
  \draw[->] (-2,0) -- (2.5,0) node[right] {$\boldsymbol{\lambda}$};
  \draw[->] (-2,-1.5) -- (-2,1.5) node[above] {$\boldsymbol{q}$};

  \draw[line width=2pt, red, domain=-2:2.5, samples=100, smooth] plot (\x, {0});
  \draw[line width=2pt, blue, domain=0:2.5, samples=100, smooth] plot (\x, {sqrt(\x)});
  \draw[line width=2pt, green, domain=0:2.5, samples=100, smooth] plot (\x, {-sqrt(\x)});

  \node[below] at (0.5,0) {$\lambda_i$};
  \draw[dashed] (1,-1.7) -- (1,1.7);
  \draw[dashed] (1.5,-1.7) -- (1.5,1.7);
  \draw[dashed] (2,-1.7) -- (2,1.7);

  \foreach \x in {-1.5, -1, -0.5, 0.5, 1, 1.5, 2} {
    \filldraw[black] (\x,0) circle (1.5pt);
  }
  \foreach \x in {0.5, 1, 1.5, 2} {
    \filldraw[black] (\x,{sqrt(\x)}) circle (1.5pt);
    \filldraw[black] (\x,{-sqrt(\x)}) circle (1.5pt);
  }
\end{tikzpicture}
\caption{Deflation strategy}
\label{fig:pitchfork3}
\end{subfigure}
\caption{Sketch of pitchfork bifurcation discovery.}
\label{fig:pitchfork_comparison}
\end{figure}

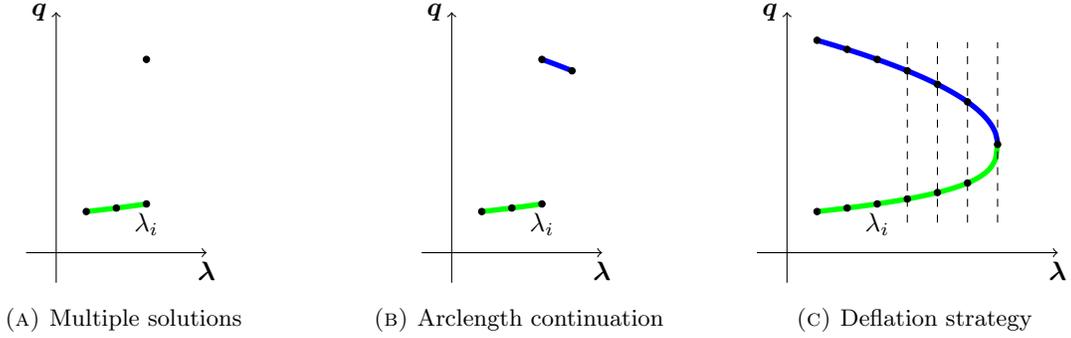
\begin{figure}[htbp]
\centering

\begin{subfigure}[b]{0.31\textwidth}
\centering
\begin{tikzpicture}[scale=0.8] 
    \draw[->] (-0.5,0) -- (2.5,0) node[below] {$\boldsymbol{\lambda}$};
    \draw[->] (0,-0.5) -- (0,4) node[left] {$\boldsymbol{q}$};
    \draw[line width=2pt, green, domain=0.5:1.5, samples=100, smooth]
        plot (\x, {-pow(3.5 - \x, 0.25) + 2});
    \foreach \x in {0.5, 1, 1.5} {
        \filldraw[black] (\x, {-pow(3.5 - \x, 0.25) + 2}) circle (1.5pt);
    }
    \foreach \x in {1.5} {
        \filldraw[black] (\x, {sqrt(3.5 - \x)+1.8}) circle (1.5pt);
    }
    \node[below] at (1.5,0.8108) {$\lambda_i$};
\end{tikzpicture}
\caption{Multiple solutions}
\label{fig:saddle1}
\end{subfigure}
\hfill
\begin{subfigure}[b]{0.31\textwidth}
\centering
\begin{tikzpicture}[scale=0.8]
    \draw[->] (-0.5,0) -- (2.5,0) node[below] {$\boldsymbol{\lambda}$};
    \draw[->] (0,-0.5) -- (0,4) node[left] {$\boldsymbol{q}$};
    \draw[line width=2pt, green, domain=0.5:1.5, samples=100, smooth]
        plot (\x, {-pow(3.5 - \x, 0.25) + 2});
    \draw[line width=2pt, blue, domain=1.5:2.0, samples=100, smooth]
        plot (\x, {sqrt(3.5 - \x)+1.8});
        \foreach \x in {0.5, 1, 1.5} {
        \filldraw[black] (\x, {-pow(3.5 - \x, 0.25) + 2}) circle (1.5pt);
    }
    \foreach \x in {1.5, 2.0} {
        \filldraw[black] (\x, {sqrt(3.5 - \x)+1.8}) circle (1.5pt);
    }
    \node[below] at (1.5,0.8108) {$\lambda_i$};
\end{tikzpicture}
\caption{Arclength continuation}
\label{fig:saddle2}
\end{subfigure}
\hfill
\begin{subfigure}[b]{0.31\textwidth}
\centering
\begin{tikzpicture}[scale=0.8]
    \draw[->] (-0.5,0) -- (4.5,0) node[below] {$\boldsymbol{\lambda}$};
    \draw[->] (0,-0.5) -- (0,4) node[left] {$\boldsymbol{q}$};
    \draw[line width=2pt, green, domain=0.5:3.5, samples=100, smooth]
        plot (\x, {-pow(3.5 - \x, 0.25) + 2});
    \draw[line width=2pt, blue, domain=0.5:3.5, samples=100, smooth]
        plot (\x, {sqrt(3.5 - \x)+1.8});
        \foreach \x in {0.5, 1, 1.5,2.0,2.5,3.0} {
        \filldraw[black] (\x, {-pow(3.5 - \x, 0.25) + 2}) circle (1.5pt);
    }
    \foreach \x in {0.5,1.0,1.5,2.0,2.5,3.0,3.5} {
        \filldraw[black] (\x, {sqrt(3.5 - \x)+1.8}) circle (1.5pt);
    }
    \node[below] at (1.5,0.8108) {$\lambda_i$};
  \draw[dashed] (2.0,0.5) -- (2.0,3.5);
  \draw[dashed] (2.5,0.5) -- (2.5,3.5);
  \draw[dashed] (3.0,0.5) -- (3.0,3.5);
  \draw[dashed] (3.5,0.5) -- (3.5,3.5);
\end{tikzpicture}
\caption{Deflation strategy}
\label{fig:saddle3}
\end{subfigure}

\caption{Sketch of saddle-node bifurcation discovery.}
\label{fig:saddle_comparison}
\end{figure}

For clarity, we report the main steps of the deflated arclength continuation strategy in Algorithm \ref{alg:defl-arc}, returning $n$ qualitatively different coexisting branches. 
Moreover, schematic illustrations of the proposed deflated arclength continuation algorithm for both pitchfork and saddle-node bifurcations are provided in Figures~\ref{fig:pitchfork_comparison} and~\ref{fig:saddle_comparison}, respectively.

In the pitchfork bifurcation case, multiple solutions are first obtained near the bifurcation point $\lambda_i$ (Figure~\ref{fig:pitchfork1}). The arclength continuation method is then applied to further trace the upper branch (Figure~\ref{fig:pitchfork2}). Meanwhile, deflation is used at each arclength continuation step to identify the solution on some branch using the previously computed solutions as initial guess at that corresponding branch, thereby recovering the complete bifurcation structure (Figure~\ref{fig:pitchfork3}).

Similarly, in the saddle-node bifurcation case, multiple solutions are first computed at a parameter value $\lambda_i$ (Figure~\ref{fig:saddle1}). The arclength continuation is then used to extend the upper branch (Figure~\ref{fig:saddle2}), while deflation is applied simultaneously to uncover the solutions on the lower branch, ultimately reconstructing the entire bifurcation diagram (Figure~\ref{fig:saddle3}). An alternative approach to constructing the saddle-node bifurcation is to proceed backward in the parameter space. By identifying two solution points near the fold, one can then continue along both branches using arclength continuation.

\begin{algorithm}
\caption{Deflated arclength continuation algorithm}
\label{alg:defl-arc}
\begin{algorithmic}[1]
\Statex \textbf{Input:} 
Initial guess $u_0^1$, arclength step $ds$, maximum value $\lambda_{\max}$,
nonlinear operator $G(u,\lambda)$
\Statex \textbf{Output:} 
Solution branches $\{\mathcal U^j\}_{j=1}^n$
\State $u_1^1 \gets$ Solve Equation \eqref{sec1eq2} with initial guess $u_0^1$ for $\lambda_1 = \lambda_{\min}$
\State Set $i \gets 2$, and $n \gets 1$

\While{$\lambda_{i-1} < \lambda_{\max}$ \text{and} $n = 1$}
    \State $(u_{i}^1,\lambda_{i}) \gets \text{Arclength}(u_{i-1}^1,\lambda_{i-1},G,ds)$ \Comment{Arclength step on a branch}
    \State $\{u_{i}^j\}_{j=2}^n \gets $\text{Deflation}$\big(u_i^1, \{u_{i}^k\}_{k<j},\lambda_i,G\big)$ \Comment{Discover new branches via deflation}
        \State $i \gets i+1$
    \EndWhile
    \While{$\lambda_i < \lambda_{\max}$}
    \State $(u_{i+1}^j,\lambda_{i+1}) \gets 
    \text{Arclength}(u_{i}^j,\lambda_{i},G,ds)$ 
    \Comment{Arclength on $j$-th branch}
    \State $\{u_{i+1}^r\}_{r\neq j}\gets \text{Deflation}\big(u_i^r, \{\{u_{i+1}^k\}_{k<r},u_{i+1}^j\},\lambda_{i+1},G\big)$\Comment{Solutions on known branches}
    \State $\{u_{i+1}^k\}_{k> n}\gets\text{Deflation}\big(u_i^1, \{u_{i+1}^j\}_{j\leq k},\lambda_{i+1},G\big)$ \Comment{Discover new branches}
    \State $i \gets i+1$

\EndWhile
\end{algorithmic}
\end{algorithm}

\subsection{Multiparametric bifurcation reconstruction}
In this section, we begin by extending the classical arclength continuation method to accommodate multiparametric settings, enabling the tracking of solution paths in a higher-dimensional parameter space. Building on this foundation, we introduce the multiparametric deflated arclength continuation algorithm, which integrates the deflation technique with the extended arclength continuation framework to robustly and efficiently reconstruct bifurcation diagrams involving multiple parameters. Finally, we present a novel methodology for detecting bifurcation curves in the multiparametric setting. This is achieved through the development of a zigzag path-following strategy, specifically designed to accurately trace bifurcation curves and surfaces, including both saddle-node and pitchfork bifurcations.

\subsubsection{Multiparametric arclength continuation}\label{sec3.4.1}
Here, we extend the arclength continuation technique to a multiparametric setting, which is essential for developing the methodology used to detect bifurcation curves and to construct the bifurcation diagrams in such contexts. To simplify the presentation, we first consider two parameters $\boldsymbol{\lambda}=(\lambda_1, \lambda_2) \in \mathbb{R}^2$ and then generalize the approach to the case of multiple parameters, i.e., $\boldsymbol{\lambda} \in \mathbb{R}^p$ with $p\geq 2$.\\

\noindent\textbf{Case $p=2$}: To reconstruct solution branches in the two parameter context, we let the arclength continuation to follow a smooth curve $g(\lambda_1,\lambda_2)$ in $(\lambda_1,\lambda_2)$-plane.

As before, instead of parameterizing the solution $u$ by $\lambda_1$ and $\lambda_2$, we parameterize the solution branch using the arclength parameter $s$, and by differentiating the system in Equation \eqref{sec1eq1} we obtain
\begin{equation}\label{sec34.eq1}
    \frac{dG(u(s),\lambda_1(s),\lambda_2(s))}{ds}=G_u \dot{u}+G_{\lambda_1} \dot{\lambda}_1+G_{\lambda_2} \dot{\lambda}_2=0.
\end{equation}
Denoting by $x(s)=(u(s),\lambda_1(s),\lambda_2(s))$ the solution vector, $\dot{x}$ represents the derivative with respect to the arclength parameter given by the unit tangent vector with
\begin{equation}\label{sec34.eq2}
    \|\dot{x}\|^2=\|\dot{u}\|^2+|\dot{\lambda}_1|^2+|\dot{\lambda}_2|^2=1.
\end{equation}
Introducing the arclength as parameter renders $(\lambda_1,\lambda_2)$ as unknowns, which requires two additional equations to obtain a consistent system. To do this, Equation \eqref{sec1eq1} can be extended as 
\begin{equation}
{\mathcal{G}}(u,\lambda_1,\lambda_2,s)=\begin{pmatrix}
G(u,\lambda_1,\lambda_2) \\
g(\lambda_1,\lambda_2)\\
\mathcal{N}(u,\lambda_1,\lambda_2,s)
\end{pmatrix}=0,\qquad\text{where}
\label{sec34.eq3}
\end{equation}
$$\mathcal{N}(u,\lambda,\mu,s)=\dot{u_0}^T(u(s)-u(s_0))+(\dot{{\lambda}}_1)_0(\lambda_1(s)-\lambda_1(s_0))+(\dot{{\lambda}}_2)_0(\lambda_2(s)-\lambda_2(s_0))-(s-s_0).$$ 
To compute the tangent vector $(\dot{u},\dot{\lambda}_1,\dot{\lambda}_2)$, we can start by differentiating $g(\lambda_1,\lambda_2)=0$ 
\begin{align}
    g_{\lambda_1}\dot{\lambda}_1+g_{\lambda_2}\dot{\lambda}_2=0 \qquad
    \Rightarrow \qquad \dot{\lambda}_2=-\frac{g_{\lambda_1}\dot{\lambda}_1}{g_{\lambda_2}}, \label{sec34.eq4}
\end{align}
and by exploiting Equations \eqref{sec34.eq1} and \eqref{sec34.eq2}, it can be derived as follows
\begin{equation}
\begin{split}
    \dot{u}&=-G_u^{-1}\left(G_{\lambda_1}-\frac{g_{\lambda_1}}{g_{\lambda_2}}G_{\lambda_2}\right)\dot{\lambda}_1,\\
    |\dot{\lambda}_1|&=\left(1+\left\lVert G_u^{-1}\left(G_{\lambda_1}-\frac{g_{\lambda_1}}{g_{\lambda_2}}G_{\lambda_2}\right)\right\rVert^2+\left|\frac{g_{\lambda_1}}{g_{\lambda_2}}\right|\right)^{-\frac{1}{2}}.  
\end{split}\label{sec34.eq55}
\end{equation}

\noindent\textbf{Case $p\geq 2$}: We now extend the arclength continuation method for general multiparametric case, i.e.\ for $\boldsymbol{\lambda}=(\lambda_1,\lambda_2,\ldots,\lambda_p)\in \mathbb{R}^p,~p>1$. Suppose, $\mathcal{C}$ is a smooth curve in the $p$-dimensional hyperplane, then it can be defined as the following set
\begin{equation*}
    \mathcal{C}=\{\boldsymbol{\lambda} \in \mathbb{R}^p\ |\ g_i(\boldsymbol{\lambda})=0,\, i=1,2,\ldots p-1\},
\end{equation*}
where $g_i(\boldsymbol{\lambda}):\mathbb{R}^p\rightarrow \mathbb{R}$ is a smooth hypersurface in $p$-dimensional space. Parameterize the solution $x=(u,\boldsymbol{\lambda})$ by the arclength parameter $s$ and differentiate the parametric system in Equation \eqref{sec1eq1} to obtain
\begin{equation}\label{sec34.eq5}
    \frac{dG(u(s),\boldsymbol{\lambda}(s))}{ds}=G_u \dot{u}+\sum_{i=1}^p G_{\lambda_i} \dot{\lambda_i}=0.
\end{equation}
Furthermore, suppose that $x(s)=(u(s),\boldsymbol{\lambda}(s))$ is the solution vector then $\dot{x}$ represents the unit tangent vector which implies that
\begin{equation}\label{sec34.eq6}
    \|\dot{x}\|^2=\|\dot{u}\|^2+\sum_{i=1}^p \dot{\lambda_i}=1.
\end{equation}
Since, $\boldsymbol{\lambda}=(\lambda_1,\lambda_2,\ldots,\lambda_p)$ will be the $p$ unknowns in the arclength continuation to determine. Hence, to have a consistent system, the parametric system $G(u,\boldsymbol{\lambda})=0$ can be extended as 
\begin{equation}
\mathcal{G}(u,\boldsymbol{\lambda},s)=\begin{pmatrix}
G(u,\boldsymbol{\lambda}) \\
g_1(\boldsymbol{\lambda})\\
\vdots\\
g_{p-1}(\boldsymbol{\lambda})\\
\mathcal{N}(u,\boldsymbol{\lambda},s)
\end{pmatrix}=0,
\label{sec34.eq7}
\end{equation}
where $$\mathcal{N}(u,\boldsymbol{\lambda},s)=\dot{u_0}^T(u(s)-u(s_0))+\sum_{i=1}^p{(\dot{{\lambda_i}})}_0(\lambda_i(s)-\lambda_i(s_0))-(s-s_0).$$ To compute the tangent vector $(\dot{u},\dot{\boldsymbol{\lambda}})$, the hypersurfaces $g_i(\boldsymbol{\lambda})$ can be differentiated as follows
\begin{equation}\label{sec34.eq8}
    \sum_{j=1}^{p}(g_i)_{\lambda_j}\dot{\lambda}_j=0,\quad i=1,2,\ldots,p-1.
\end{equation}
Finally, one can solve Equations \eqref{sec34.eq5}, \eqref{sec34.eq6} and \eqref{sec34.eq8} to obtain the tangent vector $(\dot{u},\dot{\boldsymbol{\lambda}})$, so that the system $\mathcal{G}(u,\boldsymbol{\lambda},s)$ can be solved using the iterative method with the initial guess $(u_{\text{in}},\boldsymbol{\lambda}_{\text{in}})=(u_{0},\boldsymbol{\lambda}_{0})+ds(\dot{u}_{0},\dot{\boldsymbol{\lambda}}_{0})$.

\subsubsection{Multiparametric deflated arclength continuation algorithm} \label{sec3.4.2}
The set of bifurcating solutions of the system in Equation \eqref{sec1eq1} is commonly explored by fixing all but one parameter for which the one-dimensional bifurcation diagram is computed. 
While such strategy is straightforward, it completely overlooks the multiparametric nature of the bifurcation, neglecting the coupled effects of varying all the parameters simultaneously. To address this, we have developed a strategy for computing bifurcation diagrams by extending the deflated arclength continuation technique introduced in Section \ref{sec3.3} in the multiparameteric context, in combination with the arclength continuation method presented in Section \ref{sec3.4.1}.
For ease of exposition, we restrict the discussion in what follows to the case of two parameters, $\boldsymbol{\lambda} = (\lambda_1,\lambda_2)$, but the proposed methodology naturally generalizes to higher-dimensional parameter spaces.

\begin{figure}[tbp]
    \centering
    \begin{subfigure}[b]{0.3\textwidth}
        \includegraphics[width=\textwidth]{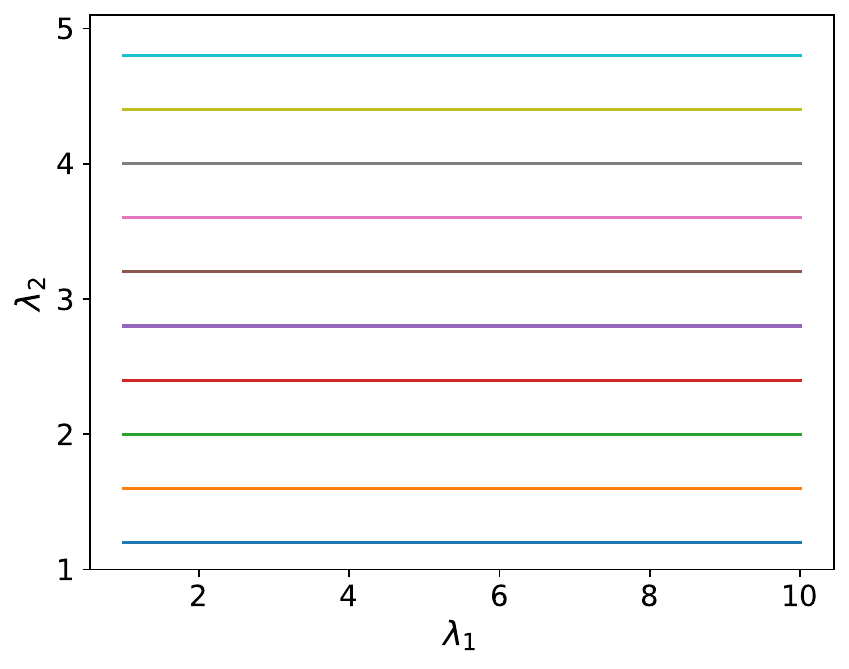}
        \caption{Horizontal lines}
        \label{subfigure211g}
    \end{subfigure}
    \quad
    \begin{subfigure}[b]{0.3\textwidth}
        \includegraphics[width=\textwidth]{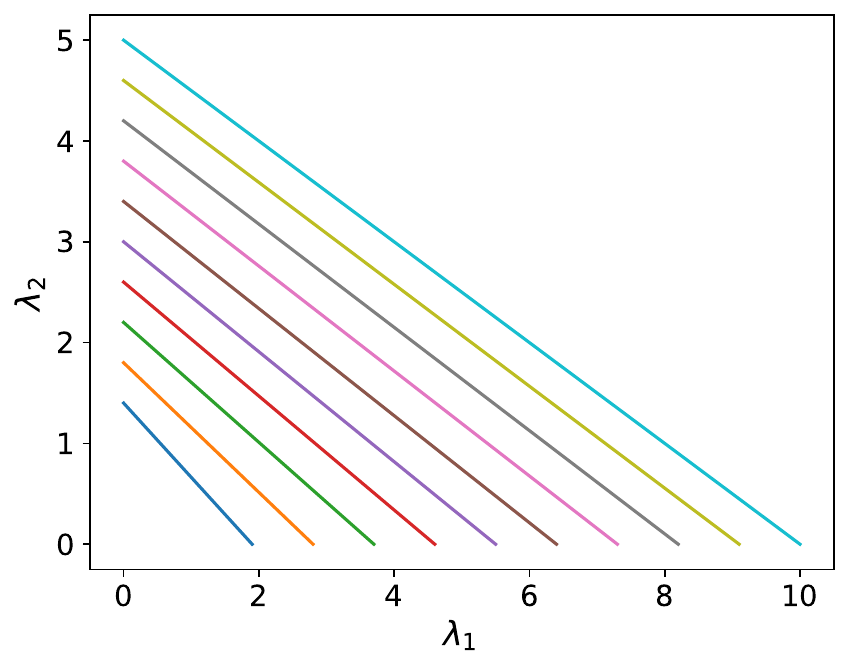}
        \caption{Diagonal lines}
        \label{subfigure212g}
    \end{subfigure}
    \begin{subfigure}[b]{0.3\textwidth}
        \includegraphics[width=\textwidth]{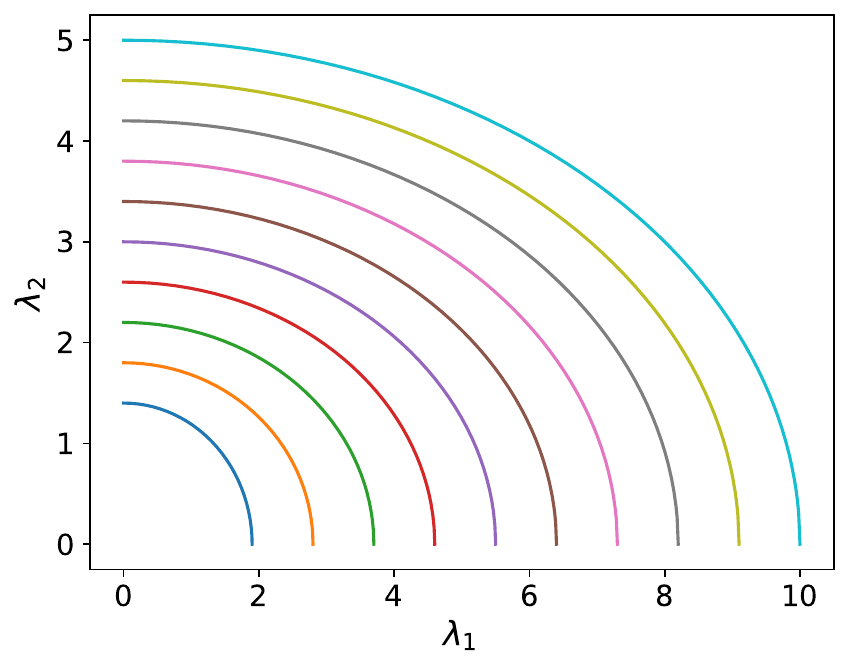}
        \caption{Elliptic curves}
        \label{subfigure213g}
    \end{subfigure}
    \caption{Different paths for the multiparametric arclength continuation.}
    \label{fig:curves}
\end{figure}

The application of multiparametric arclength continuation naturally raises the question of how to select the curve $g(\lambda_1, \lambda_2)$ in a suitable way. To compare different choices, we consider the parameter domain $\mathcal{P} = [a,b]\times [c,d]$ and three representative paths described by the following sets of $n+1$ curves depicted in Figure \ref{fig:curves} for $a=1$, $b=10$, $c=1$, and $d=5$:
\begin{enumerate}
    \item Horizontal: $g(\lambda_1,\lambda_2)=\left\{\lambda_2-ih_1=0\ |\ h_1=\frac{d-c}{n},\ i=0,\ldots,n\right\}$,
    \item Diagonal: $g(\lambda_1,\lambda_2)=\left\{\frac{\lambda_2}{c+ih_1}+\frac{\lambda_1}{a+ih_2}-1=0\ |\ h_1=\frac{d-c}{n},\ h_2=\frac{b-a}{n},\ i=0,\ldots,n\right\}$,
    \item Elliptic: $g(\lambda_1,\lambda_2)=\left\{\frac{\lambda_2^2}{(c+ih_1)^2}+\frac{\lambda_1^2}{(a+ih_2)^2}-1=0\ |\ h_1=\frac{d-c}{n},\ h_2=\frac{b-a}{n},\ i=0,\ldots,n\right\}$.
\end{enumerate}


We aim at using these paths in combination with deflation to construct bifurcation diagrams for multiparametric systems.
In particular, the proposed methodology follows the same conceptual framework as the deflated arclength continuation method discussed in Section \ref{sec3.3}, where the key distinction is the exploitation of the multiparametric arclength continuation along a prescribed path $g(\lambda_1,\lambda_2)$, rather than the standard arclength continuation with respect to a single parameter. In this setting, continuation is performed in the extended $(u,\boldsymbol{\lambda})$-space, where the parameter vector $\boldsymbol{\lambda}$ evolves along the chosen path.
In this way, deflated arclength continuation provides a unified framework for exploring bifurcation structures in multiparametric problems without resorting to parameter grid sampling, while retaining the robustness of arclength continuation and the deflation capability to detect coexisting solutions.


\begin{remark}\label{rem1}
The continuation path should not exhibit sharp directional changes potentially hindering convergence.
Indeed, the arclength continuation method relies on the tangent vector at the previously computed solution to generate an initial guess by intersecting a hyperplane at a fixed arclength $ds$, as illustrated in Figure \ref{fig:arc_leng}. However, near these sharp turns, the initial guess may lie significantly far from the actual solution, often leading to a convergence failure of the solver.
Moreover, using a very small arclength $ds$ does not imply a successful tracking of the path, it may still lead to divergence, convergence to a distant point on the path or skipping over the actual path segment.
\end{remark}



\subsubsection{Bifurcation curve detection}\label{sec3.4.3}
The loss of uniqueness for the solution is linked with the presence of bifurcation points, curves and surfaces, respectively for one-, two- and three-dimensional parameter spaces. 
Even when interested solely in the accurate detection of these bifurcating regions, rather than the complete diagrams, the computational cost is still challenging, especially for multiparametric contexts. In this section, we propose a hierarchical strategy tracing the bifurcation curves in the plane spanned by the first two parameters $(\lambda_1,\lambda_2)$, and then the procedure is repeated for each fixed choice of parameters. 



\begin{figure}[htbp]
\centering
    \includegraphics[width=0.7\linewidth]{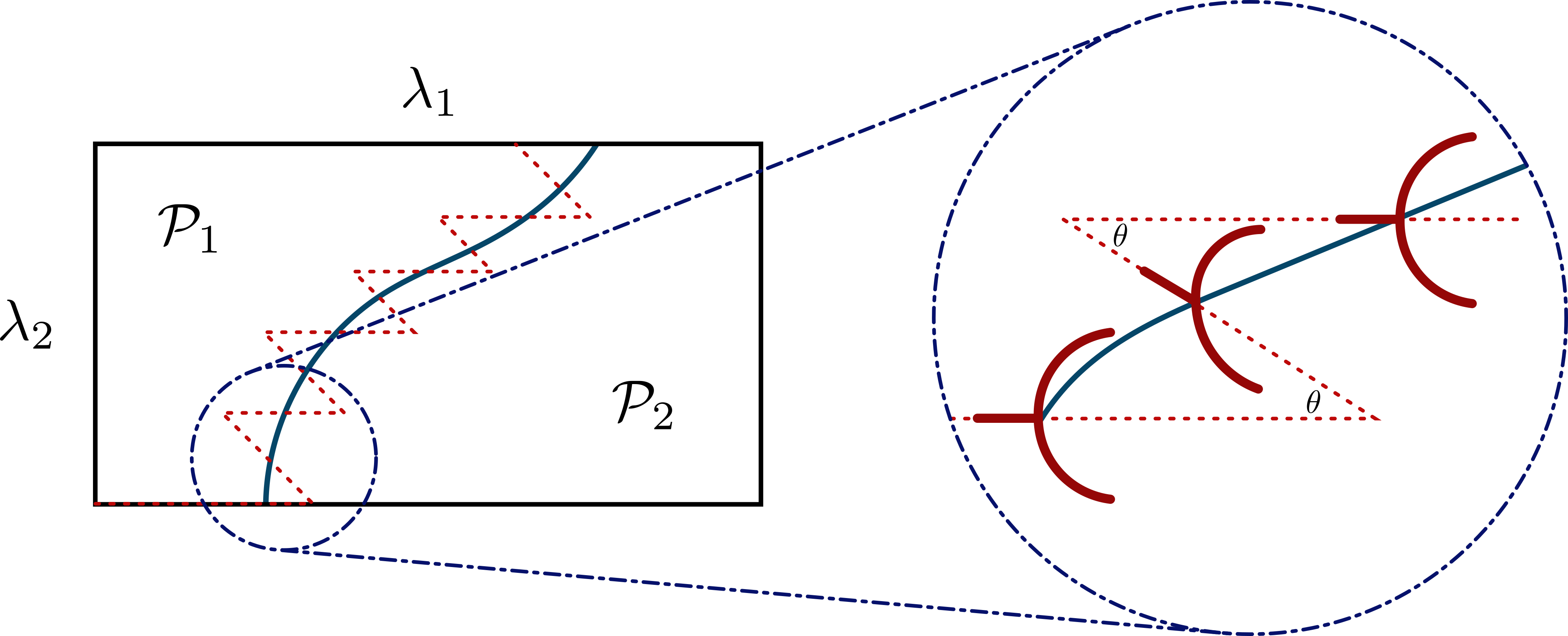}
    \caption{Bifurcation curve and the approximating zigzag tracking path (red) on the $(\lambda_1,\lambda_2)$-domain with $p=2$ with the zoom-in showing the local bifurcation diagram.}
    \label{fig:comparison}
\end{figure}

Thus, when a comprehensive approximation of all solution branches is not needed, we propose a \textit{zigzag tracking path}, illustrated in Figure \ref{fig:comparison} for $p=2$, to shift the focus from the complete reconstruction of the bifurcation diagram to the detection of the bifurcation regions, providing an accurate description of the interface between uniqueness and non-uniqueness.

Let us assume that the parameter ranges are given by $\lambda_1 \in [\lambda_{1, \text{min}}, \lambda_{1, \text{max}}]$ and $\lambda_2 \in [\lambda_{2, \text{min}}, \lambda_{2,\text{max}}]$, with $\mathcal{P}_1$ denoting the region of unique solutions and $\mathcal{P}_2$ the region of multiple solutions. 
Thanks to the methodology developed in the previous section, the idea is to proceed the exploration parallel to one of the two parametric directions and invert the path with a diagonal turn when the bifurcation is detected.
In particular, we start by applying the multiparametric deflated arclength continuation algorithm along the curve $g(\lambda_1, \lambda_2) = \lambda_2 - \lambda_{2, \text{min}}$.
We apply the deflation process until multiple solutions $\{u_i^j\}_j$ are discovered for some $\lambda_{1,i}$. After having computed $k$ additional solutions in the same horizontal direction, where $k$ is usually small but large enough to avoid being too close to the critical region, the continuation direction is updated to follow the zig-zag path defined by
\begin{equation}\label{sec343.eq1}
g(\lambda_1, \lambda_2) = (\lambda_2 - \tilde{\lambda}_2)  - \tan(\theta)( \lambda_1  - \tilde{\lambda}_1),
\end{equation}
where $g(\lambda_1, \lambda_2)$ represents a straight line inclined at an angle $\theta$ with respect to the $\lambda_1$-axis, and the continuation exploits $\{u_{i+k}^j\}_j$ the initial guesses. The quantities $\tilde{\lambda}_1$ and $\tilde{\lambda}_2$ denote the current parameter values, which at the initial stage correspond to $\lambda_{1,i+k}$ and $\lambda_{2,\mathrm{min}}$, respectively.
The angle $\theta$ is a hyperparameter that should be chosen carefully. If a large angle is used, the resulting zigzag path may skip portions of the bifurcation curve, particularly in regions where the curve exhibits sharp turns or rapid changes in direction. A smaller angle leads to a finer sampling of the curve, improving robustness and accuracy in capturing its geometry, but it increases the computational cost. 
For example, when employing the aforementioned strategy for bifurcation curves which contain nearly horizontal segments with respect to the bifurcating parameter $\lambda_1$-axis, such features can still be effectively detected, but a sufficiently small value of $\theta$ needs to be chosen for the continuation path. 

When the algorithm crosses again the bifurcation curve from the non-uniqueness region $\mathcal{P}_2$ to the uniqueness region $\mathcal{P}_1$,  we again continue the computation in the same direction for $k$ successive steps, finally updating again the continuation path via Equation \eqref{sec343.eq1}.
By repeating this alternating strategy for different fixed configurations of the remaining parameters when $p\geq 2$, crossing the bifurcation curve back and forth, the algorithm effectively tracks the bifurcating region with a zigzag path regardless of the dimension of the multiparametric domain.

%% file: sections/results.tex
\section{Numerical Results}
In this section, we discuss the performance of the proposed methodology for constructing bifurcation diagrams in general multiparametric settings. We analyze the developed deflated arclength continuation on three benchmark problems of increasing complexity: starting from the Bratu and Allen--Cahn equations in 1D/2D spatial domains, and finally testing the methodology on the complex Rayleigh--Benard convection problem for both 2D and 3D cavities. Furthermore, we utilize the zigzag path-following strategy to detect bifurcation regions in two- and three-dimensional parametric spaces. 


The nature of the bifurcation differs among these examples. The Bratu equation exhibits a saddle-node bifurcation, while the Allen--Cahn equation undergoes a pitchfork bifurcation. Notably, already in the one-parameter setting, the Allen--Cahn equation is characterized by challenging multiple bifurcation points, where several bifurcating states arise at the same time, with even greater complexity in the two-parameter case. The Rayleigh--Bénard convection problem represents a considerably more challenging large-scale nonlinear PDE system, especially in $3D$ space, characterized by coupled velocity, pressure, and temperature fields and rich bifurcation structures in the Rayleigh--Prandtl parameter space.



For the spatial discretization of the Bratu and Allen--Cahn benchmark problems, we employ the \textit{wavelet collocation method} described in Appendix~\ref{appendix1}. However, we emphasize that the proposed methodologies for bifurcation reconstruction and bifurcation curve detection are independent of the underlying discretization scheme. Indeed, for the Rayleigh--Bénard convection problem in Section \ref{sub_RBC}, we exploit the \text{Finite Element method}. The Bratu and Allen--Cahn computations are carried out in MATLAB using the \textsc{fsolve} routine, whereas the Rayleigh--Bénard simulations are performed in Python using FEniCSx \cite{baratta2023dolfinx}. This demonstrates the applicability and versatility of the proposed strategy across fundamentally different numerical discretization techniques and software environments.

\subsection{Bifurcating benchmarks} 
Here, we present the two models under investigation, the Bratu and the Allen--Cahn equations, in their most general multiparametric form, which will be numerically investigated for 1D and 2D domains via the multiparametric deflated arclength continuation algorithm and bifurcation curve detection.

\subsubsection{Bratu equation}  
 The Bratu equation \cite{boyd1986analytical,fabiani2021numerical,shahab2024neural} is a nonlinear parabolic PDE that exhibits a saddle-node bifurcation with respect to a positive parameter $\lambda_1$. In this work, we consider the following \emph{multiparametric} formulation of the Bratu equation:
\begin{equation}\label{bratu_multi}
\begin{cases}
    \lambda_2 \Delta u + \lambda_1 e^u = 0, & x \in \Omega, \\
    u = \lambda_3, & x \in \partial \Omega,
\end{cases}
\end{equation}
where $x = (x_1,\ldots,x_d)$, $\Omega = [0,1]^d \subset \mathbb{R}^d$, $d$ denotes the spatial dimension, and $\boldsymbol{\lambda}=(\lambda_1,\lambda_2,\lambda_3)$ is the multiparameter, with $\lambda_1$ the component responsible for the primary bifurcation measuring the magnitude of the reaction term, while $\lambda_2$ and $\lambda_3$ represent, respectively, a diffusion coefficient and a non-homogeneous Dirichlet boundary condition. In particular,  we have chosen the parametric ranges as $\mathcal{P} = [0,4]\times[0,10]\times[0,1.5]$ for $d=1$, and $\mathcal{P} = [0,7]\times[0,6]\times[0,1.5]$ for $d=2$, but we note that the effective intervals depend on the chosen continuation path, with $\lambda_1$ serving as the primary bifurcation parameter.
%
In the classical Bratu problem, for fixed values of $\lambda_2$ and $\lambda_3$, the equation exhibits a saddle-node bifurcation at a critical value $\lambda_1 = \lambda_1^*$, beyond which no solution exists.
For $0 < \lambda_1 < \lambda_1^*$, two distinct solutions coexist, commonly referred to as the lower and upper solution branches, while at $\lambda_1 = \lambda_1^*$, these two branches coalesce into a single solution. For the classical case with $p=1$, i.e.\ $\lambda_2 = 1$ and $\lambda_3 = 0$, the approximate critical bifurcation points are $\lambda_1^* \approx 3.5$ for $d=1$, and $\lambda_1^* \approx 6.8$ for $d=2$.

\subsubsection{Allen--Cahn equation} 
We consider the Allen--Cahn equation \cite{kuehn2024uncertainty,kuehn2015numerical,allen1979microscopic}, a classical nonlinear reaction--diffusion model that exhibits multiple pitchfork bifurcations and widely used in materials science to describe phase separation processes, interface motion, and pattern formation in binary alloys. 

In this work, we consider the following \emph{multiparametric} formulation of the Allen--Cahn equation:
\begin{equation}\label{allenCahn_mult}
\begin{cases}
    \lambda_2 \Delta u - u(u^2 - \lambda_1) = 0, & x \in \Omega_{\lambda_3},\\
    u = 0, & x \in \partial \Omega_{\lambda_3},
\end{cases}
\end{equation}
where $x = (x_1,\ldots,x_d)$, $\Omega_{\lambda_3} = [0,\lambda_3]^d \subset \mathbb{R}^d$, and $\boldsymbol{\lambda}=(\lambda_1,\lambda_2,\lambda_3)$ is the multiparameter, with $\lambda_1$ the bifurcating component for the nonlinear term, while $\lambda_2$ and $\lambda_3$ represent, respectively, a diffusion coefficient and a geometric scaling parameter. Here, we have chosen as parametric ranges $\mathcal{P}=[0,14]\times[1,10]\times[\pi,3.8]$ for $d=1$, and $\mathcal{P}=[0,12]\times[1,8]\times[\pi,3.8]$ for $d=2$.


For fixed values of $\lambda_2$ and $\lambda_3$, a stable trivial branch of solution loses its stability and originates two additional nontrivial branches. Moreover, the equation also exhibits a sequence of pitchfork bifurcations as $\lambda_1$ increases further, giving rise to a complex pattern of coexisting solutions and multiple bifurcations.

Unlike the Bratu equation, where the third parameter is introduced through the boundary condition, imposing a nonhomogeneous boundary parameter in the Allen--Cahn equation would destroy the intrinsic symmetry of the problem and thereby prevent the occurrence of a pitchfork bifurcation. For this reason, we consider a geometric parameter $\lambda_3$ scaling the computational domain and preserving the symmetry required for pitchfork-type bifurcations.

\subsection{Reconstructing bifurcation diagrams via deflated arclength continuation}
In this section, we will present the behavior of the solution branches for both benchmarks in 1D/2D settings with an increasing dimensionality of the parametric space.
 
\subsubsection{Bratu equation with $p=1$} 
\begin{figure}[tbp]
    \centering
    \includegraphics[width=0.4\textwidth]{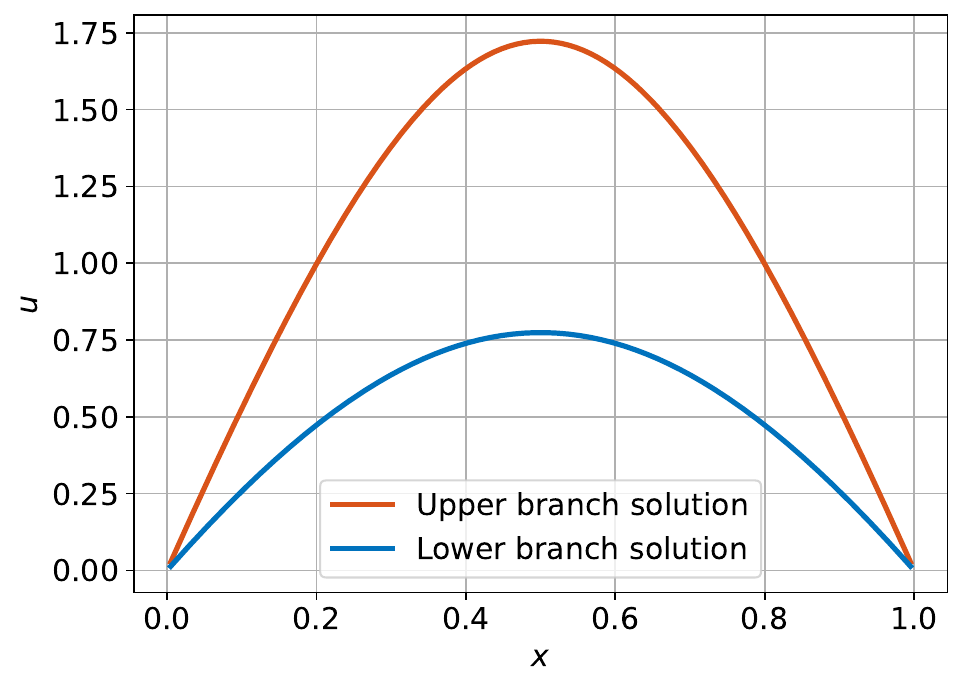} 
    \qquad    
    \includegraphics[width=0.4\textwidth]{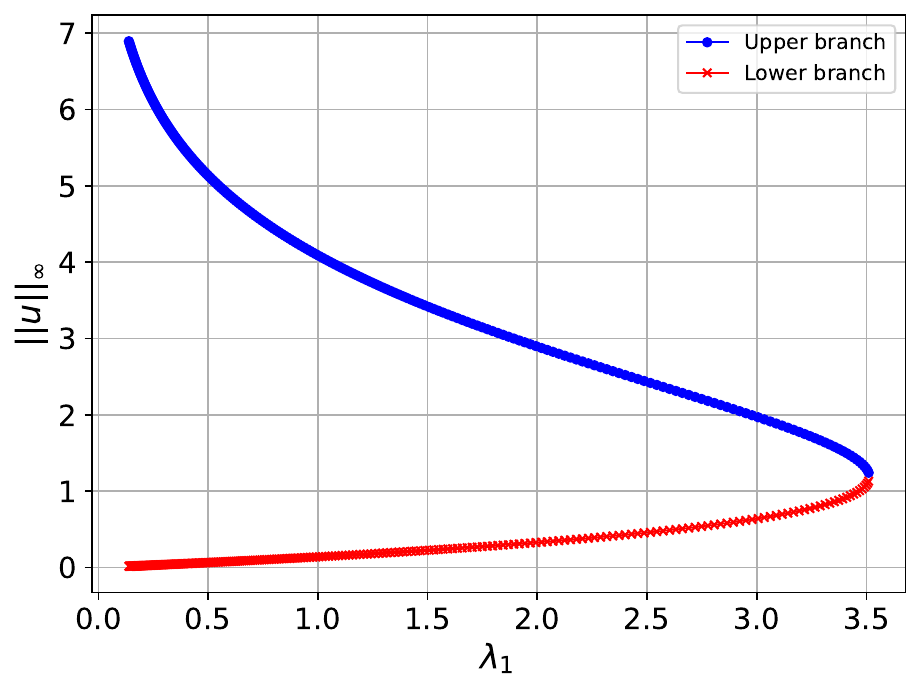}
    \caption{Bratu solutions with $d=1$ for $\lambda_1=3.24$, and saddle-node bifurcation diagram.}
    \label{fig:problem11}
\end{figure}

\begin{figure}[tbp]
    \centering
        \includegraphics[width=0.3\textwidth]{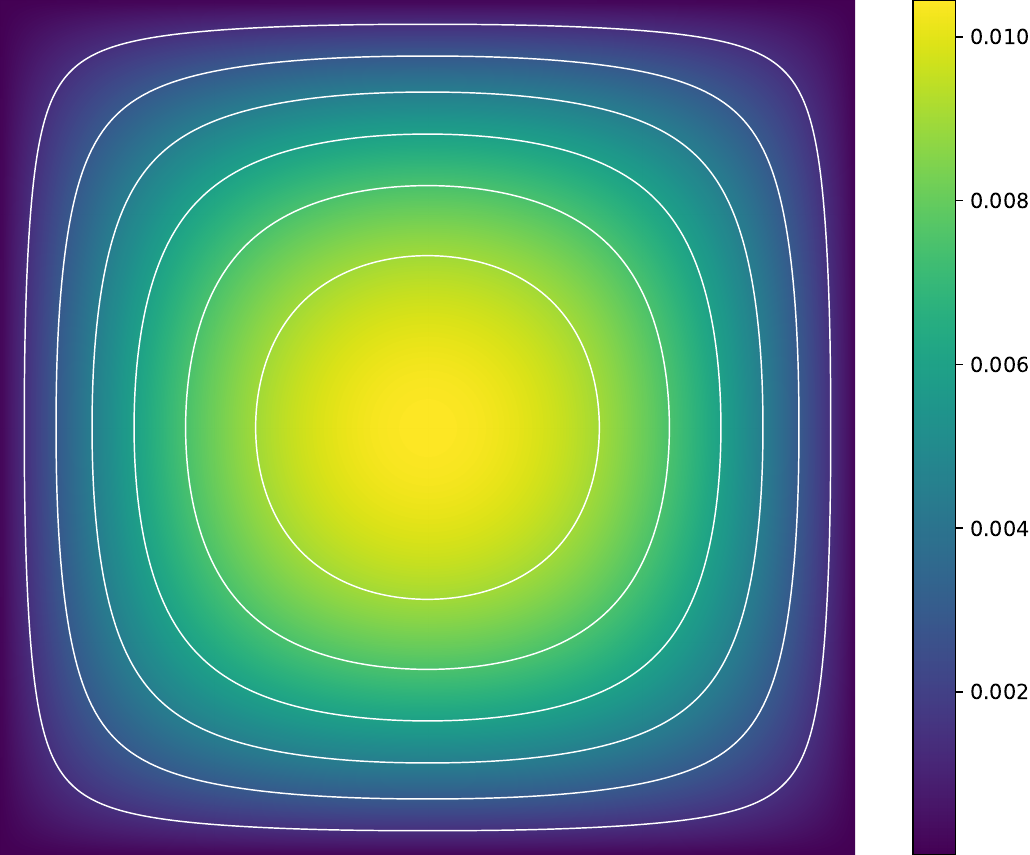}
    \hfill
        \includegraphics[width=0.29\textwidth]{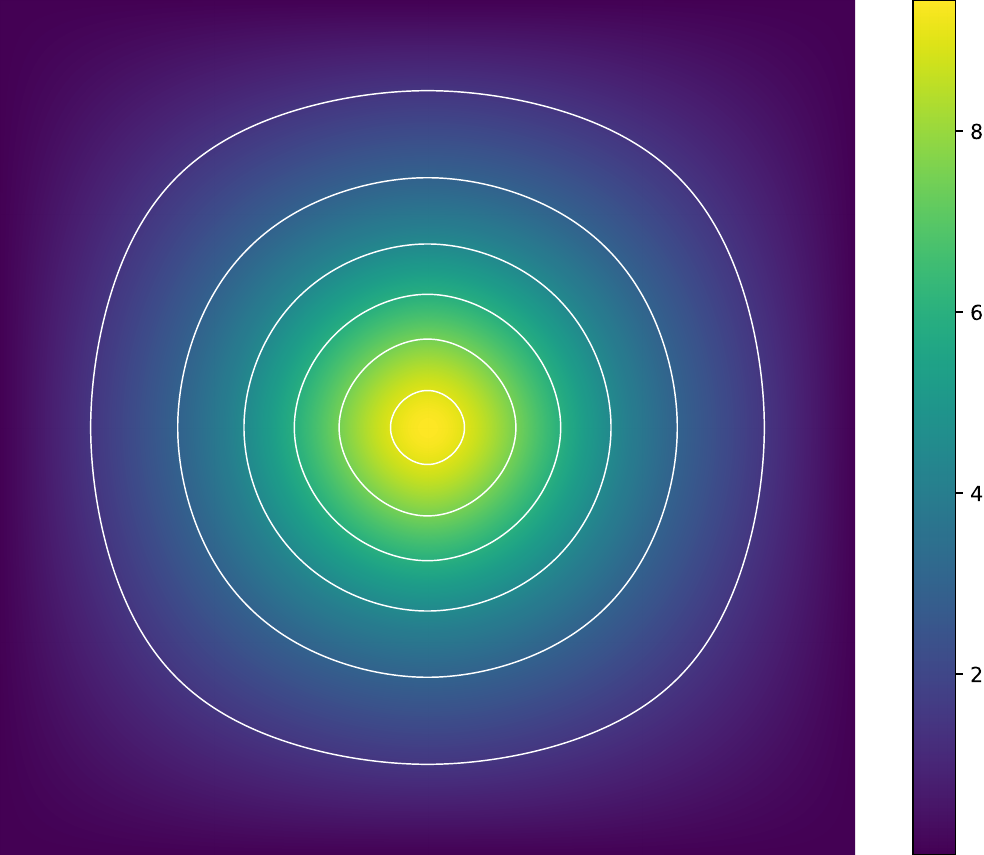}
    \hfill
        \includegraphics[width=0.35\textwidth, trim={0 0.3cm 0 0},clip]{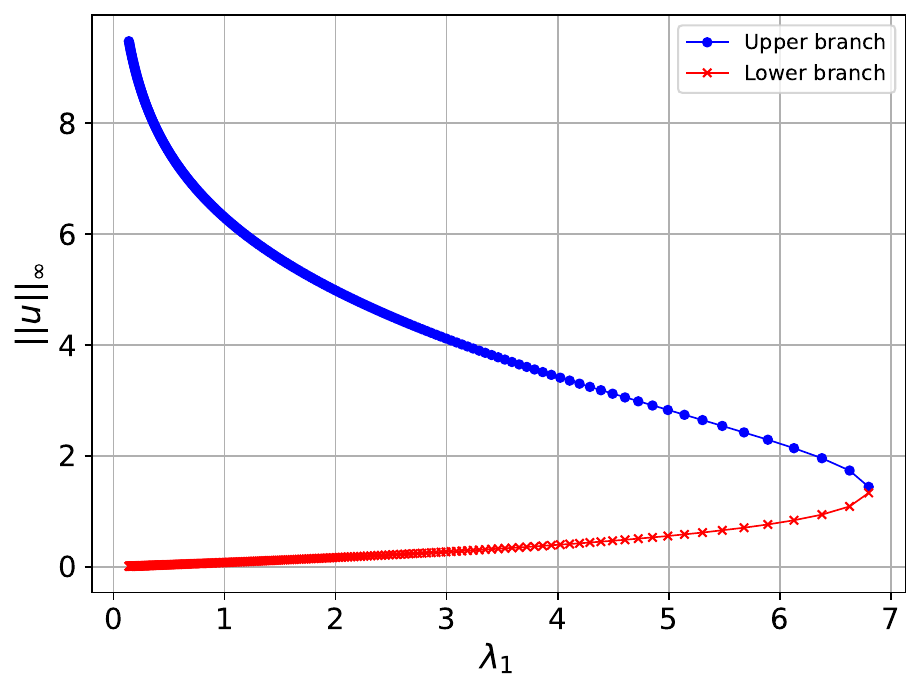}
    \caption{Bratu solutions with $d=2$ for $\lambda_1=0.14$, and saddle-node bifurcation diagram.}
    \label{fig:problem12}
\end{figure}

Let us start our discussion with the one-dimensional parameter space framework for the Bratu problem in Equation \eqref{bratu_multi}, with fixed values $\lambda_2=1$ and $\lambda_3=0$. 
In particular, we exploited the proposed strategy to obtain the bifurcation diagram by varying the parameter $\lambda_1$ with $d=1$ and $d=2$. 
Having obtained the two solutions near $\lambda_1^*$, we continue the solution path through deflated arclength with step $ds=0.2$ to simultaneously obtain the upper branch, and deflate it to retrieve the lower one. 
Figures \ref{fig:problem11} and \ref{fig:problem12} show the reconstruction of the bifurcation diagrams for the Bratu equation in the $\lambda_1$–$\|u\|_{\infty}$ plane, respectively for $d=1$ and $d=2$, with the corresponding coexisting lower and upper solutions.

The bifurcation behavior of the Bratu equation possesses the same qualitative properties between the one- and two-dimensional settings, with some notable differences. In one dimension, the saddle-node bifurcation occurs at a relatively small critical value $\lambda_1^*\approx 3.5$, while 
for $d=2$ the critical value increases to $\lambda_1^*\approx 6.8$, reflecting a stronger contribution needed to counterbalance the diffusive effect of the higher-dimensional Laplacian.

In this case, deflation is essential for recovering the upper solution branch, and arclength continuation significantly improves robustness compared to simple parameter continuation. Being able to automatically adjust the stepsize via arclength is even more important in the two-dimensional problem, requiring substantially higher computational effort, with increased degrees of freedom, larger numbers of solver iterations, and longer overall runtimes. For these reasons, subsequent explorations will be devoted to combining the developed strategy with Reduced Order Model approaches for the efficient reconstruction of the bifurcation diagram \cite{pichi_artificial_2023,pichi_deflation-based_2025}.  

\subsubsection{Allen--Cahn equation with $p=1$} 
\begin{figure}[tbp]
    \centering
        \includegraphics[width=0.45\textwidth]{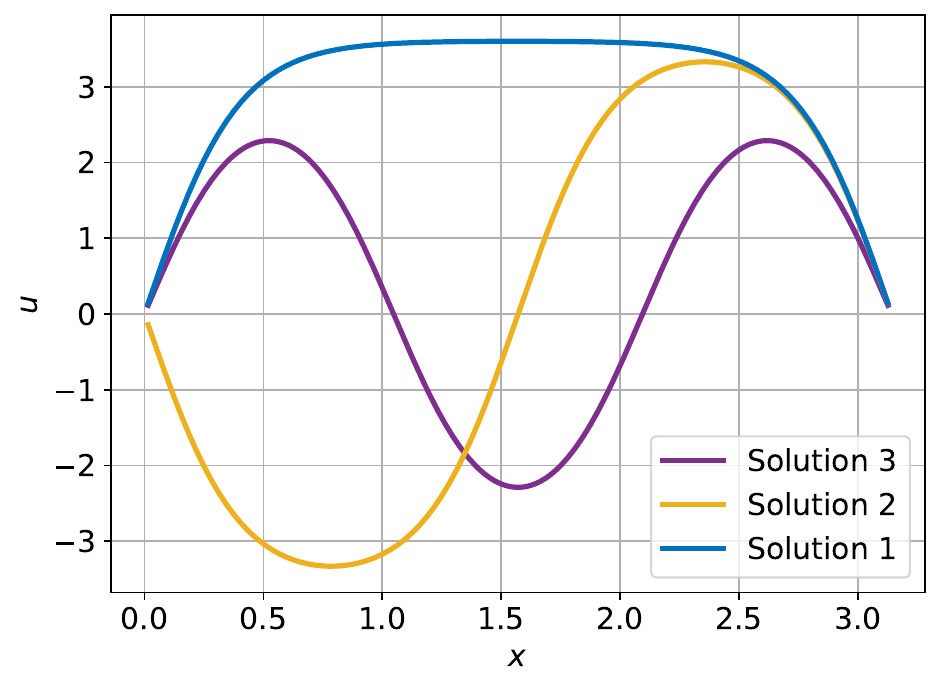} 
    \quad        
        \includegraphics[width=0.45\textwidth]{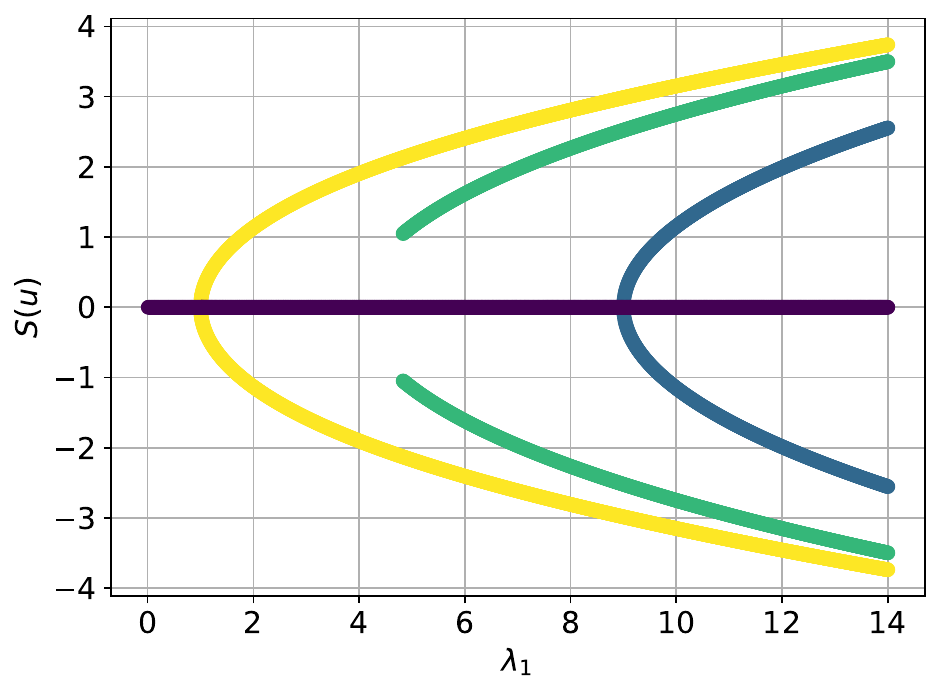}
    \caption{Allen--Cahn solutions with $d=1$ for $\lambda_1=13$, and pitchfork bifurcation diagram.}
    \label{fig:problem21}
\end{figure}

\begin{figure}[tbp]
    \centering
    \begin{minipage}{0.59\textwidth}
        \centering
        \includegraphics[width=0.32\textwidth]{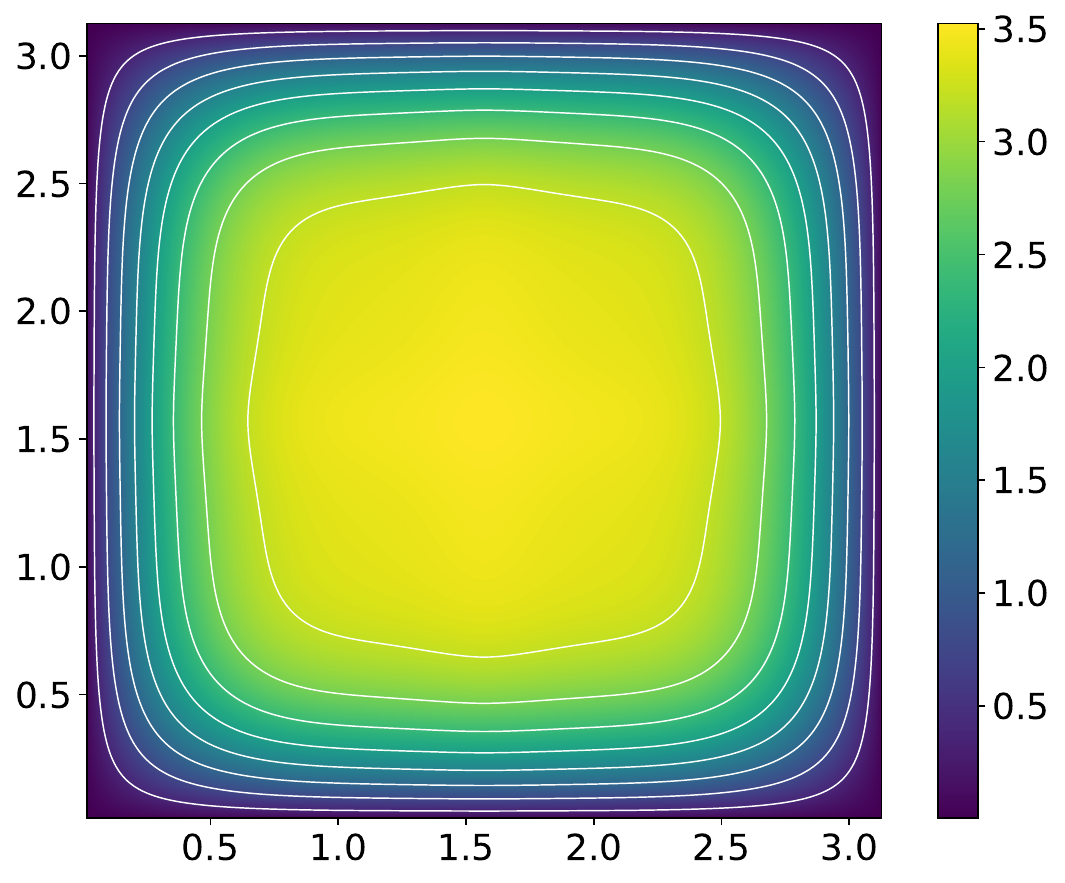}
        \hfill
        \includegraphics[width=0.32\textwidth]{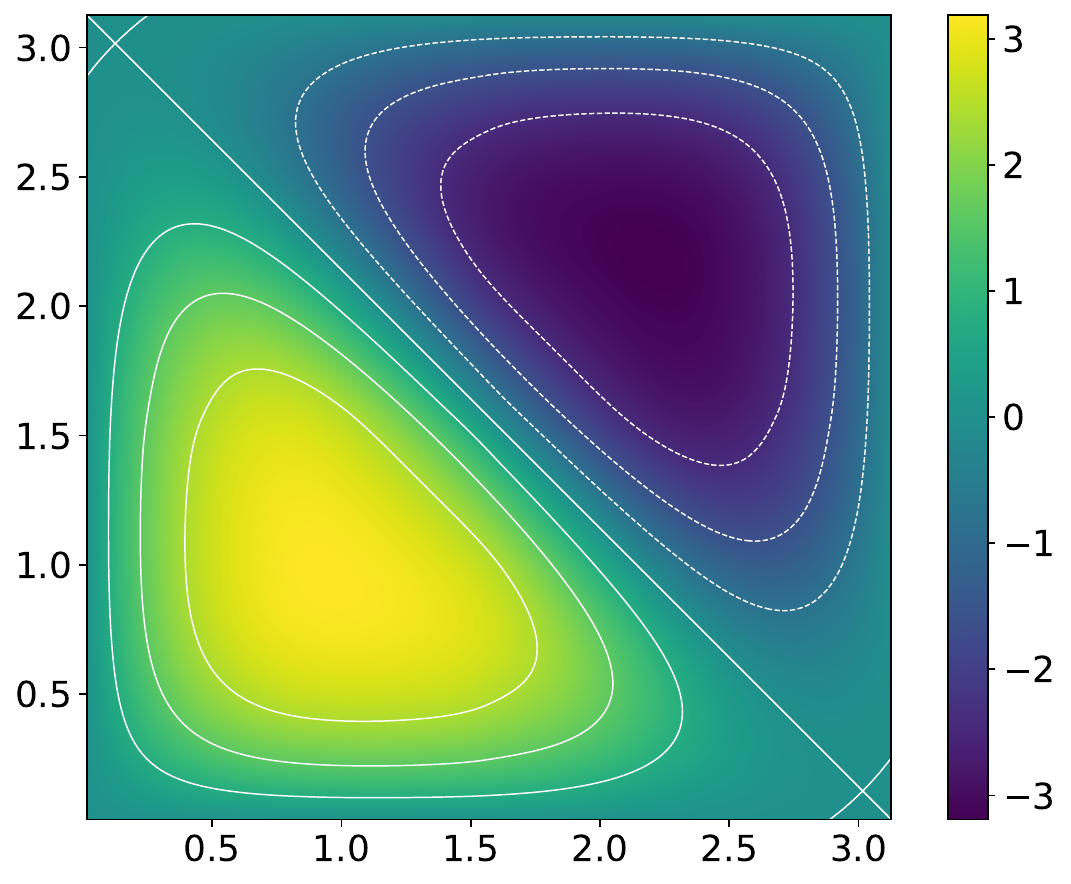}
        \hfill
        \includegraphics[width=0.32\textwidth]{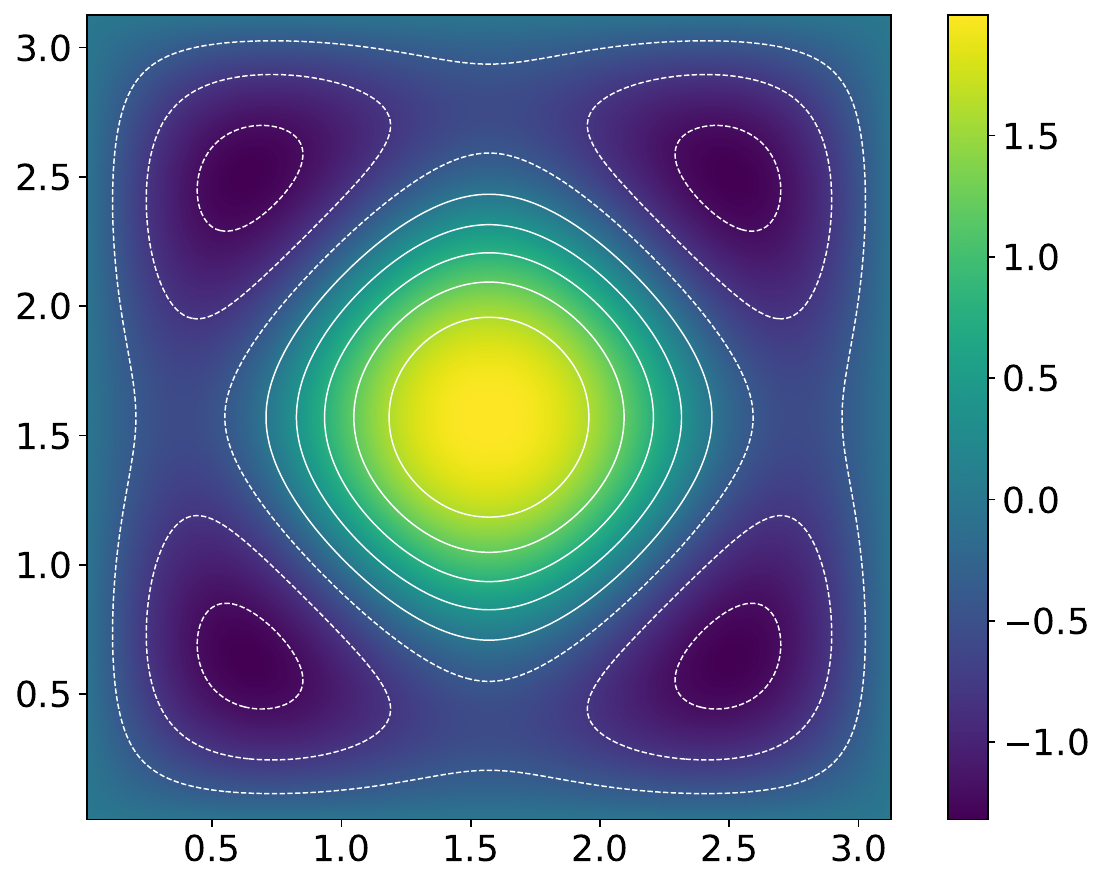}
        
        \includegraphics[width=0.32\textwidth]{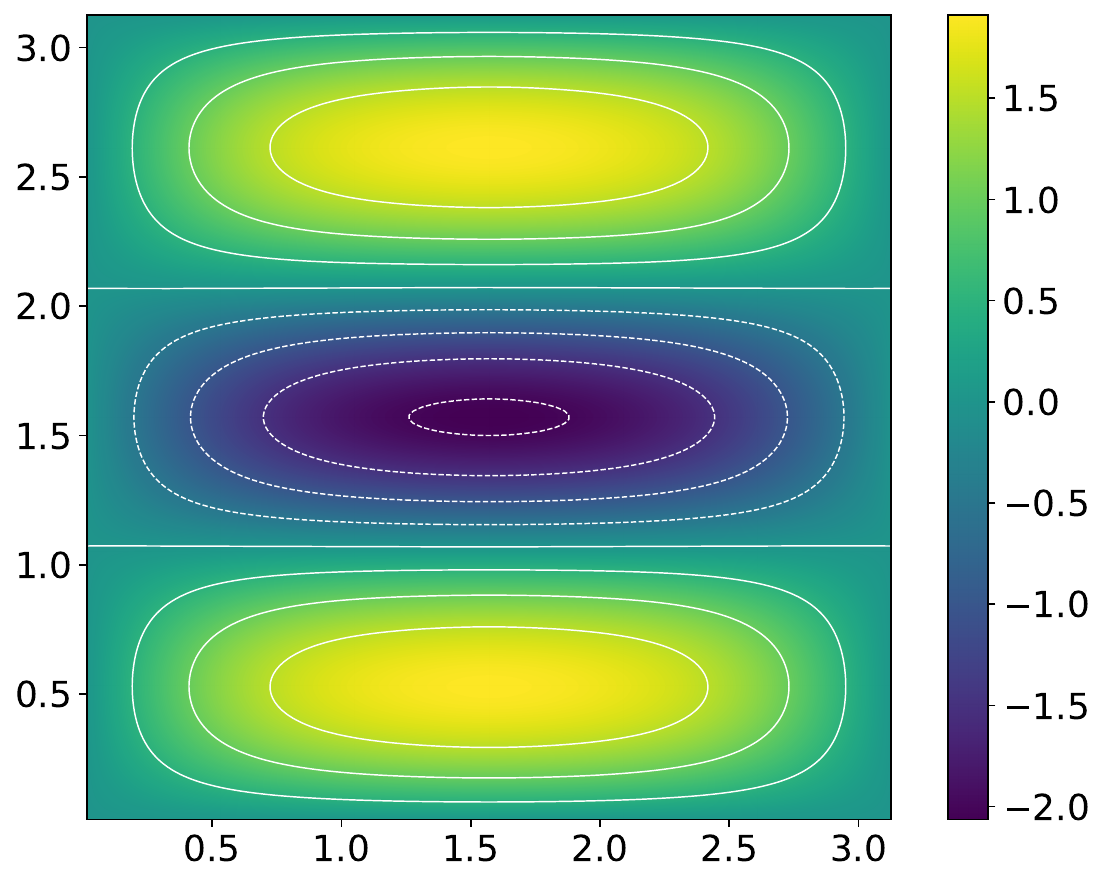}
        \quad
        \includegraphics[width=0.32\textwidth]{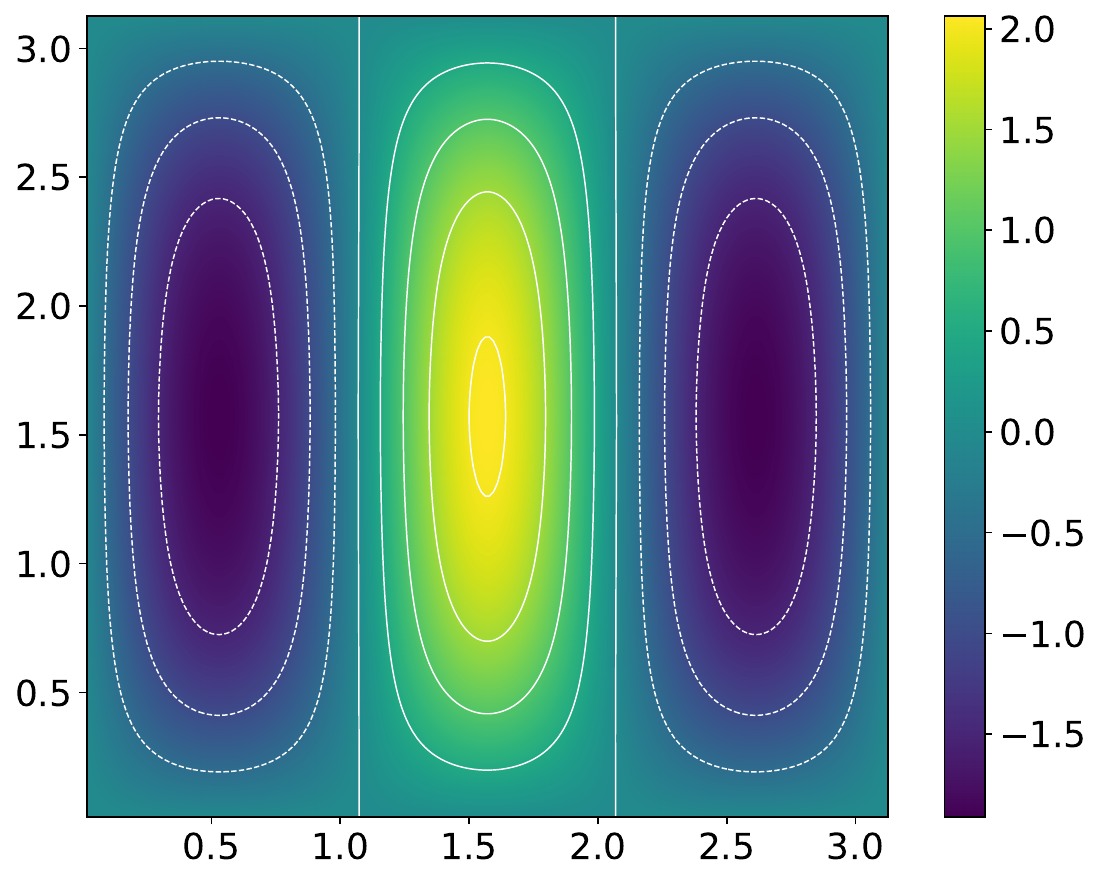}
    \end{minipage}\hfill
    \begin{minipage}{0.4\textwidth}
        \centering
        \includegraphics[width=\textwidth]{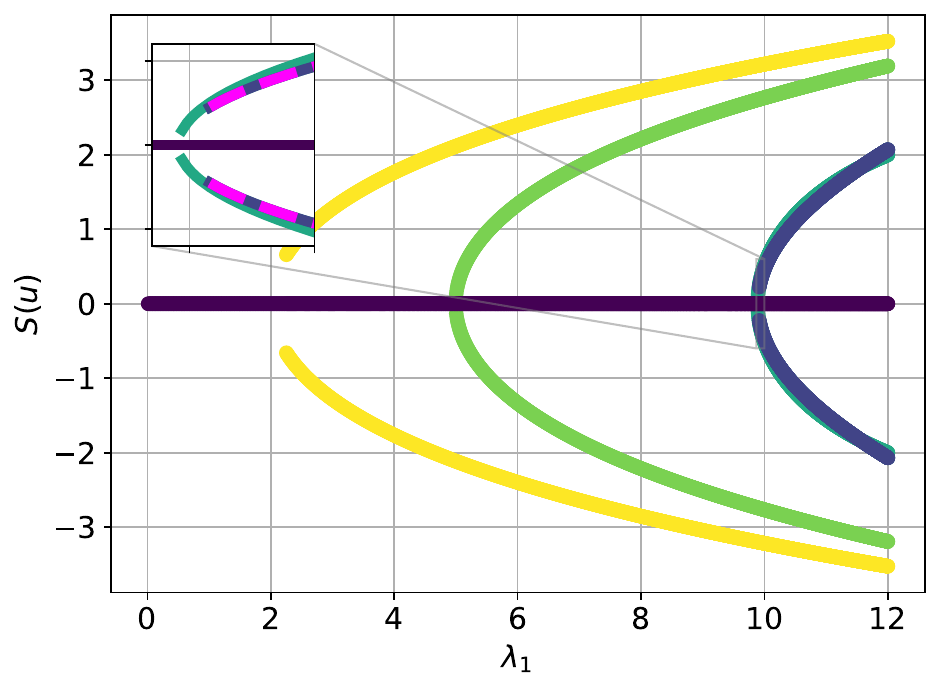}
    \end{minipage}        
    \caption{Allen–Cahn solutions with $d=2$ at $\lambda_1 = 12$, and pitchfork bifurcation diagram.}
    \label{fig:problem22}
\end{figure}

We now consider the Allen--Cahn equation with fixed values $\lambda_2=1$ and $\lambda_3=1$.
The deflated arclength continuation strategy is employed to compute the bifurcation diagrams by varying the parameter $\lambda_1$ with $d=1$ and $d=2$, depicting them in Figures \ref{fig:problem21} and \ref{fig:problem22} in the $\lambda_1$--$S(u)$ plane, where the scalar measure $S(u)$ is defined as $S(u) = \mathrm{sign}(u(x))\,\|u\|_{\infty}$, respectively for $x=2.19$ and $x=(0.02,\,2.19)$.

For the one-dimensional Allen--Cahn equation, given the higher complexity of the model, we employed  the deflated arclength continuation algorithm using an arclength step size of $ds = 0.01$. Figure \ref{fig:problem21} shows the recovered bifurcation diagram and three symmetric coexisting solutions for $\lambda_1 = 13$ belonging to branches with positive values of $S(u)$. Starting from the trivial solution branch, the algorithm detects the first bifurcation point at approximately $\lambda_1^* \approx 1$, and continuing along the solution path, additional bifurcation points are identified for $\lambda_1^* \approx 4$ and $\lambda_1^* \approx 9$.  We remark that, due to the intrinsic symmetry of the Allen--Cahn equation, the lower solution branches consist of symmetric solutions, which appear in pairs related by sign symmetry. These symmetric solutions persist across all bifurcation levels and form the lower branches of the bifurcation diagram.

A similar behavior is also observed in the two-dimensional case, where multiple bifurcation points have been obtained and depicted in the bifurcation diagram in Figure \ref{fig:problem21}. Here, the first two bifurcation points are located approximately at $\lambda_1^*\approx2$ and $\lambda_1^*\approx5$. Moreover, this setting is even more difficult since we see three branches originating from the third bifurcation point at approximately $\lambda_1^*\approx10$. 


In the two-dimensional case shown in Figure \ref{fig:problem22}, the bifurcation structure is even richer, with five coexisting solutions with positive values of $S(u)$ for $\lambda_1=12$. The corresponding symmetric solutions from the pitchfork phenomena with negative values of $S(u)$ also exist, reflecting the preserved symmetry of the problem.

In this setting, although the one- and two-dimensional Allen--Cahn equations exhibit bifurcation points at comparable parameter values, their resulting bifurcation diagrams differ significantly in structure and complexity. Indeed, three primary bifurcation points are observed in both cases as the parameter $\lambda_1$ increases. However, when for $d=1$, each of them gives rise to a simple pitchfork structure with two nontrivial symmetric solution branches emerging from the trivial one, while for $d=2$ the first two bifurcation points exhibit a similar qualitative behavior, but the third one originates three distinct solution branches, shown in the zoom-in of the right plot in Figure \ref{fig:problem22}. These branches are the result of a multiple pitchfork bifurcation, with two of them being characterized by the same value of $S(u)$, but correspond to solutions with different spatial patterns and symmetry properties.

From a numerical perspective, tracking all solutions branching from the third bifurcation point in the two-dimensional case is significantly more challenging. Standard continuation methods often may fail to capture all of them, as they tend to converge to the most stable branch. Consequently, robust techniques such as the one we proposed are required to reliably detect and compute multiple coexisting states. 

\subsubsection{Bratu equation with $p=2$} 
\begin{figure}[tbp]
    \centering

    \begin{subfigure}[b]{0.32\textwidth}
        \includegraphics[width=\textwidth]{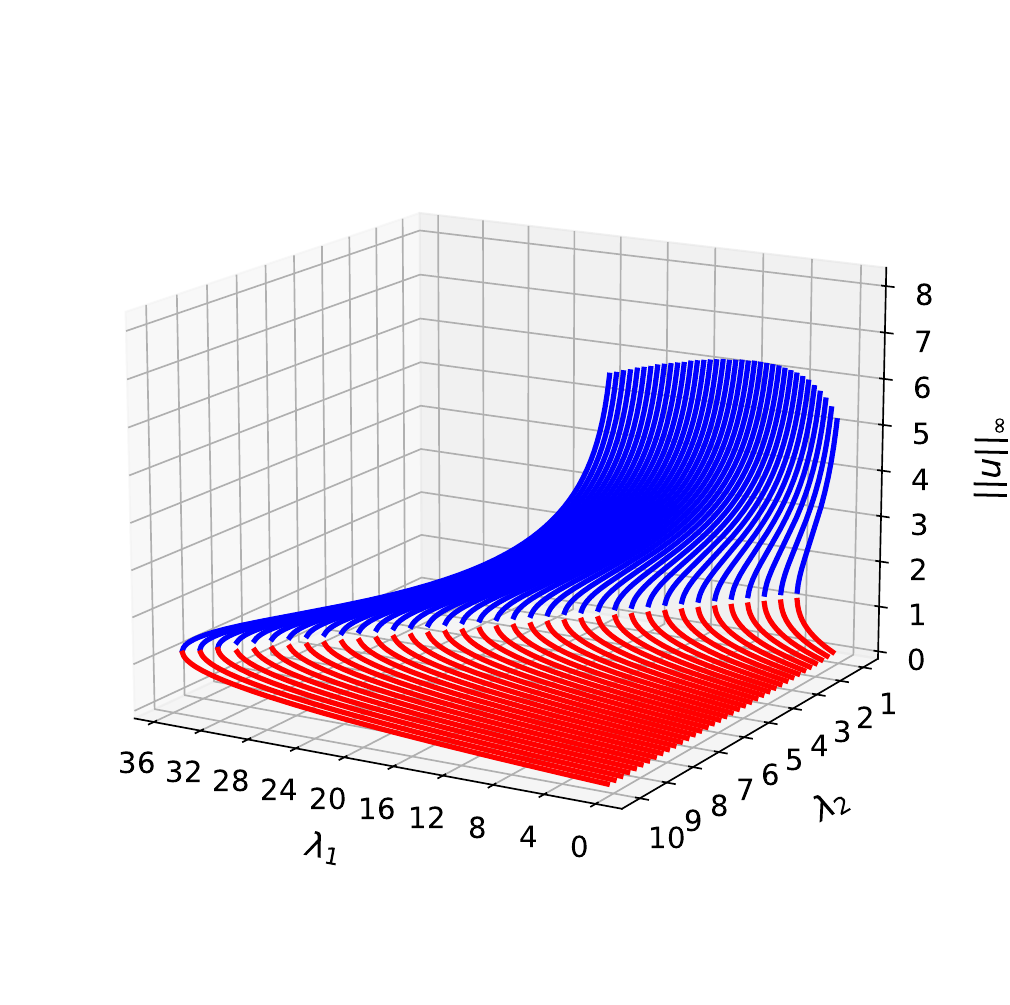}
        \caption{Horizontal lines}
        \label{bratu1d1}
    \end{subfigure}
    \hfill
    \begin{subfigure}[b]{0.32\textwidth}
        \includegraphics[width=\textwidth]{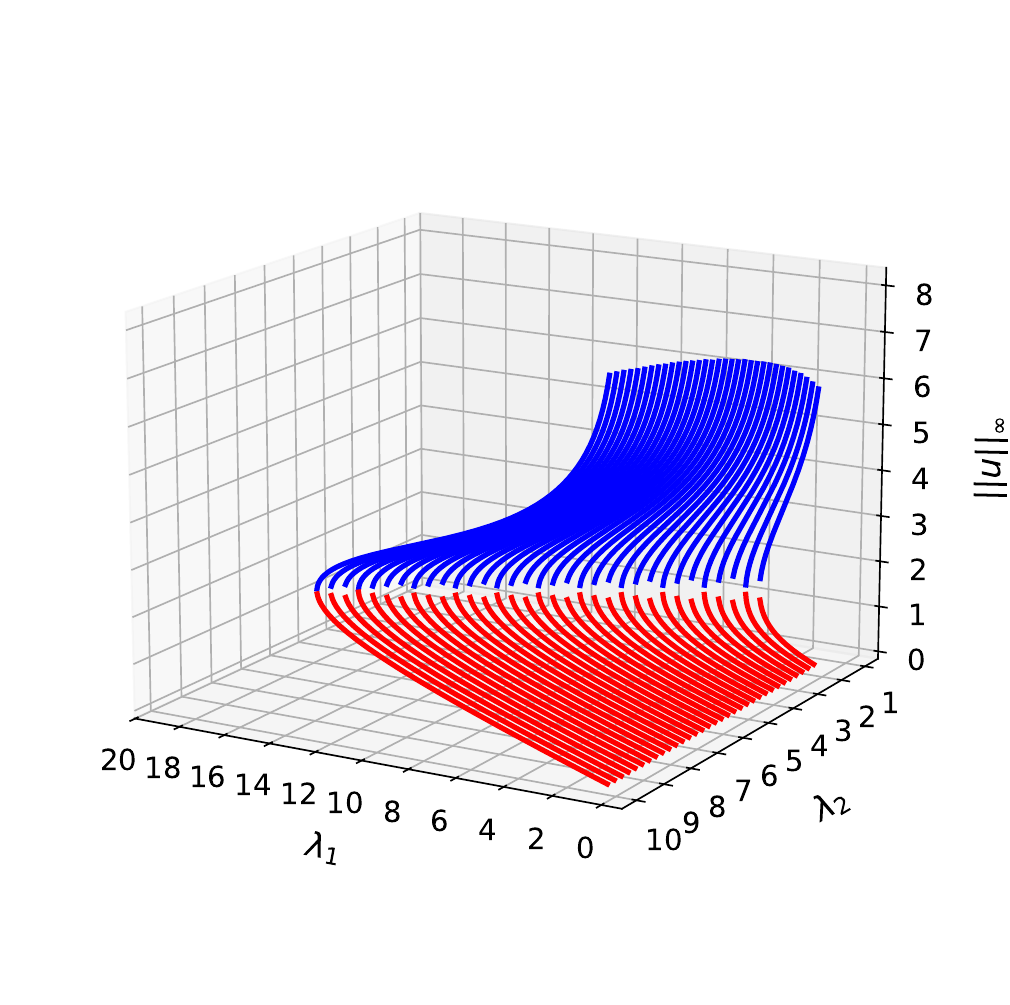}
        \caption{Diagonal lines}
        \label{bratu1d2}
    \end{subfigure}
    \hfill
    \begin{subfigure}[b]{0.32\textwidth}
            \includegraphics[width=\textwidth]{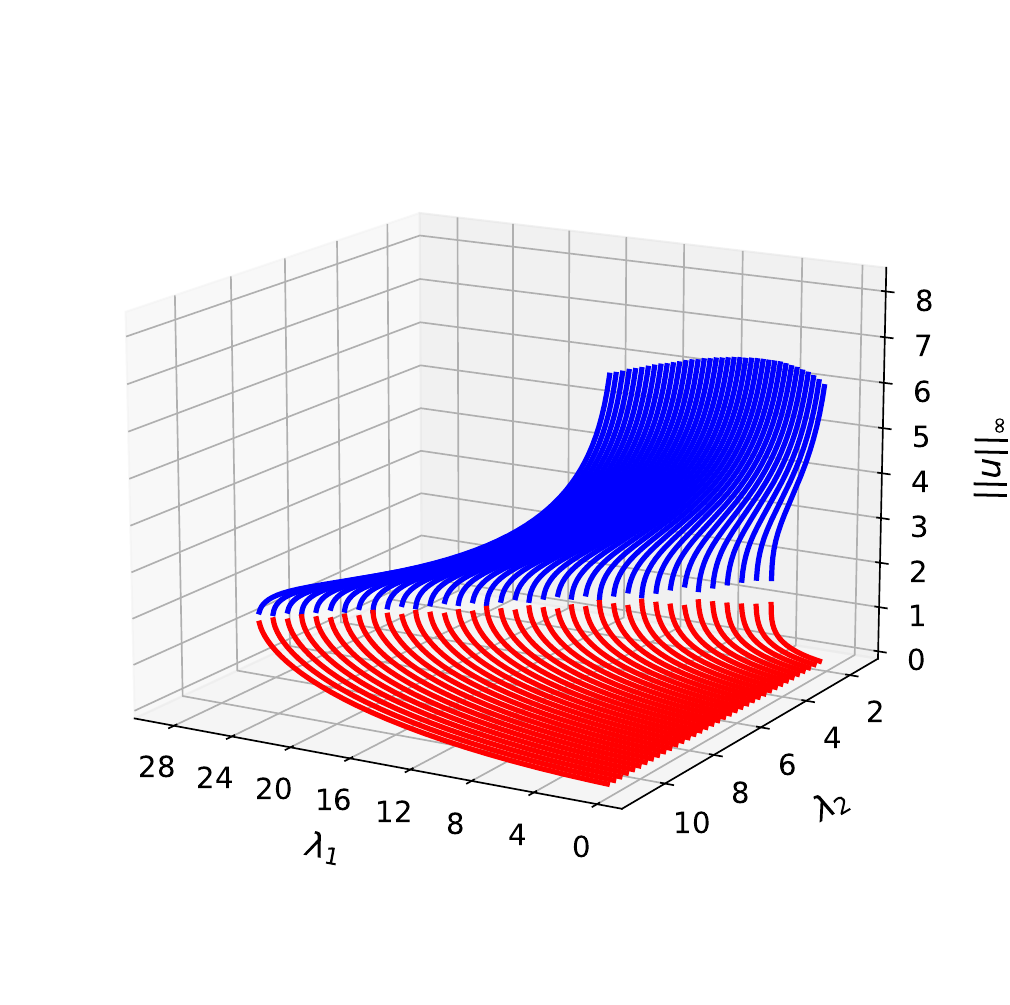}
        \caption{Elliptic curves}
        \label{bratu1d3}
    \end{subfigure}

    \caption{Bifurcation diagrams for Bratu with different paths and $d=1$.}
    \label{fig:problem31}
\end{figure}

\begin{figure}[tbp]
    \centering

    \begin{subfigure}[b]{0.32\textwidth}
        \includegraphics[width=\textwidth]{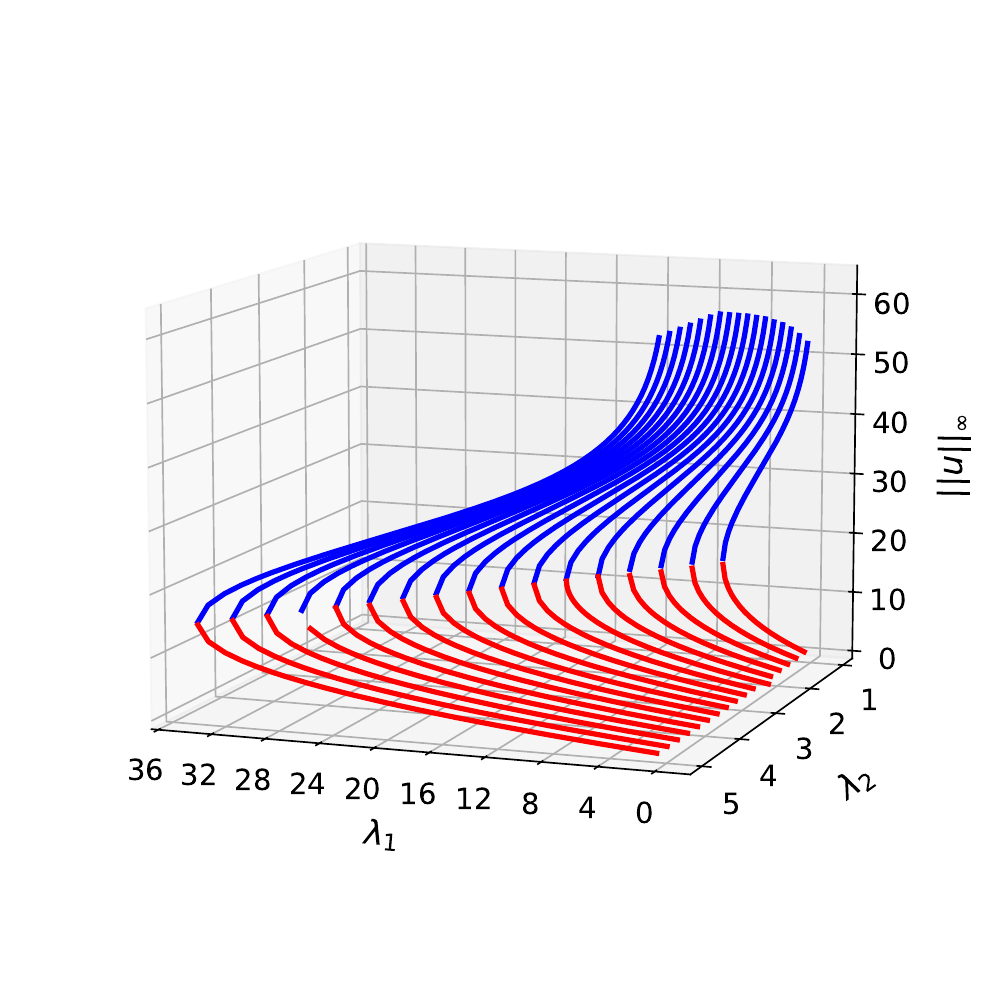}
        \caption{Horizontal lines}
        \label{bratu2d1}
    \end{subfigure}
    \hfill
    \begin{subfigure}[b]{0.32\textwidth}
        \includegraphics[width=\textwidth]{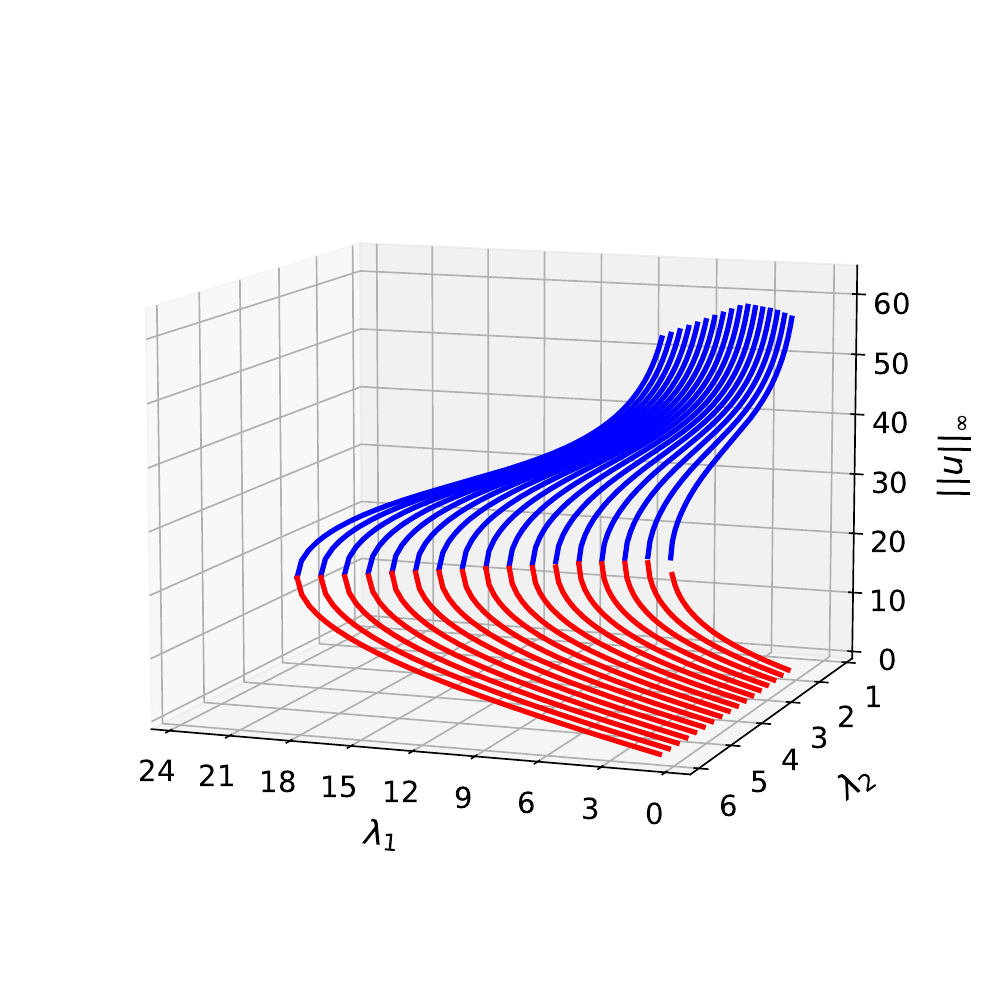}
        \caption{Diagonal lines}
        \label{bratu2d2}
    \end{subfigure}
    \hfill
    \begin{subfigure}[b]{0.32\textwidth}
        \includegraphics[width=\textwidth]{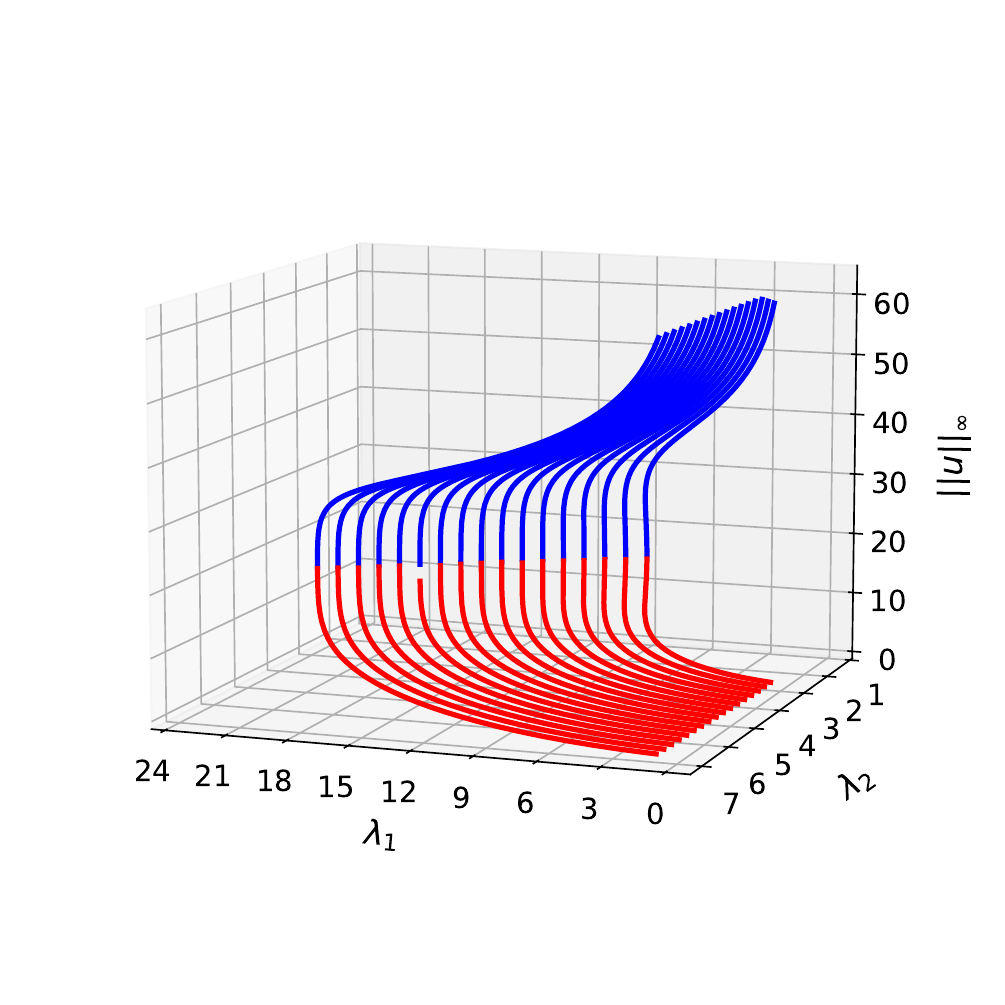}
        \caption{Elliptic curves}
        \label{bratu2d3}
    \end{subfigure}

    \caption{Bifurcation diagrams for Bratu with different paths and $d=2$.}
    \label{fig:Bratu2D}
\end{figure}

Having verified the performance of the proposed approach across several benchmarks in the single parameter scenario, here we aim at exploring the novel multiparametric setting, exploiting the deflation and arclength continuation strategies on the previously defined paths, to discuss the flexibility and robustness of the strategy under different parameter couplings. 
We start with the multiparametric Bratu problem in Equation \eqref{bratu_multi} with fixed $\lambda_3=0$, for which we exploited the multiparametric deflated arclength continuation technique, discussed in Section \ref{sec3.4.2}, with $ds = 0.2$ and different continuation paths $g(\lambda_1,\lambda_2)$. 
As illustrated in Figures \ref{fig:problem31} and \ref{fig:Bratu2D}, respectively for $d=1$ and $d=2$, all path strategies perform efficiently, successfully capturing both upper and lower solution surfaces. The computed solutions are once again seen in the $(\lambda_1, \lambda_2)$-plane, where the saddle-node bifurcation curves are clearly visible.
This further demonstrates the robustness of the proposed technique in exploring complex solution landscapes even in the multiparametric setting

We remark that different continuation paths correspond to different regions explored in the $(\lambda_1,\lambda_2)$ parameter space, but a nice property of our approach is that the resulting bifurcation diagrams are qualitatively and quantitatively similar across all paths. In particular, for every path considered, both the lower and upper solution surfaces terminate on the same oblique line, which corresponds to the bifurcation curve.

\begin{figure}[tbp]
    \centering

    \begin{subfigure}[b]{0.4\textwidth}
        \includegraphics[width=\textwidth]{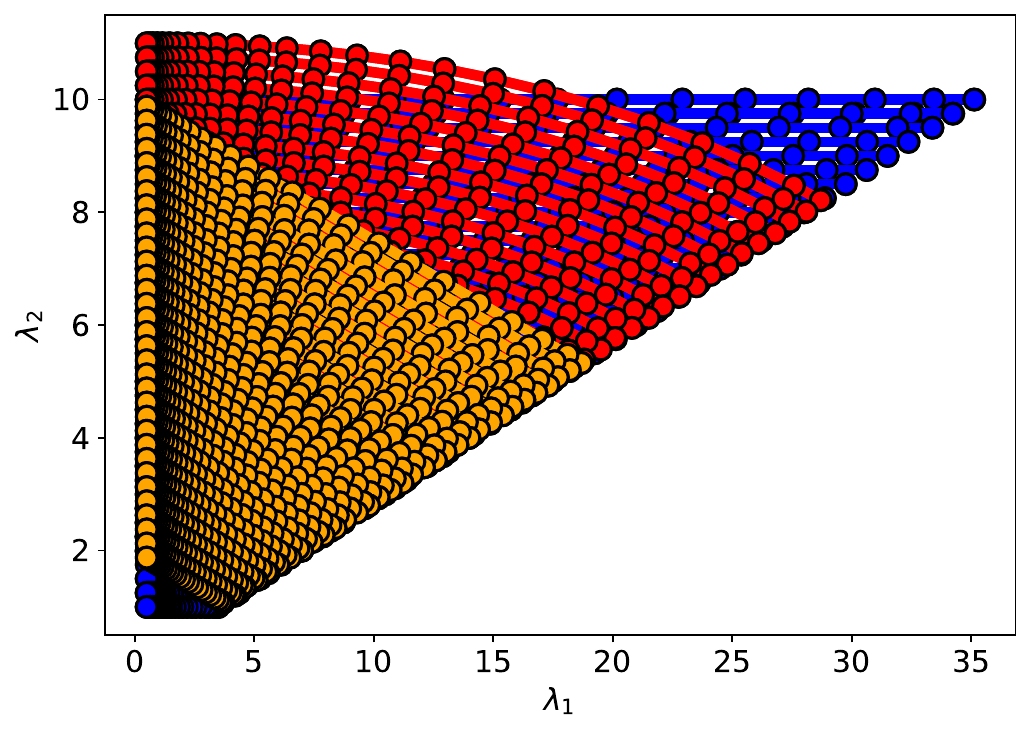}
        \caption{$d=1$}
        \label{fig:subfigure341}
    \end{subfigure}
    \hfill
    \begin{subfigure}[b]{0.4\textwidth}
        \includegraphics[width=\textwidth]{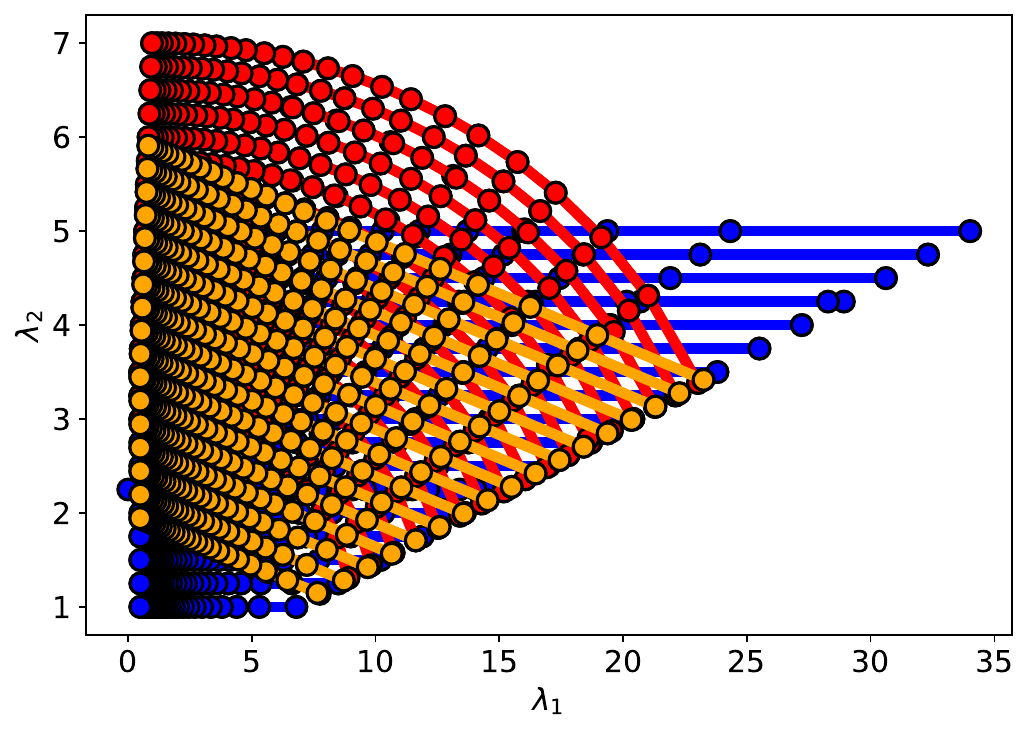}
        \caption{$d=2$}
        \label{fig:subfigure342}
    \end{subfigure}

    \caption{Paths in the parameter space with Bratu bifurcation curve.}
    \label{fig:problem3path}
\end{figure}

The corresponding projections of the paths onto the $(\lambda_1,\lambda_2)$ plane are presented in Figure \ref{fig:problem3path}, for the one- and two-dimensional Bratu equations, left and right respectively. These projections clearly show that all continuation paths identify the same critical bifurcation region, as the detected bifurcation points consistently align along the oblique bifurcation curve. This path-independence demonstrates the robustness and consistency of the proposed multiparametric deflated arclength strategy, while also allowing flexible exploration of different regions of the parameter space without loss of accuracy.

\subsubsection{Allen--Cahn equation with $p=2$} 
Here we consider the multiparametric Allen--Cahn problem in Equation \eqref{bratu_multi} with fixed $\lambda_3=1$, and verify the performance of the deflated arclength strategy to reconstruct the bifurcation diagrams in the multiparametric setting with different continuation paths. Once again, all path strategies perform effectively, even dealing with the more complex pitchfork phenomena, with the methodology being able to identify multiple bifurcation curves from which serveral coexisting solutions originate. 

Exploiting the same paths described before, we explore the parameter space detecting the first bifurcation point and the corresponding upper solution branch with an arclength step size $ds=0.01$. The deflation is then allowing for the computation of the remaining solution branches, yielding an accurate and comprehensive reconstruction of the bifurcation diagrams for both $d=1$ and $d=2$.

\begin{figure}[tbp]
    \centering

    \begin{subfigure}[b]{0.32\textwidth}
        \includegraphics[width=\textwidth]{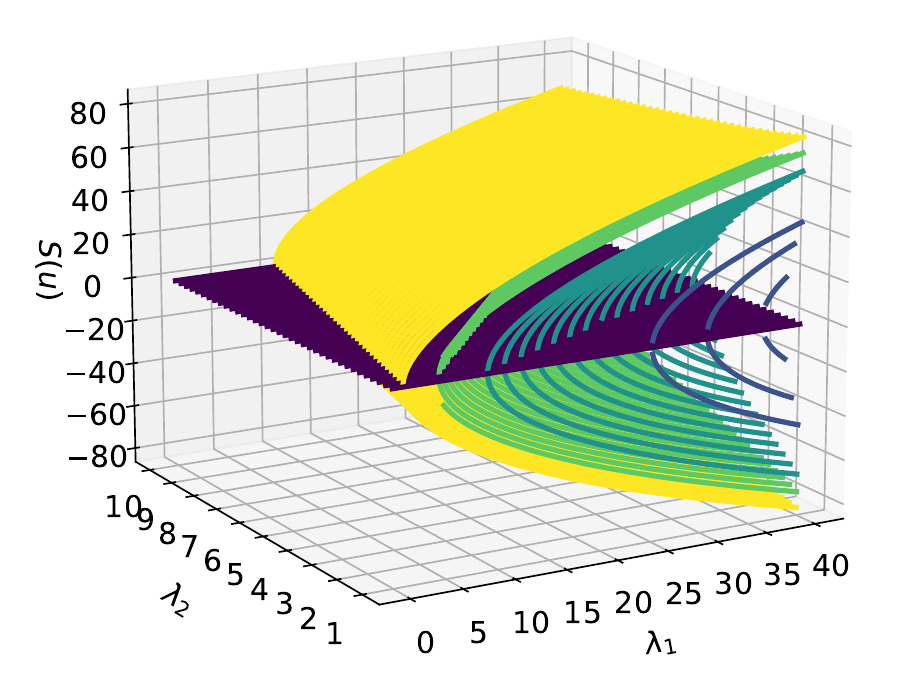}
        \caption{Horizontal lines}
        \label{allen1d1}
    \end{subfigure}
    \hfill
    \begin{subfigure}[b]{0.32\textwidth}
        \includegraphics[width=\textwidth]{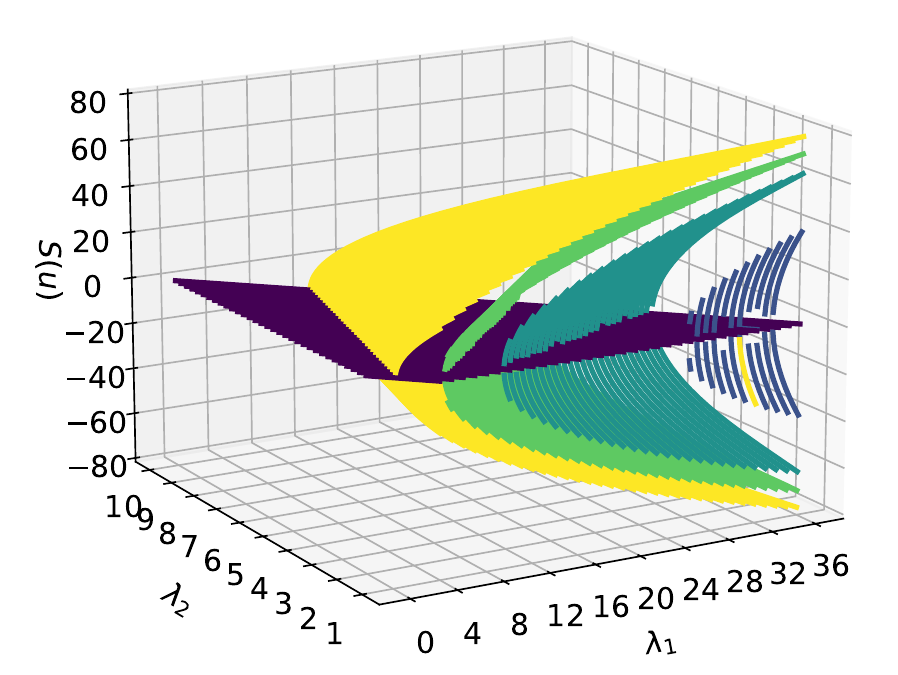}
        \caption{Diagonal lines}
        \label{allen1d2}
    \end{subfigure}
    \hfill
    \begin{subfigure}[b]{0.32\textwidth}
        \includegraphics[width=\textwidth]{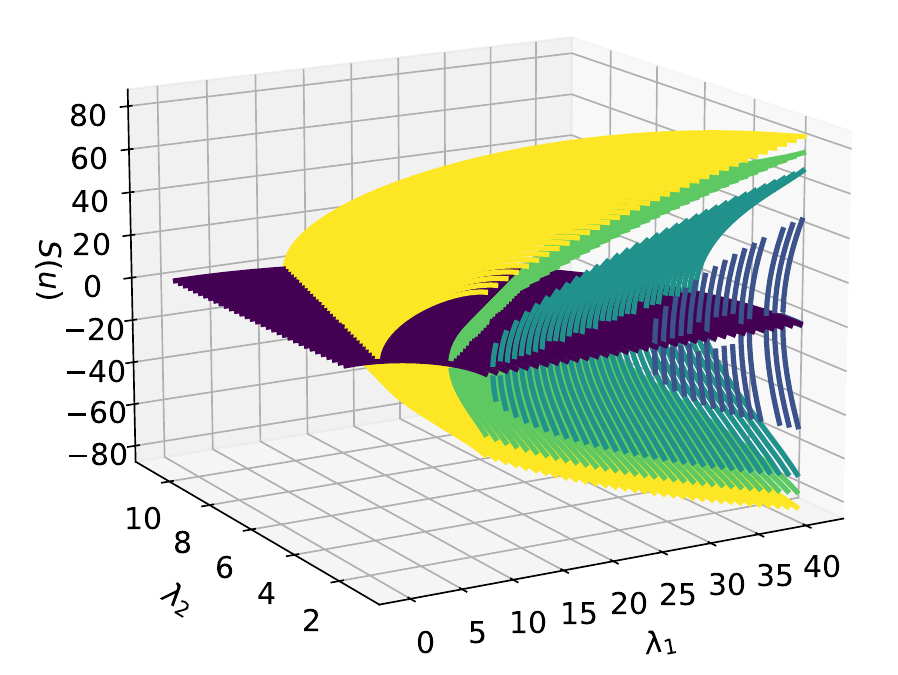}
        \caption{Elliptic curves}
        \label{allen1d3}
    \end{subfigure}

    \caption{Bifurcation diagrams for Allen--Cahn with different paths and $d=1$.}
    \label{fig:problem41}
\end{figure}

In particular, we show in Figure \ref{fig:problem41} the bifurcation diagrams for the one-dimensional Allen--Cahn equation with different continuation paths, where multiple bifurcation curves are detected, each giving rise to distinct solution surfaces. As expected, each successive bifurcation curve has an output value $S(u)$ smaller than the previous ones, accounting for the "higher frequency" of the following modes. 

\begin{figure}[tbp]
    \centering

    \begin{subfigure}[b]{0.32\textwidth}
        \includegraphics[width=\textwidth]{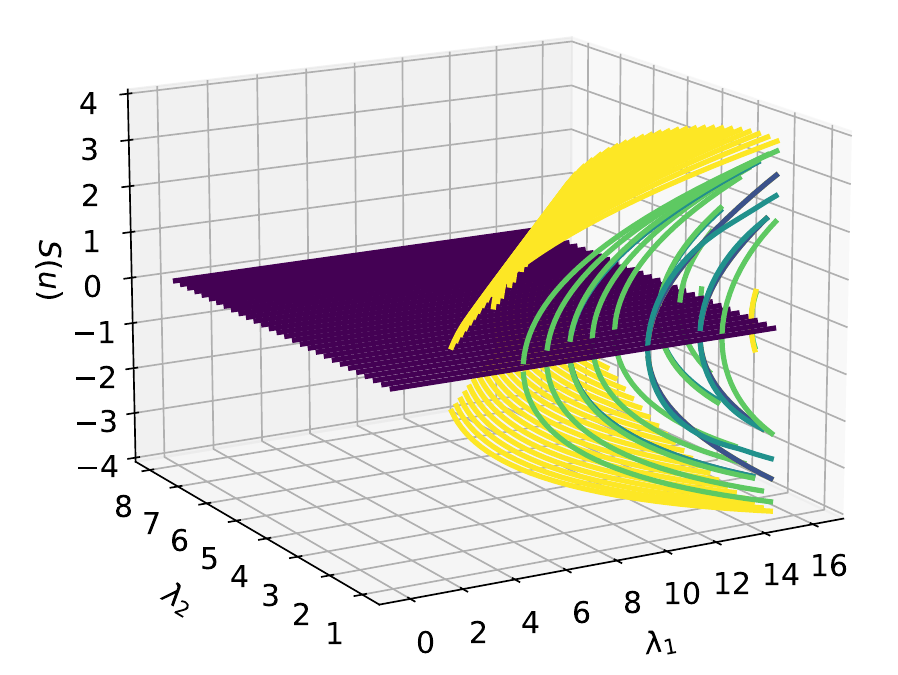}
        \caption{Horizontal lines}
        \label{allen2d1}
    \end{subfigure}
    \hfill
    \begin{subfigure}[b]{0.32\textwidth}
        \includegraphics[width=\textwidth]{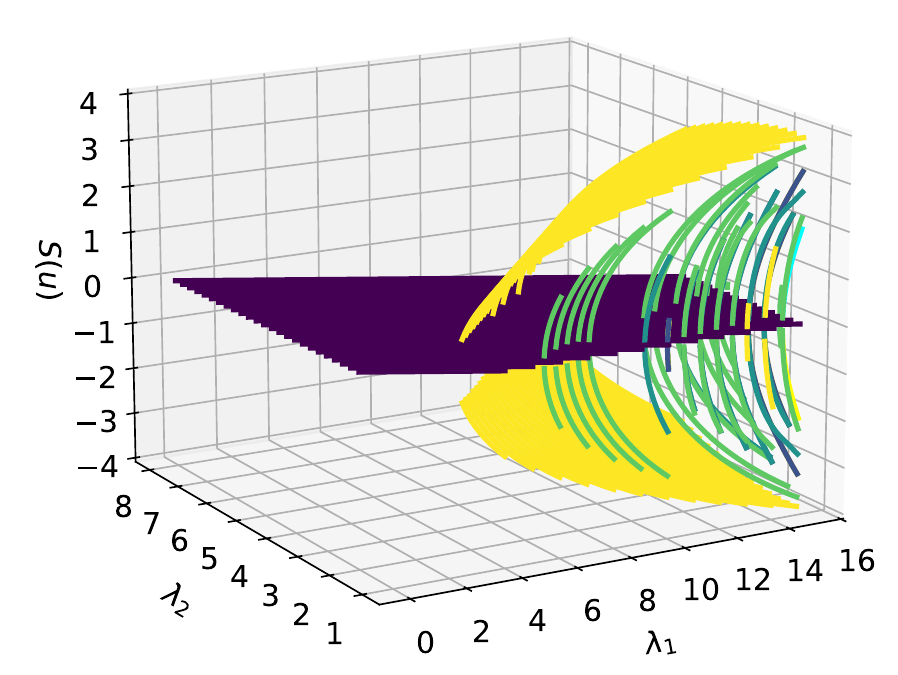}
        \caption{Diagonal lines}
        \label{allen2d2}
    \end{subfigure}
    \hfill
    \begin{subfigure}[b]{0.32\textwidth}
            \includegraphics[width=\textwidth]{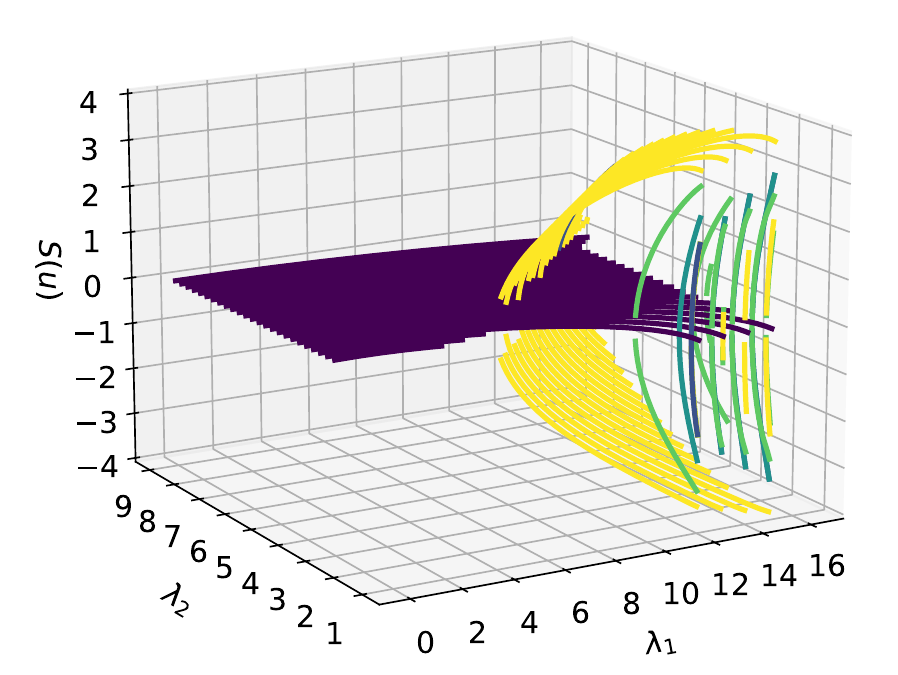}
        \caption{Elliptic curves}

        \label{allen2d3}
    \end{subfigure}

    \caption{Bifurcation diagrams for Allen--Cahn with different paths and $d=2$.}
    \label{fig:allen2D}
\end{figure}

Similarly, Figure \ref{fig:allen2D} shows the bifurcation diagrams for the two-dimensional case. Given the higher complexity of the multiparametric case for $d=2$, also exhibiting multiple bifurcation, the methodology shows slight less robustness especially for larger values of the critical parameter $\lambda_1$. 
There, three solution branches emerge from the third bifurcation point, with similar output values $S(u)$ that make it more challenging to visually distinguish the different surfaces and their reconstruction through deflation.

\subsection{Multiparametric bifurcation curve detection}
Finally, in this section we aim at focusing solely on the reconstruction of the bifurcating regions, rather than the entire bifurcation diagram, exploiting the zigzag methodology outlined in Section \ref{sec3.4.3}.
In particular, we use the Bratu and Allen--Cahn problems defined before as illustrative examples for $p=2$ and $p=3$, highlighting the effectiveness of the proposed approach in identifying transition regions in multiparametric systems. 

\subsubsection{Detecting bifurcating curves for $p=2$} 
Let us start by considering the saddle-node bifurcating phenomena for Bratu problem with two parameters and $d=1,2$. To trace the bifurcation path with the zigzag strategy, we set the angle in the path function, as defined in Equation \eqref{sec343.eq1}, as $\theta = \frac{\pi}{20}$.

\begin{figure}[tbp]
    \centering

    \begin{subfigure}[b]{0.32\textwidth}
        \includegraphics[width=\textwidth]{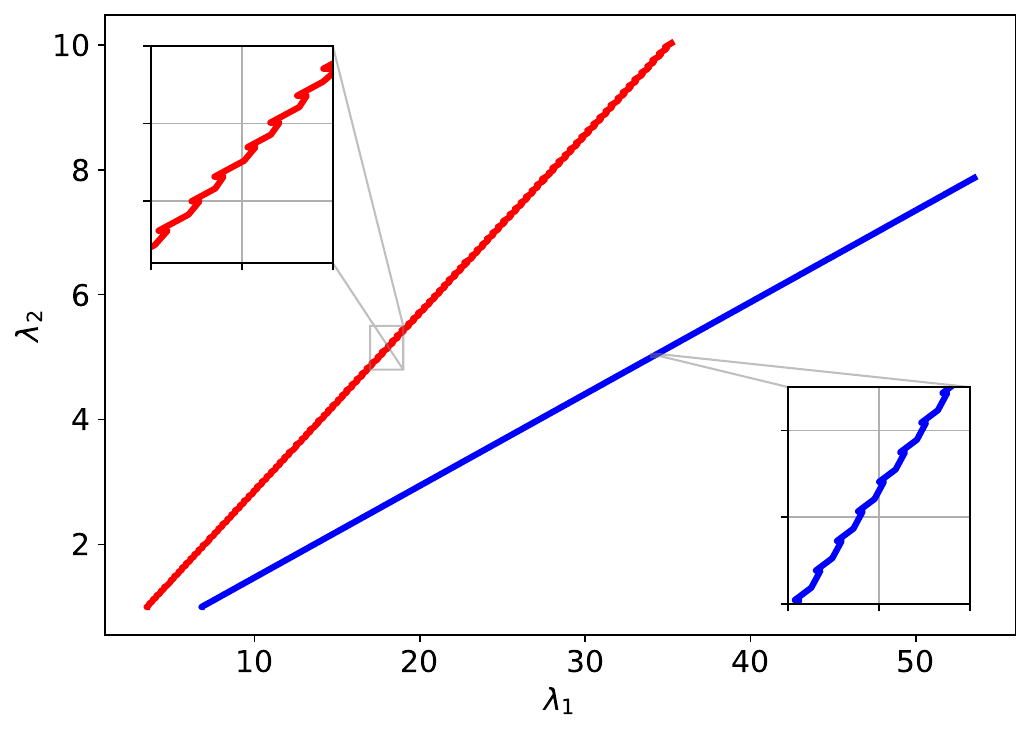}
        \caption{Bratu}
        \label{fig:det_bratu}
    \end{subfigure}
    \hfill
    \begin{subfigure}[b]{0.32\textwidth}
        \includegraphics[width=\textwidth]{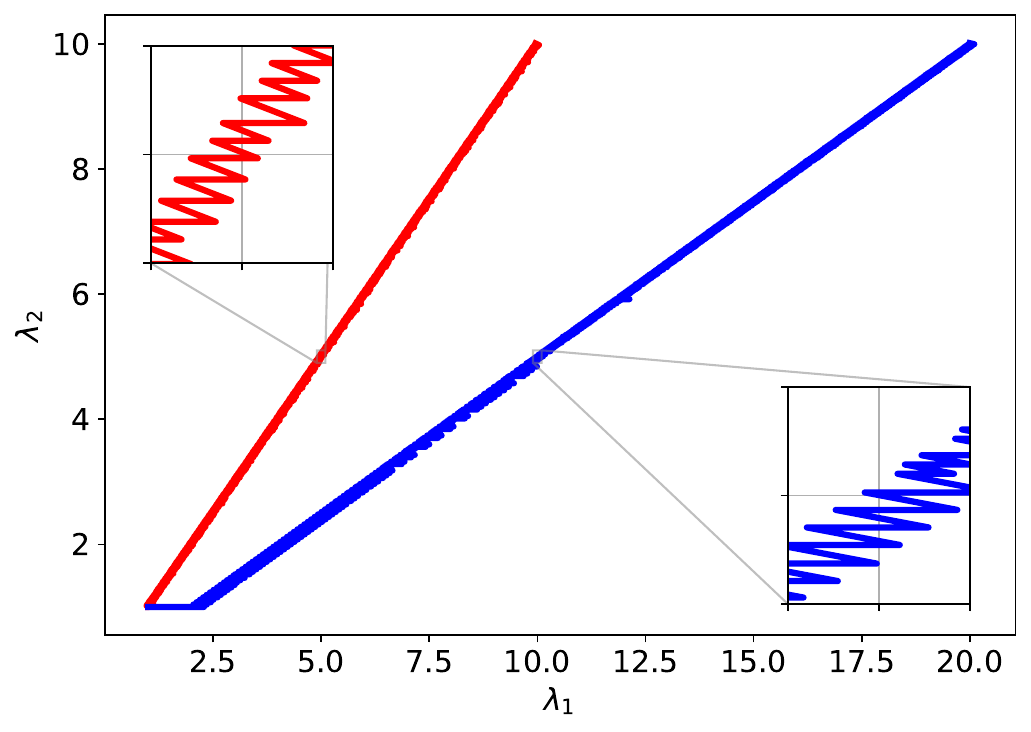}
        \caption{Allen–Cahn}
        \label{fig:det_Allen}
    \end{subfigure}
    \hfill
    \begin{subfigure}[b]{0.32\textwidth}
        \includegraphics[width=\textwidth]{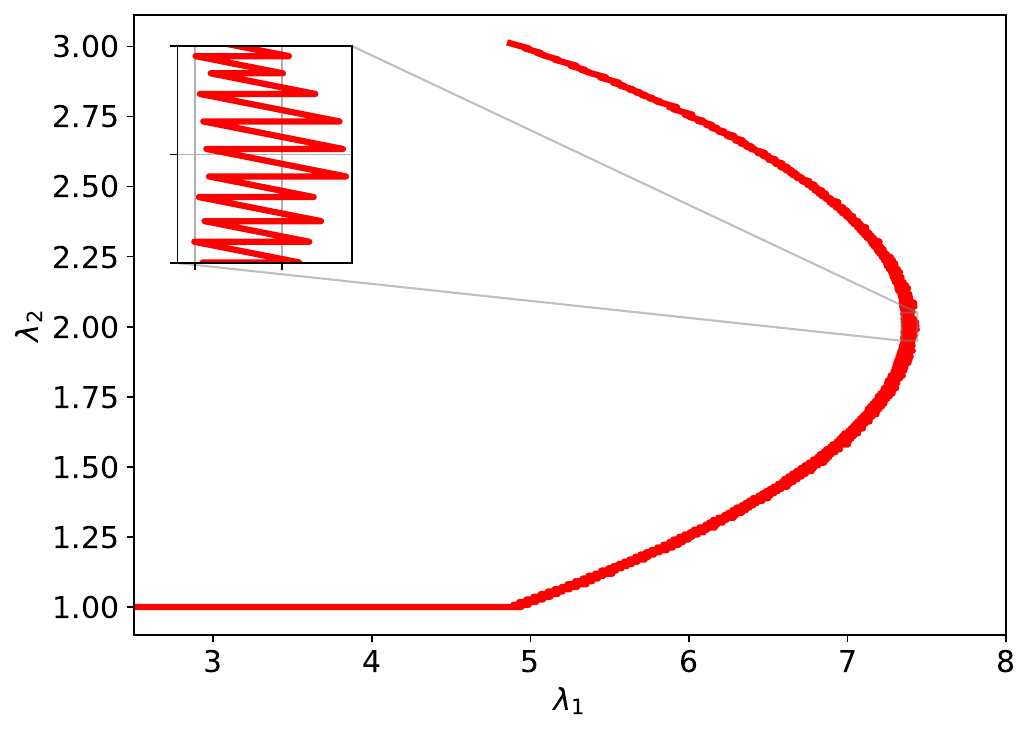}
        \caption{Modified Allen--Cahn}
        \label{fig:det_non_bif1D}
    \end{subfigure}

    \caption{Bifurcation curve detection for different benchmarks and spatial dimension, i.e.\ $d=1$ (red) and $d=2$ (blue).}
    \label{fig:detection1}
\end{figure}

The procedure begins with the choice of the initial path $g(\lambda_1, \lambda_2) = \lambda_2 - 1$. Once the two solutions in the vicinity of the bifurcation point are identified, we switch to arclength continuation with step size $ds = 0.01$, employing deflation and continue exploring the current direction for a fixed amount of steps $k = 5$. Then, the path is rotated by the angle $\theta$, and the deflated continuation strategy is performed to cross again the bifurcating curve identifying its behavior, i.e.\ until we reach again the non-existence region, and so on. The results of such an approach are illustrated in Figure \ref{fig:det_bratu}, where both straight bifurcation curves, for $d=1$ and $d=2$, have been correctly reconstructed with great precision.

We now continue with the detection of the pitchfork bifurcation curve for the Allen--Cahn problem with $p=2$. Similarly as before, the zigzag approach is exploited with $\theta = \pi - \frac{\pi}{20}$, considering the relative direction of the arclength continuation with $\lambda_2 = 1$. As the continuation procedure crosses from the unique solution region to the multiple solution region with three coexisting states, a bifurcation point is detected, and the initial stages of the bifurcation diagram are discovered for $ds = 0.01$ and $k=5$.
By iterative repetition of this zigzag procedure, the algorithm shows great robustness in tracing the bifurcation curves of the Allen--Cahn problem for both $d=1$ and $d=2$, as depicted in Figure \ref{fig:det_Allen}.

We remark that, until now, the chosen parameterization for both benchmarks essentially resulted in the translation of the first bifurcation point. To provide an even more challenging task to our methodology, rather than detecting linear bifurcation curves, we also consider a custom modification of the Allen--Cahn problem in $\Omega= [0,1]$, where the diffusion coefficients multiplying the laplacian is given as the quadratic term in the parameter $\rho(\lambda_2)=-(\lambda_2-1)^2+3$, and the boundary conditions are modified as $u_x(0)=0$ and $u(1)=0$ for $d=1$. In this way, we aim at manually constructing a benchmark featuring a nonlinear bifurcation curve, and test the properties of the proposed strategy even for less trivial bifurcation curves, that also characterize complex phenomena in fluid dynamics \cite{pichi_artificial_2023, PichiGraphConvolutionalAutoencoder2024}.
Exploiting the same setting as the previous Allen--Cahn case, we detect the bifurcation point and proceed with the zigzag pattern back and forth across the curved bifurcating region resembling the quadratic pattern chosen for the diffusion coefficient. Once again, as illustrated in Figure \ref{fig:det_non_bif1D}, the complete behavior of the critical curve is nicely and effectively recovered, providing an efficient way to discover the most interesting region, also available as an important insight for subsequent investigations, with no need for a-priori discretization of the parametric grid and exhaustive search on all possible system's configurations.

\subsubsection{Detecting bifurcating regions for $p=3$} 
With this last investigation, we aim at showing the capability of the zigzag strategy in combination with the deflated arclength continuation when considering a three-dimensional parameter space, i.e.\ $p=3$. Indeed, in such context, reconstructing all admissible solution branches would be unbearable from the computational point of view, while an efficient discovery of the critical surface defined by the evolution of the first bifurcation point in the 3D space could reveal as the key ingredient as a preprocessing step or for a more detailed exploration at a reduced cost.

\begin{figure}[tbp]
    \centering
    \begin{subfigure}[b]{0.45\textwidth}
        \includegraphics[width=\textwidth]{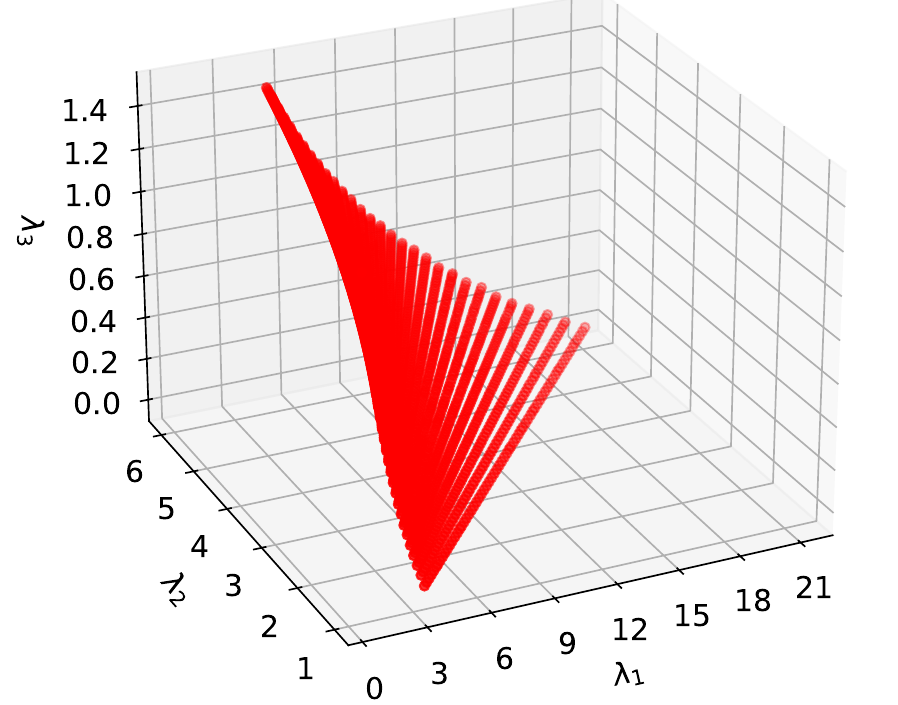}
        \label{fig:det_3par_Bratu1D}
    \end{subfigure}
    \quad
    \begin{subfigure}[b]{0.45\textwidth}
        \includegraphics[width=\textwidth]{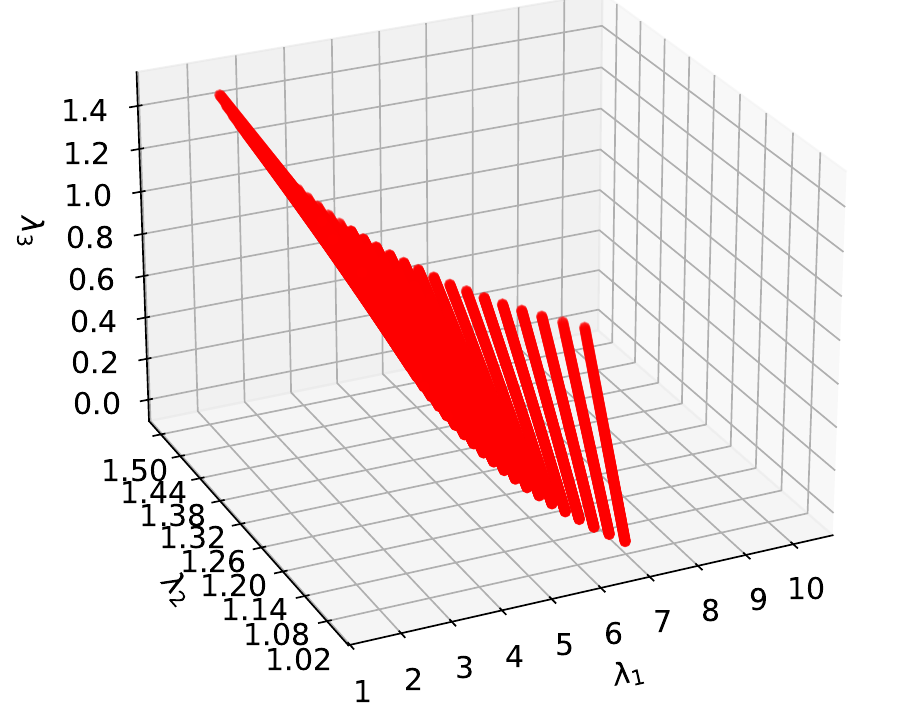}
        \label{fig:det_3par_Bratu2D}
    \end{subfigure}
    \caption{Bifurcation surface for Bratu problem with $p=3$ and different spatial dimension, $d=1$ and $d=2$, left and right respectively.}
    \label{fig:det_3par_Bratu}
\end{figure}

As usual, we start with the Bratu problem in Equation \eqref{bratu_multi}, and for each instance of $\lambda_3$ we apply the zigzag bifurcating curve detection. The points belonging to the bifurcation surface correspond to the ones for which the continuation path crosses between the multiple-solution region and the non-existence one. The plot in Figure~\ref{fig:det_3par_Bratu} illustrates the interesting evolution of the bifurcating curve, for both $d=1$ and $d=2$, while varying the value of the non-homogeneous Dirichlet boundary condition $\lambda_3$. The resulting surface exhibits the non-trivial behavior of the critical curve and the persistence of the saddle-node bifurcation


\begin{figure}[tbp]
    \centering
    \begin{subfigure}[b]{0.45\textwidth}
        \includegraphics[width=\textwidth]{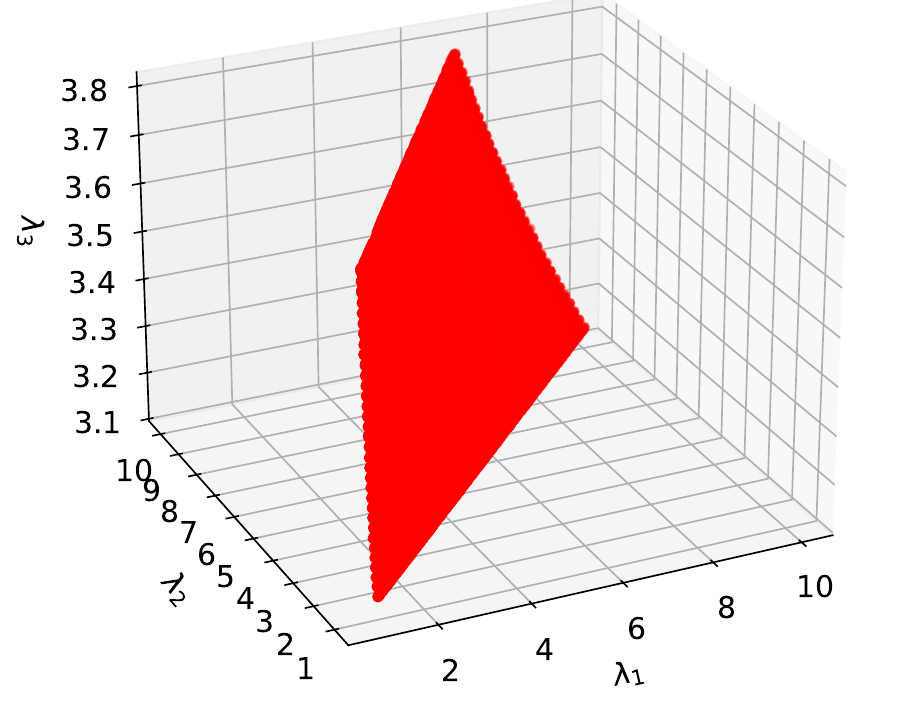}
        \label{fig:det_3par_Allen1D}
    \end{subfigure}
    \quad
    \begin{subfigure}[b]{0.45\textwidth}
        \includegraphics[width=\textwidth]{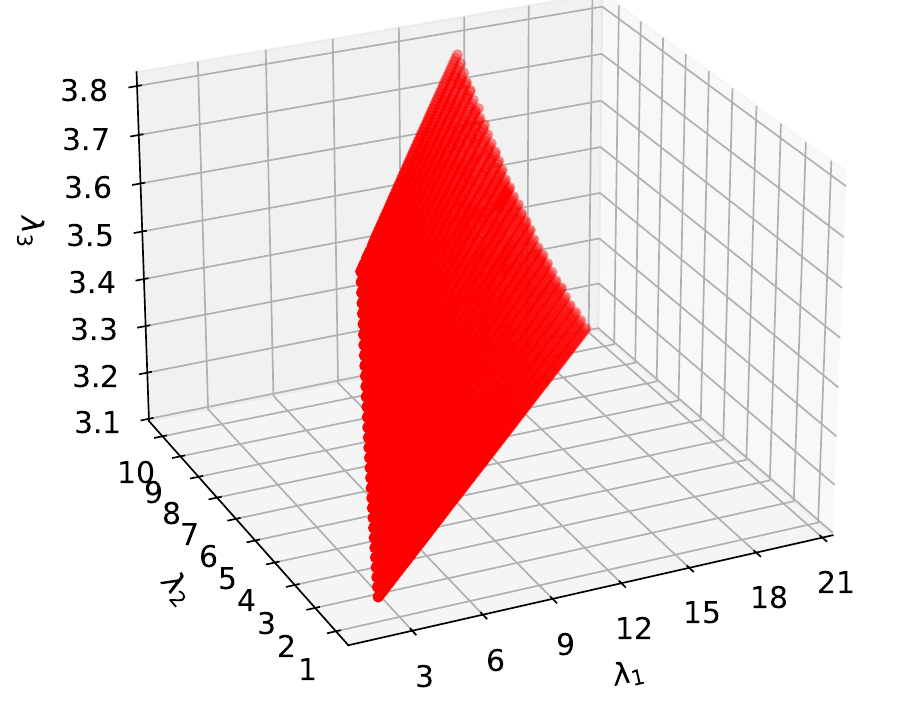}
        \label{fig:det_3par_Allen2D}
    \end{subfigure}
    \caption{Bifurcation surface for Allen--Cahn problem with $p=3$ and different spatial dimension, $d=1$ and $d=2$, left and right respectively.}
    \label{fig:det_3par_Allen}
\end{figure}

Similar comments hold for the multiparametric Allen--Cahn problem in Equation \eqref{allenCahn_mult}, where for each instance of the geometric parameter $\lambda_3$, denoting the length of the interval, we depict the bifurcating surface for $d=1$ and $d=2$, in Figure~\ref{fig:det_3par_Allen}. The robustness of the approach allowed us to obtain also in this case a smooth surface with points detected by the zigzag strategy, indicating the change of stability of the branches and separating the uniqueness and the non-uniqueness parametric regions.

\subsection{Rayleigh-Benard convection}\label{sub_RBC} 
Finally, we consider the well-known Rayleigh–Bénard convection problem \cite{boulle2022bifurcation,herrero2013rb,cortes2025certified}, describing a fluid confined in a cavity and heated from below, with a constant temperature difference imposed along the vertical direction. The fluid density is assumed constant except in the buoyancy term, while the kinematic viscosity $\nu$ and thermal diffusivity $\kappa$ are considered temperature-independent. 

The dimensionless governing equations in $\Omega=[0,1]^d$ with $d=2,3$ are given by
\begin{equation*}
\begin{cases}
-(\mathbf{u}\cdot\nabla)\mathbf{u} 
-\nabla p + Pr\nabla^2 \mathbf{u} + PrRaT\hat{\mathbf{z}}=0,\\
\nabla\cdot\mathbf{u} =0,\\
-(\mathbf{u}\cdot\nabla)T + \nabla^2 T=0,
\end{cases}
\end{equation*}
where $\mathbf{u}$ denotes the velocity field, $p$ the pressure, and $T$ the temperature. The parameters $Ra$ and $Pr$ denote the Rayleigh and Prandtl numbers, respectively, while $\hat{\mathbf{z}}$ represents the vertical unit vector. Homogeneous no-slip boundary conditions are imposed on the velocity field, $\mathbf{u}=0~\text{on }~\partial\Omega,$
and the temperature satisfies $T=1 ~\text{on } z=0,
~T=0 ~\text{on } z=1,$
with thermally insulated lateral boundaries, $\nabla T \cdot \mathbf{n}=0~\text{on } {\mathrm{sidewalls}}.$

The trivial steady state, also called the conducting state, is motionless with a negative thermal gradient through the vertical direction, i.e., $\mathbf{u} =0$ and $T=1-z$.

The spatial discretization is performed using Taylor–Hood elements $P_2$-$P_1$-$P_1$, respectively for velocity, pressure, and temperature approximation. In the $2D$ setting, a structured mesh with $50^2$ cells is considered, resulting in a nonlinear system with $\num{25604}$ degrees of freedom, while for the $3D$ simulation, we employ a $20^3$ mesh, leading to a system with $\num{225285}$ degrees of freedom. Given the significant computational complexity associated with such large-scale benchmark, a robust and efficient strategy to discover bifurcations in the multiparametric context is of utmost importance.

\begin{figure}[tbp]
    \centering
        \centering
        \includegraphics[width=0.24\textwidth,clip, trim=4cm 0cm 4cm 0cm]{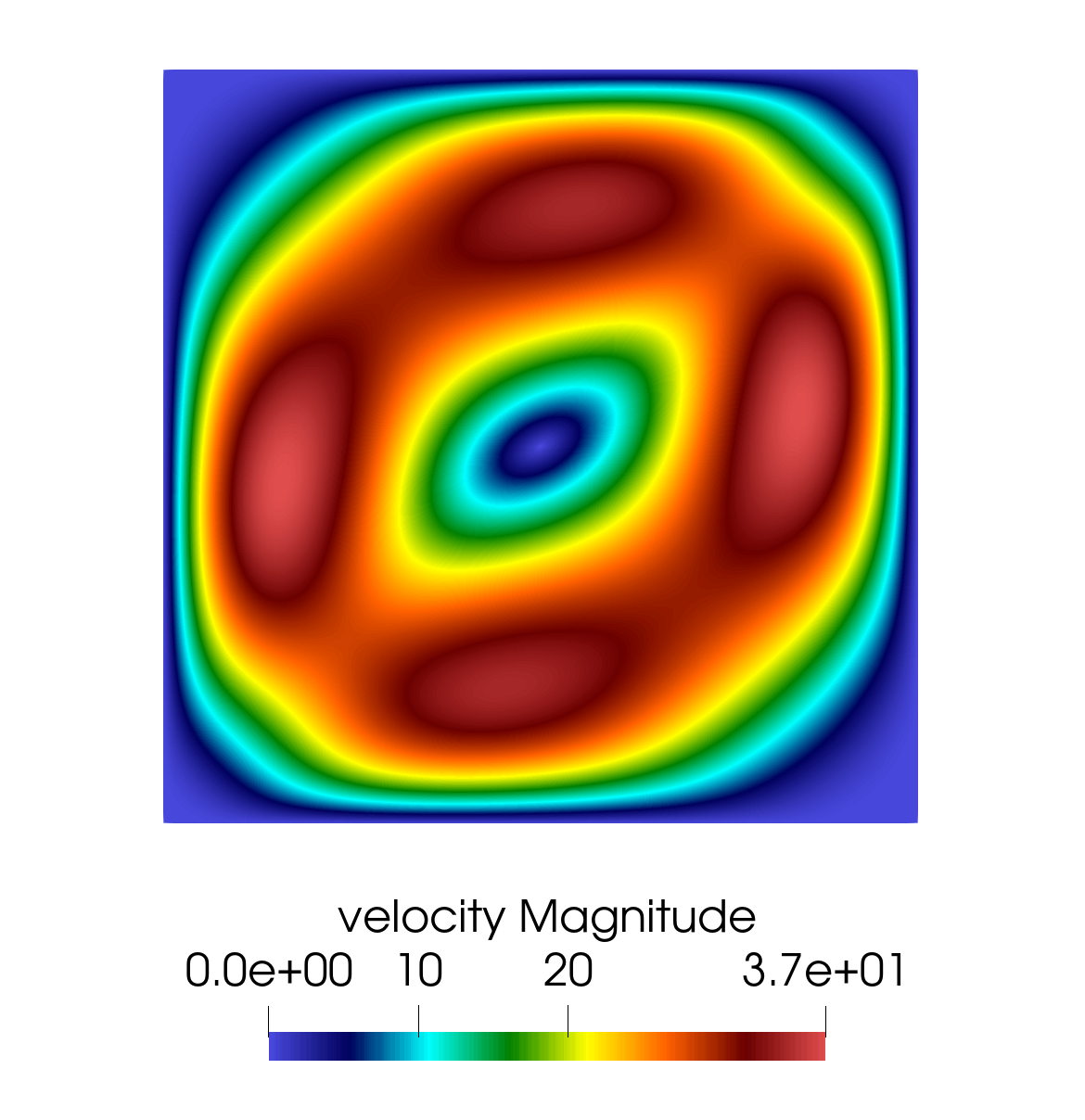}
        \hfill
        \includegraphics[width=0.24\textwidth,clip, trim=4cm 0cm 4cm 0cm]{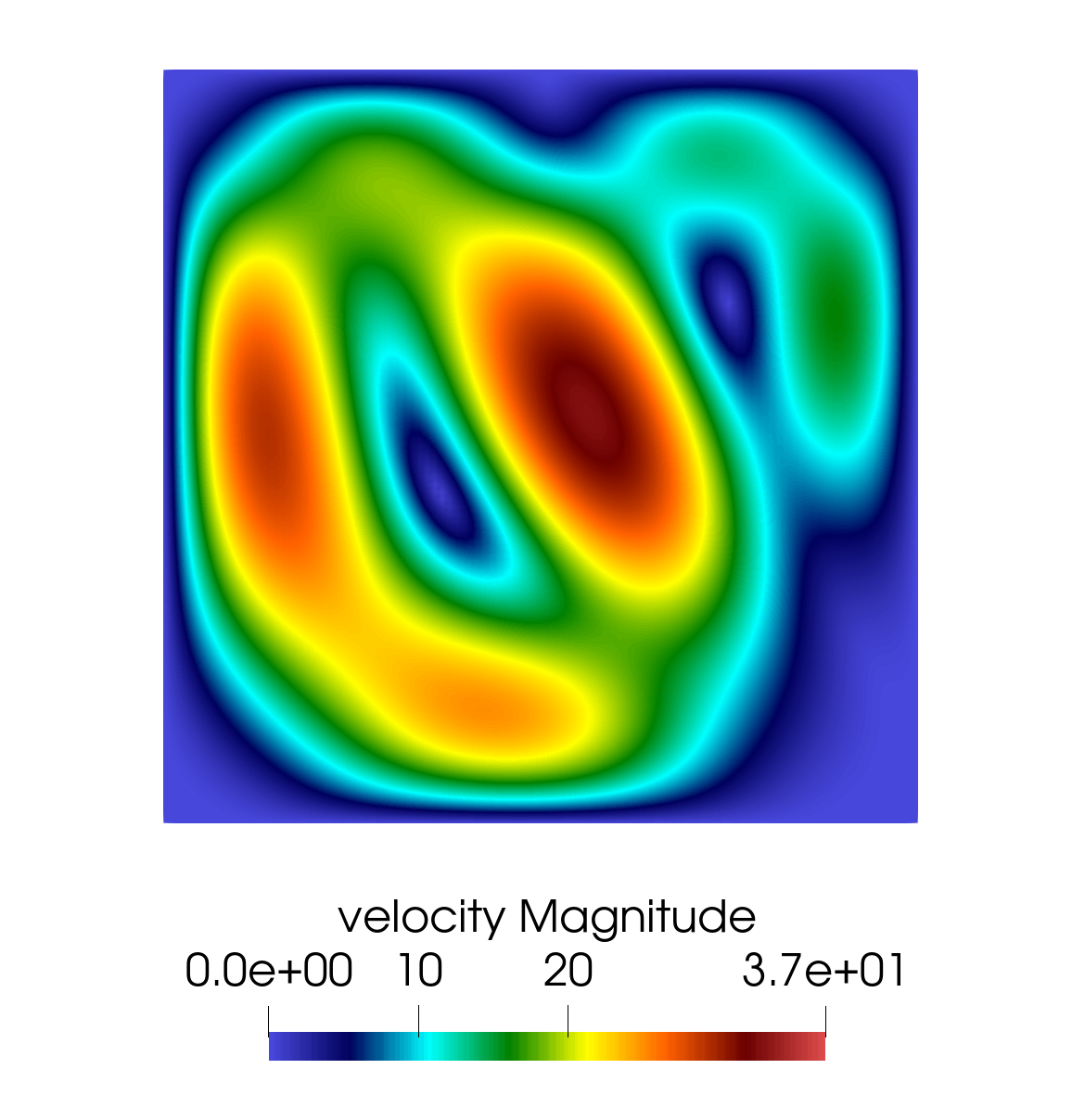}
        \hfill
        \includegraphics[width=0.24\textwidth,clip, trim=4cm 0cm 4cm 0cm]{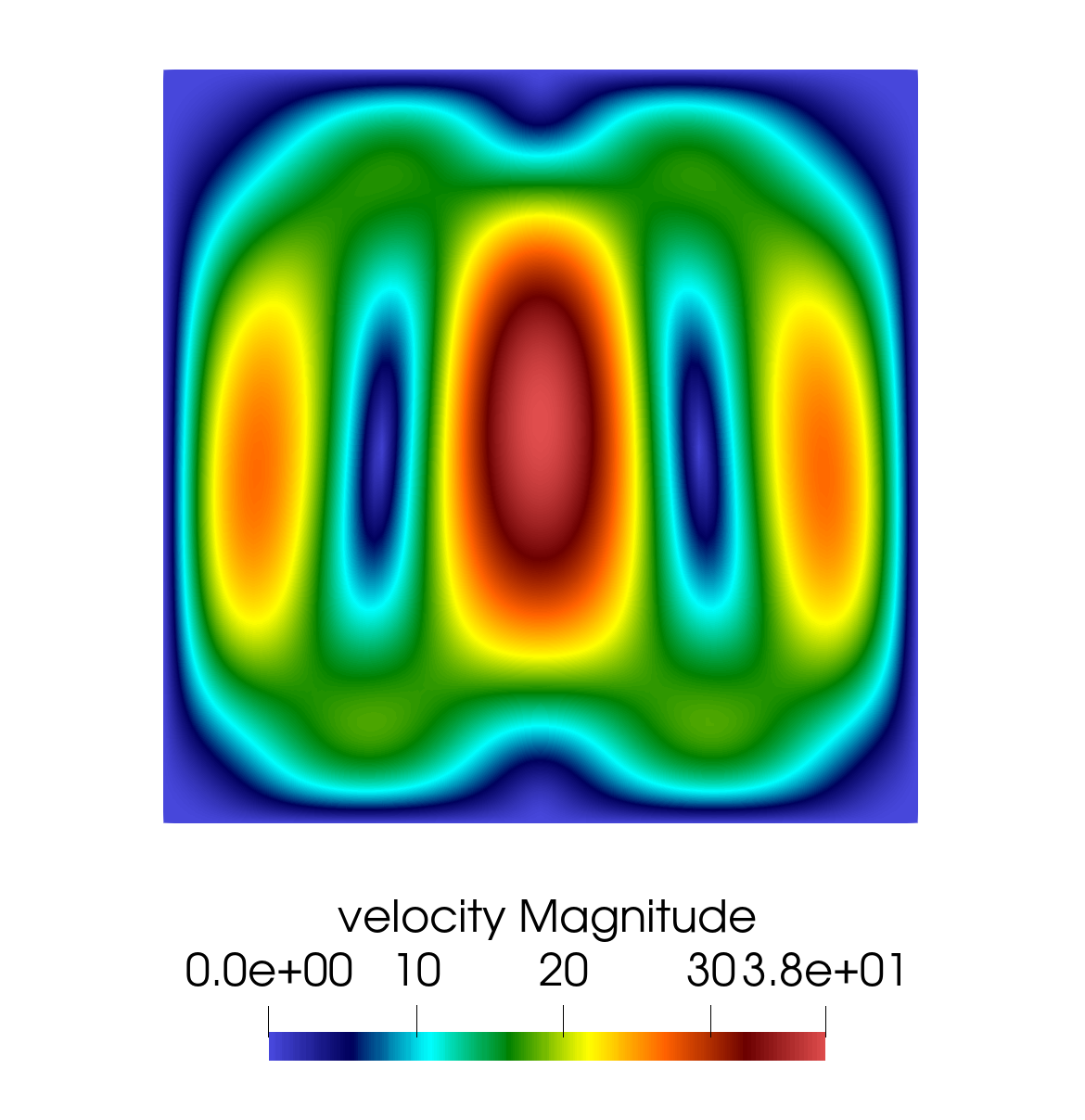}
        \hfill
        \includegraphics[width=0.24\textwidth,clip, trim=4cm 0cm 4cm 0cm]{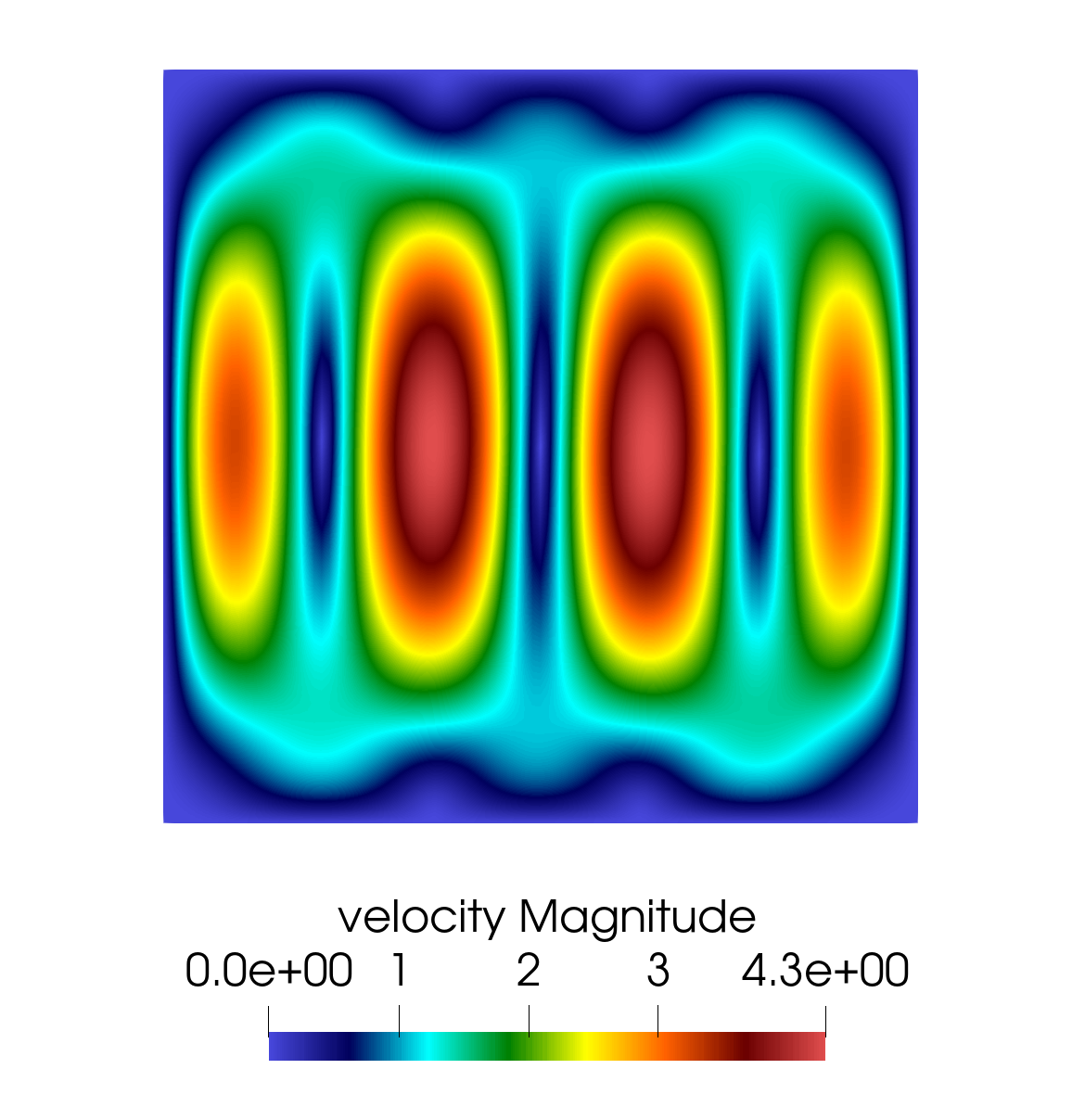}
        
        \includegraphics[width=0.24\textwidth,clip, trim=4cm 0cm 4cm 0cm]{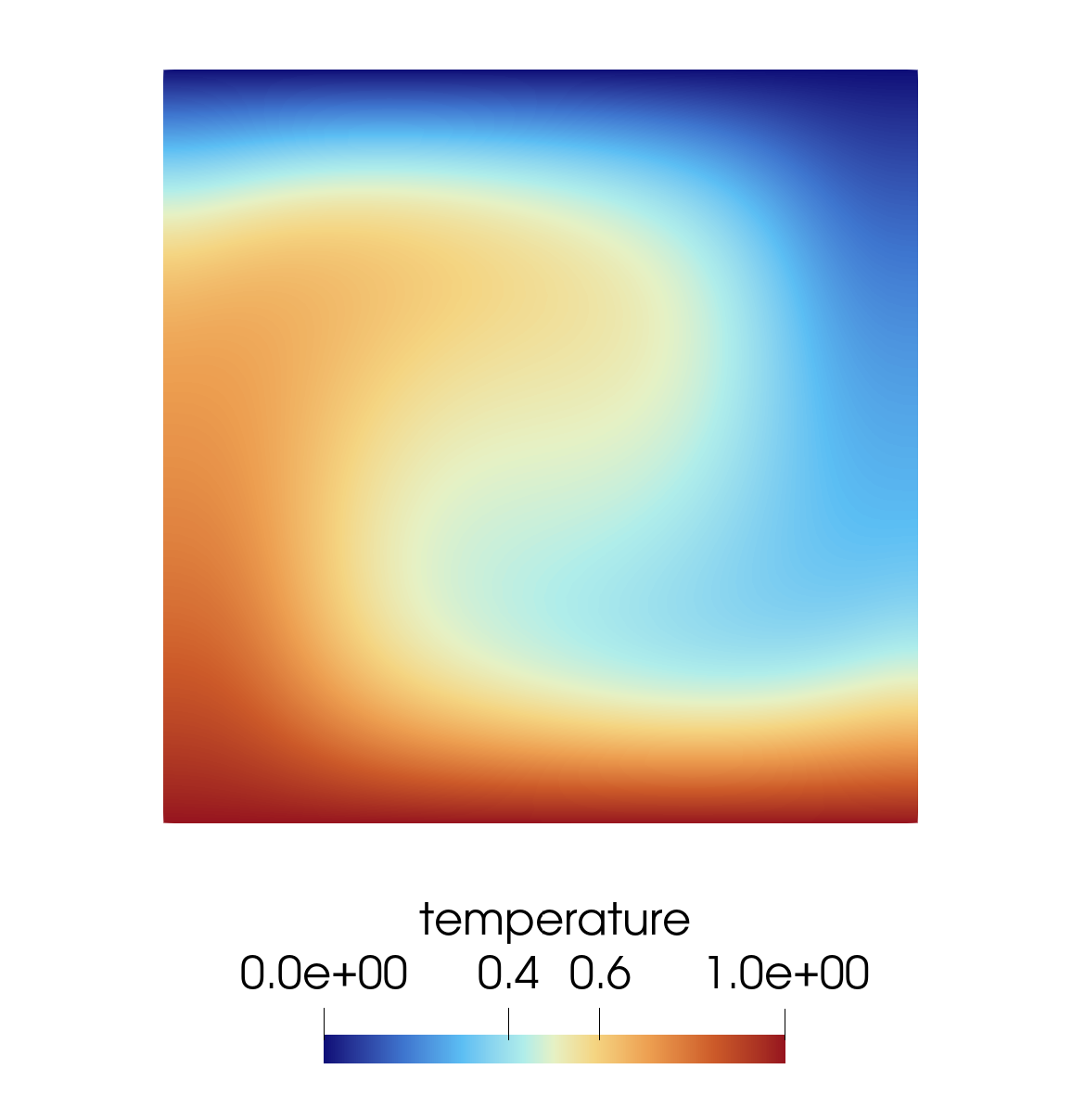}
        \hfill
        \includegraphics[width=0.24\textwidth,clip, trim=4cm 0cm 4cm 0cm]{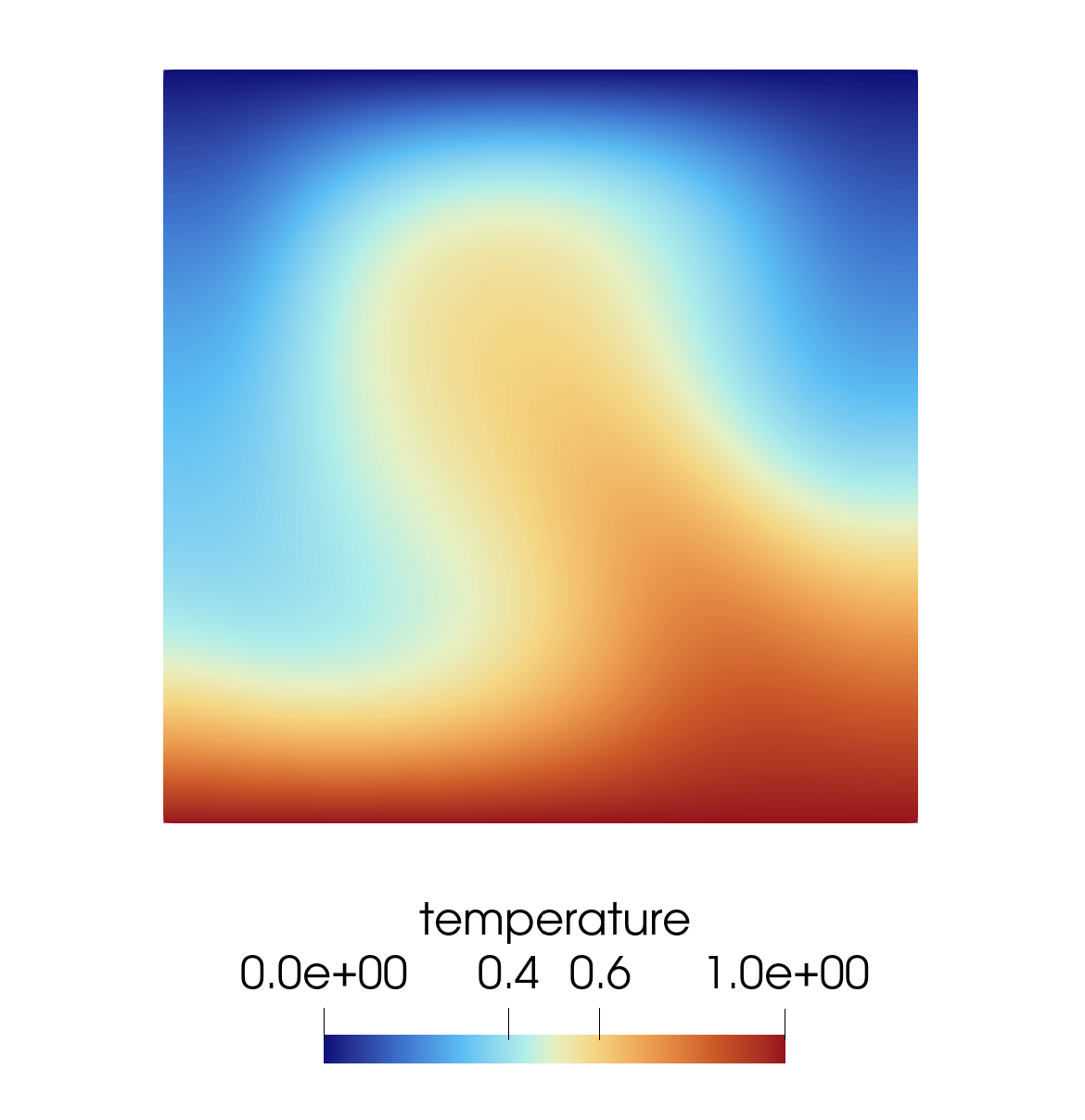}
        \hfill
        \includegraphics[width=0.24\textwidth,clip, trim=4cm 0cm 4cm 0cm]{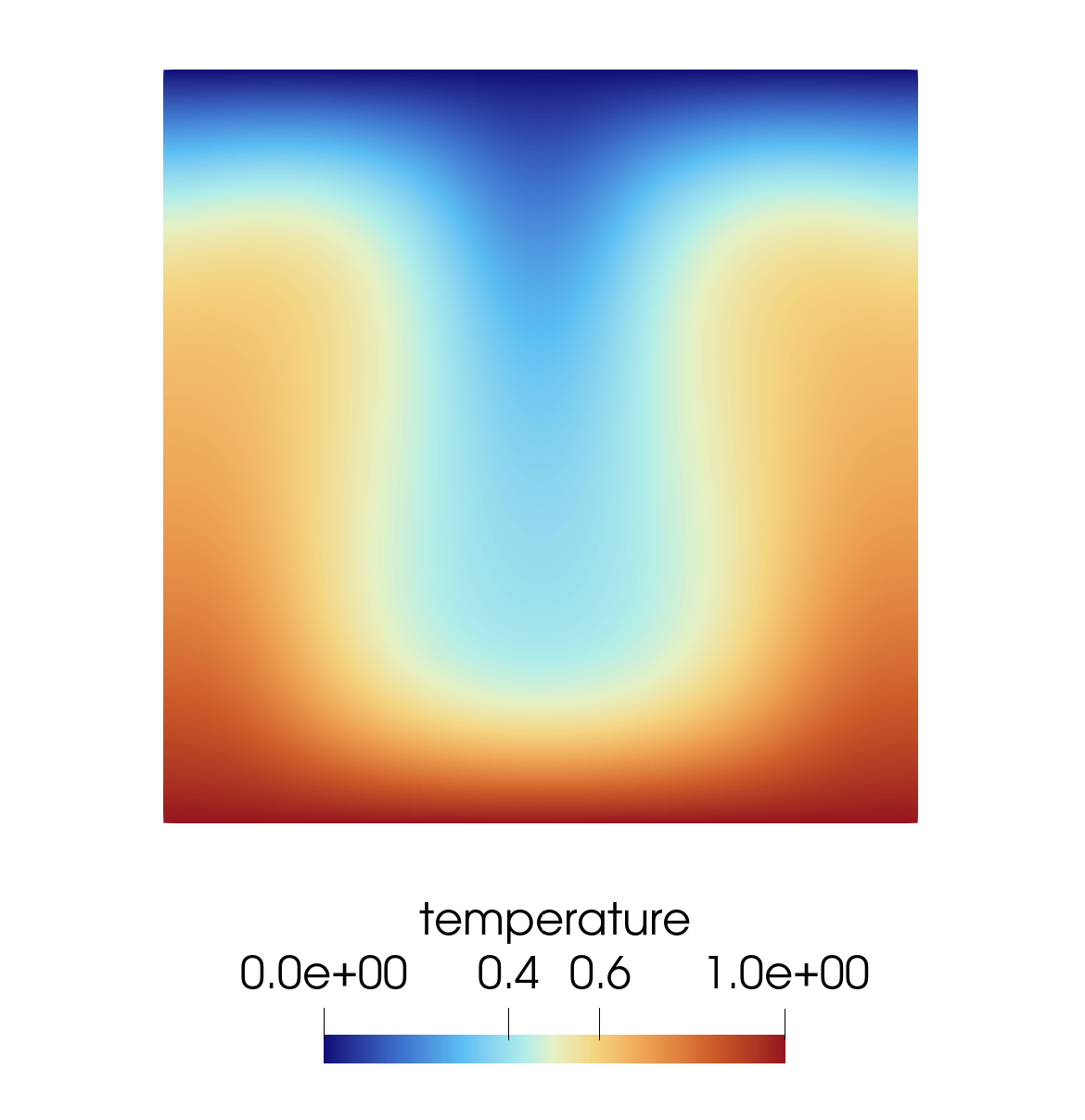}
        \hfill
        \includegraphics[width=0.24\textwidth,clip, trim=4cm 0cm 4cm 0cm]{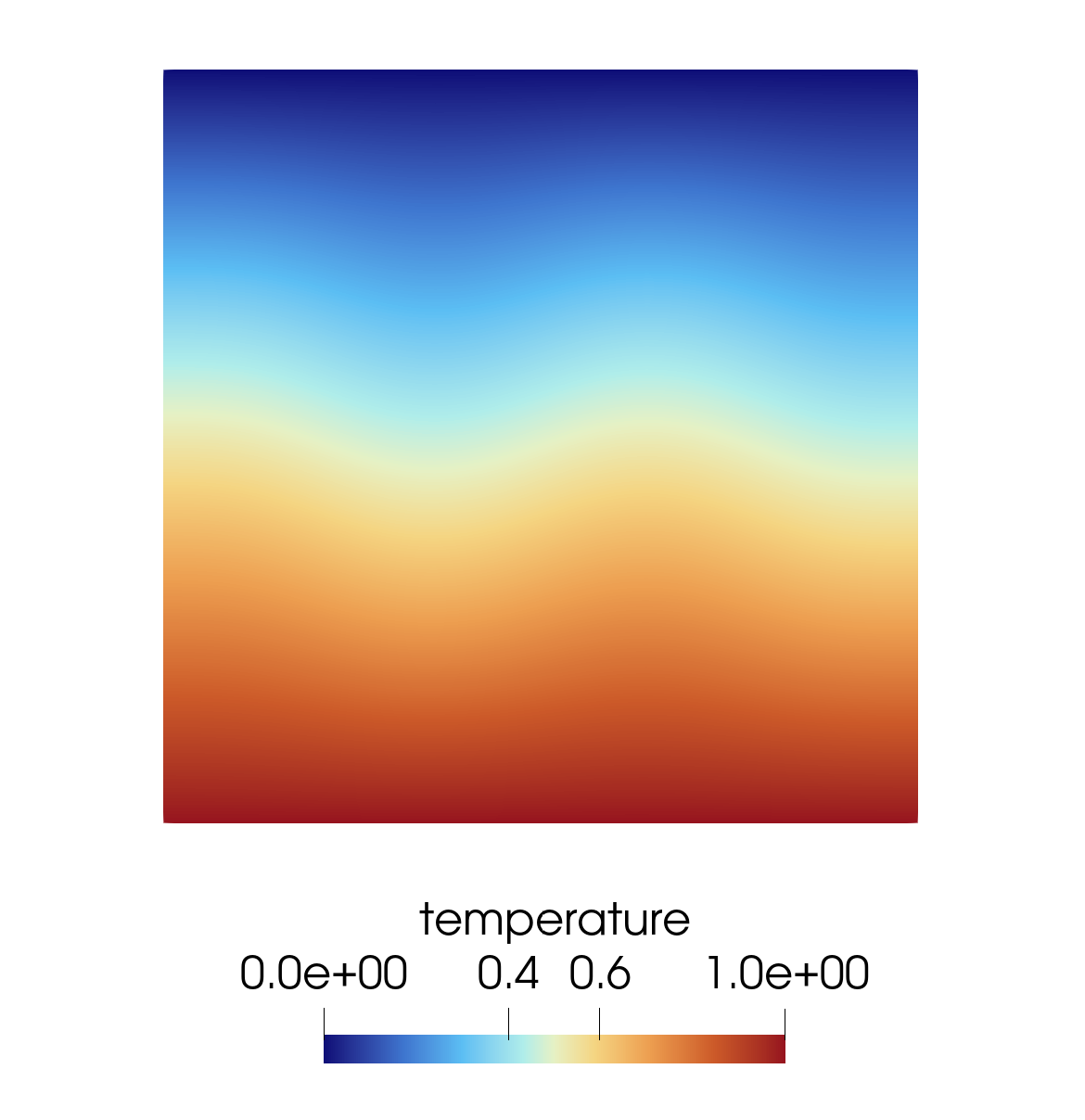}
    \caption{Velocity (upper row) and temperature (lower row) fields for the $2D$ Rayleigh--Bénard convection problem at $Ra\approx 20000$ for $Pr=1$.}
    \label{fig:problemRBC_2D_fields}
\end{figure}

\begin{figure}[tbp]
    \centering
    \begin{minipage}{0.49\textwidth}
        \centering
            \includegraphics[
            width=\textwidth,
            trim=2cm 1.cm 0.4cm 2.5cm,
            clip
        ]{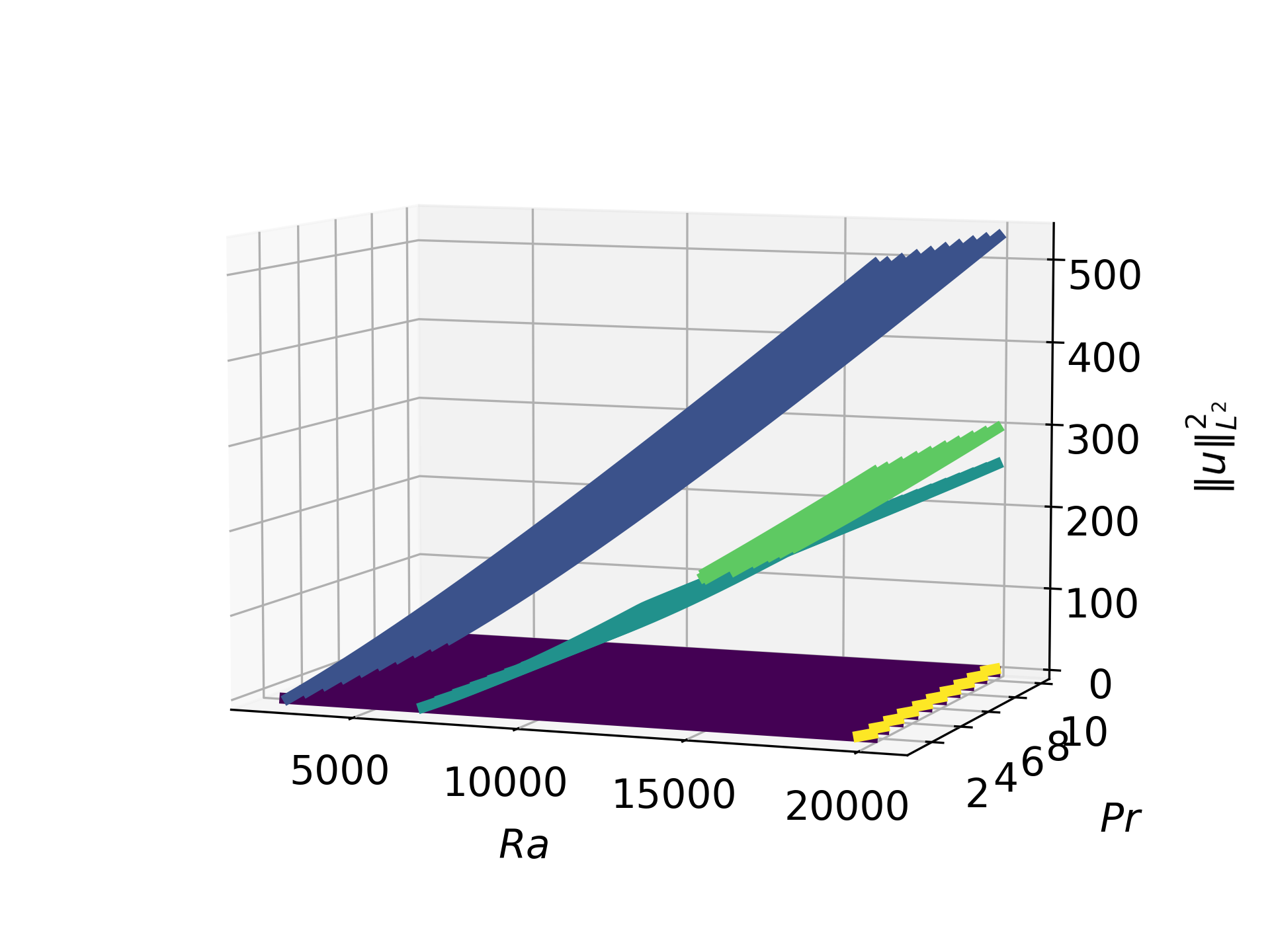}
    \end{minipage}\hfill
    \begin{minipage}{0.45\textwidth}
        \centering
        \includegraphics[
            width=\textwidth,
            trim=1.cm 0cm 0cm 0.cm,
            clip
        ]{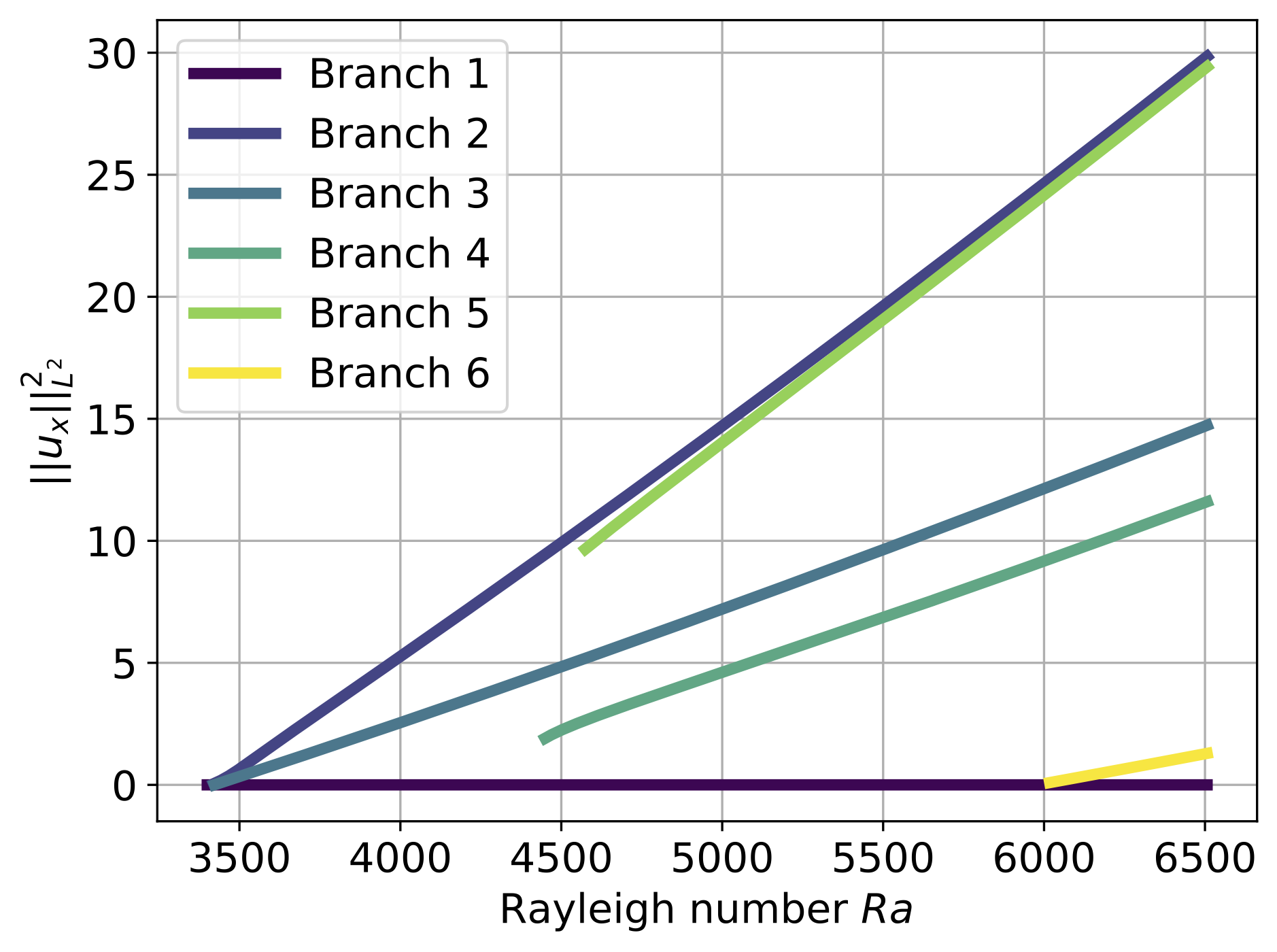}
    \end{minipage}  
    \caption{Bifurcation diagrams for the Rayleigh-Bénard problem in 2D with $p=2$ and in the 3D with $p=1$.}
    \label{fig:problemRBC_2D_diagrams}
\end{figure}

For the $2D$ Rayleigh--Bénard convection problem, we consider the parameter ranges $Ra\in[0,20000]$ and $Pr\in[1,10]$, and construct the bifurcation diagram using straight-line paths in the $(Ra,Pr)$-plane. The proposed methodology reveals three primary bifurcation points at approximately $Ra_c^1\approx\num{2588}$, $Ra_c^2\approx\num{6758}$, and $Ra_c^3\approx\num{19717}$, which, consistently with the linear stability analysis, remain unchanged throughout the considered range of $Pr$. In addition, a secondary bifurcation is observed at $Ra\approx\num{14343}$ from the branch originating from the second primary bifurcation point. Figure \ref{fig:problemRBC_2D_fields} illustrates the velocity and temperature fields at $Ra\approx\num{20000}$ for $Pr=1$, with corresponding bifurcation diagram in the $(Ra,Pr)$-plane depicted in the left plot of Figure \ref{fig:problemRBC_2D_diagrams}. As expected, the solutions exhibit only minor variations with respect to $Pr$, indicating that the overall flow and thermal structures remain qualitatively similar. We remark that, owing to the reflection symmetries of the problem with respect to $x=\tfrac{1}{2}$ and $z=\tfrac{1}{2}$, additional symmetric solution branches have been also obtained, but we have omitted them from the figures for clarity.

\begin{figure}[tbp]
    \centering
        \centering
        \includegraphics[width=0.32\textwidth]{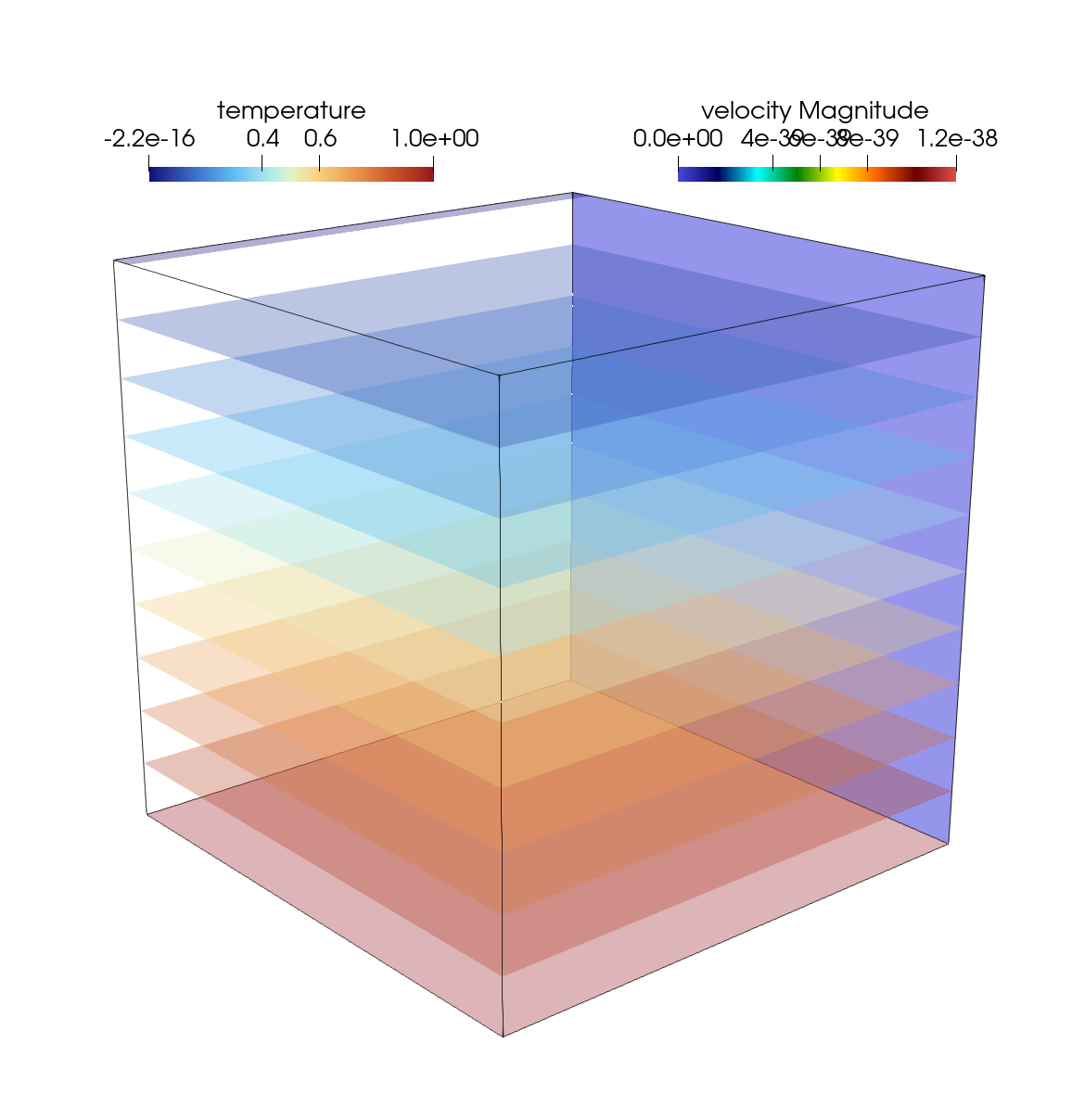}
        \hfill
        \includegraphics[width=0.32\textwidth]{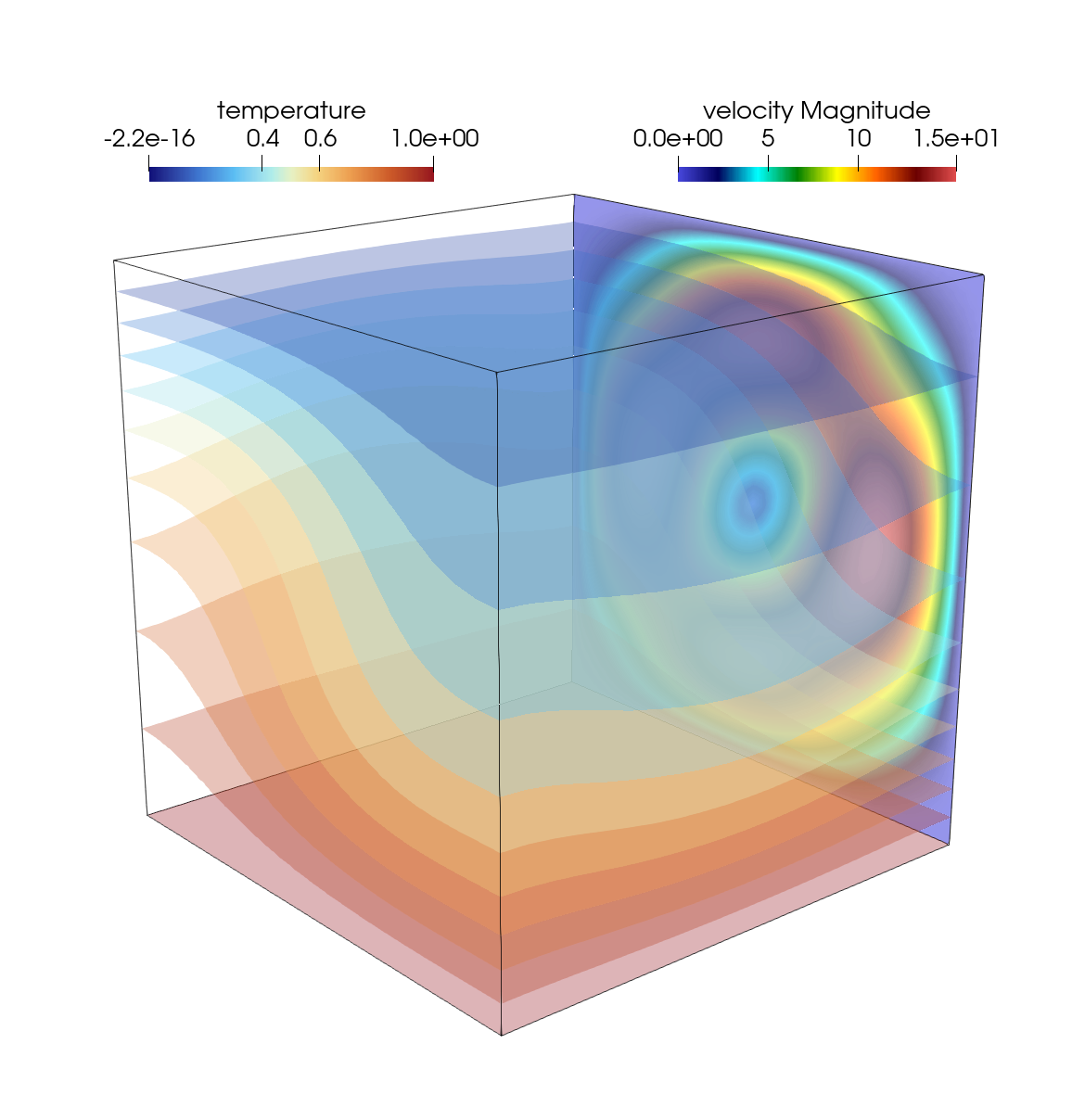}
        \hfill
        \includegraphics[width=0.32\textwidth]{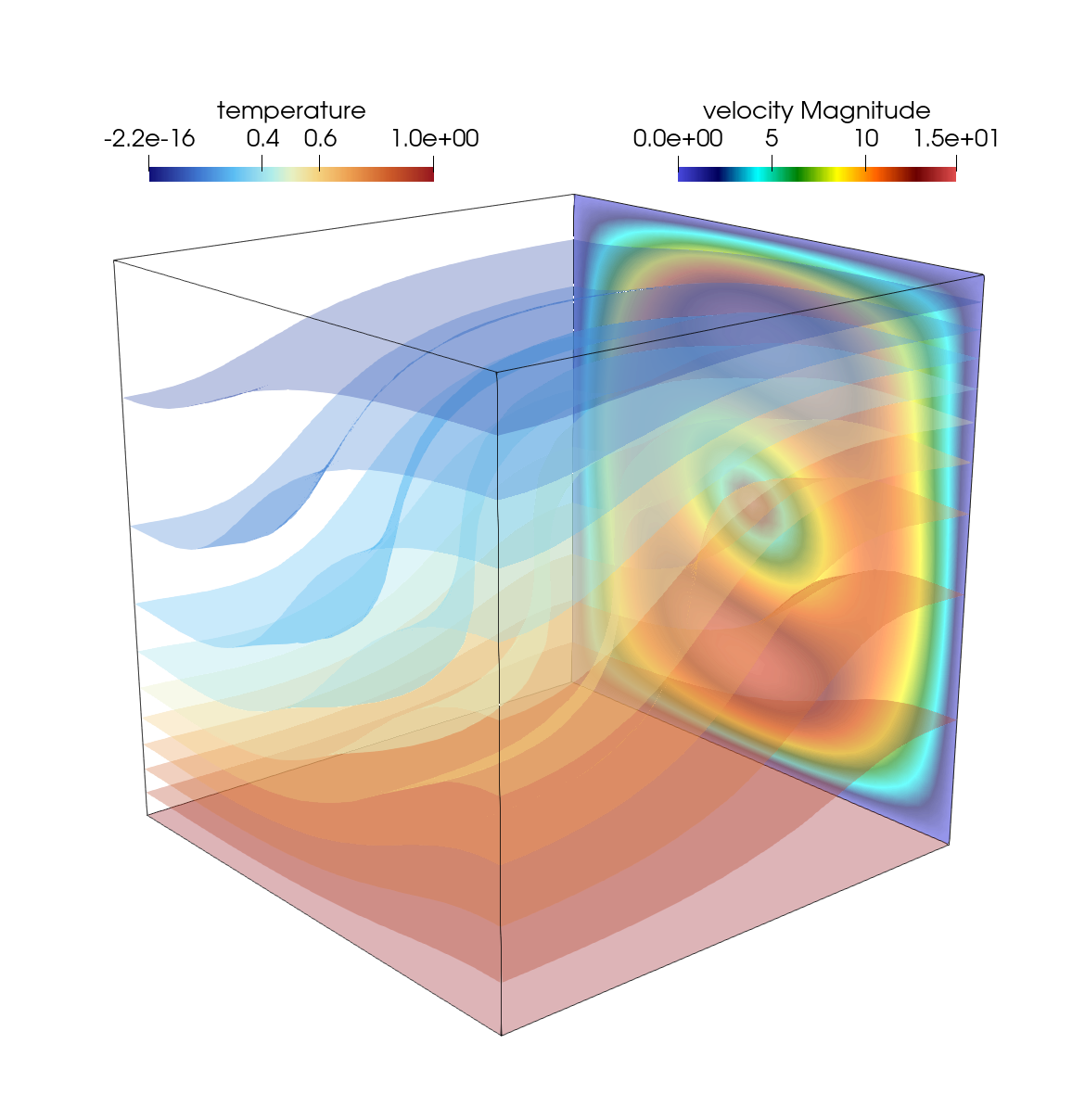}

        \includegraphics[width=0.32\textwidth]{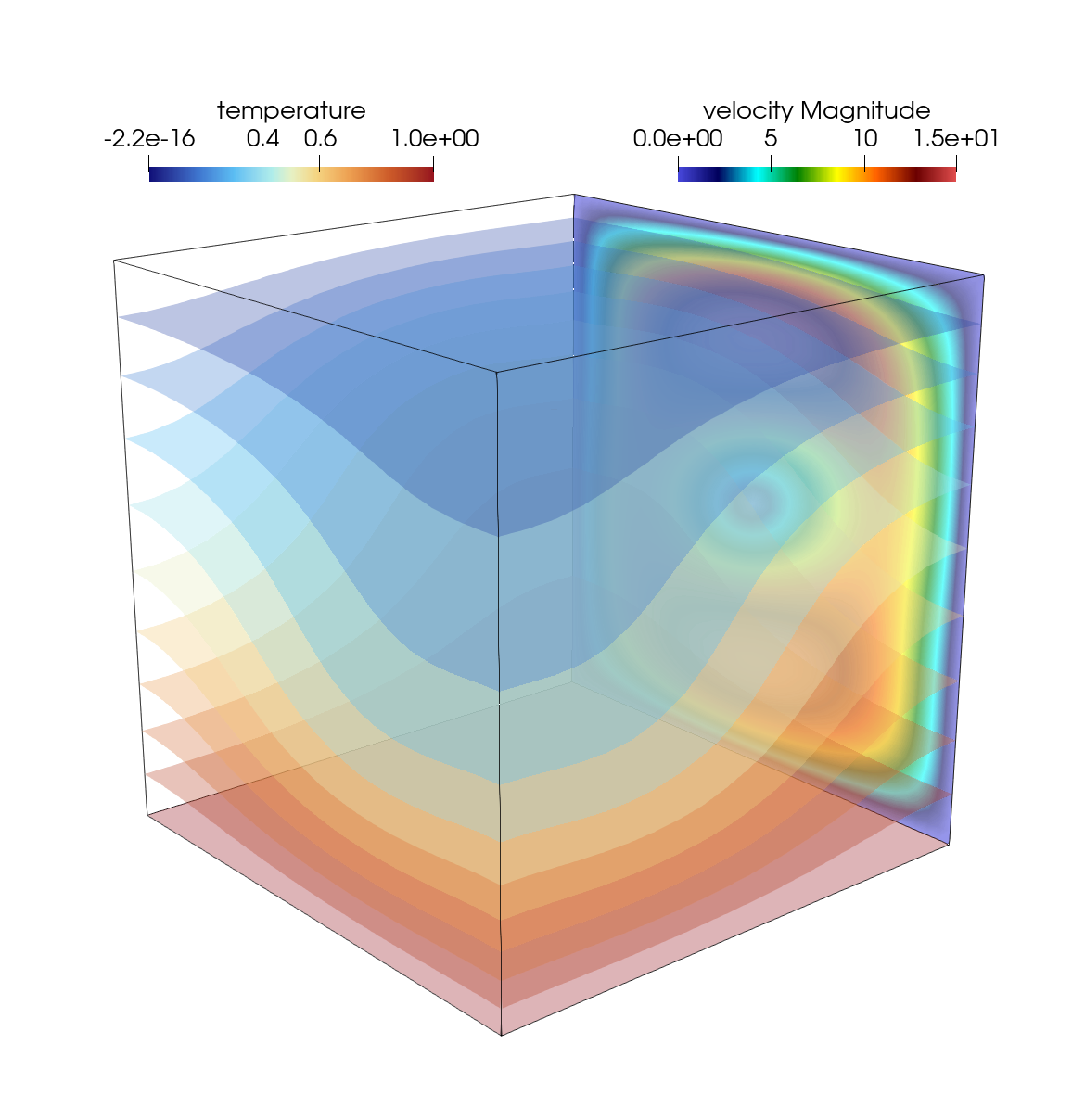}
        \hfill
        \includegraphics[width=0.32\textwidth]{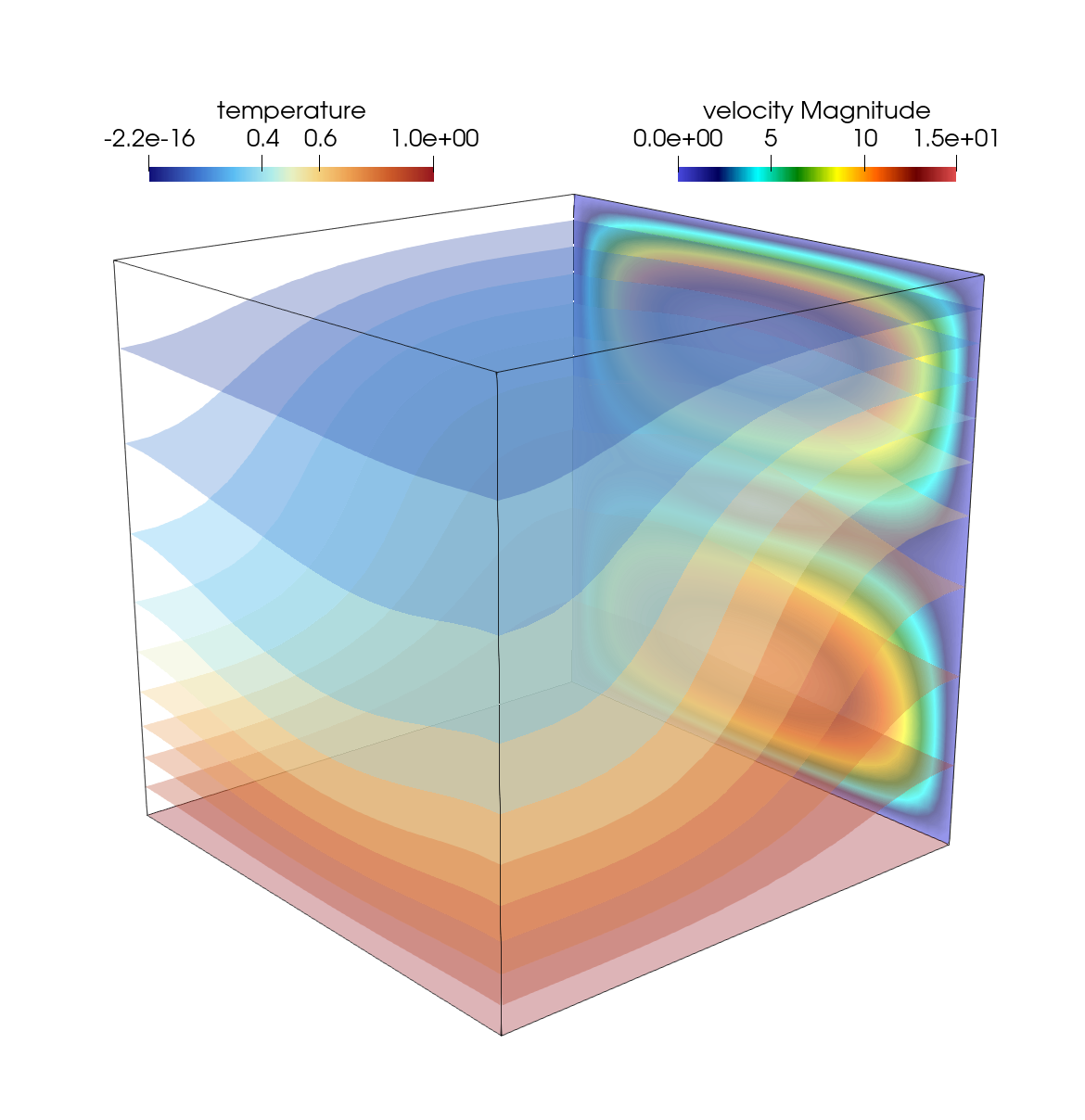}
        \hfill
        \includegraphics[width=0.32\textwidth]{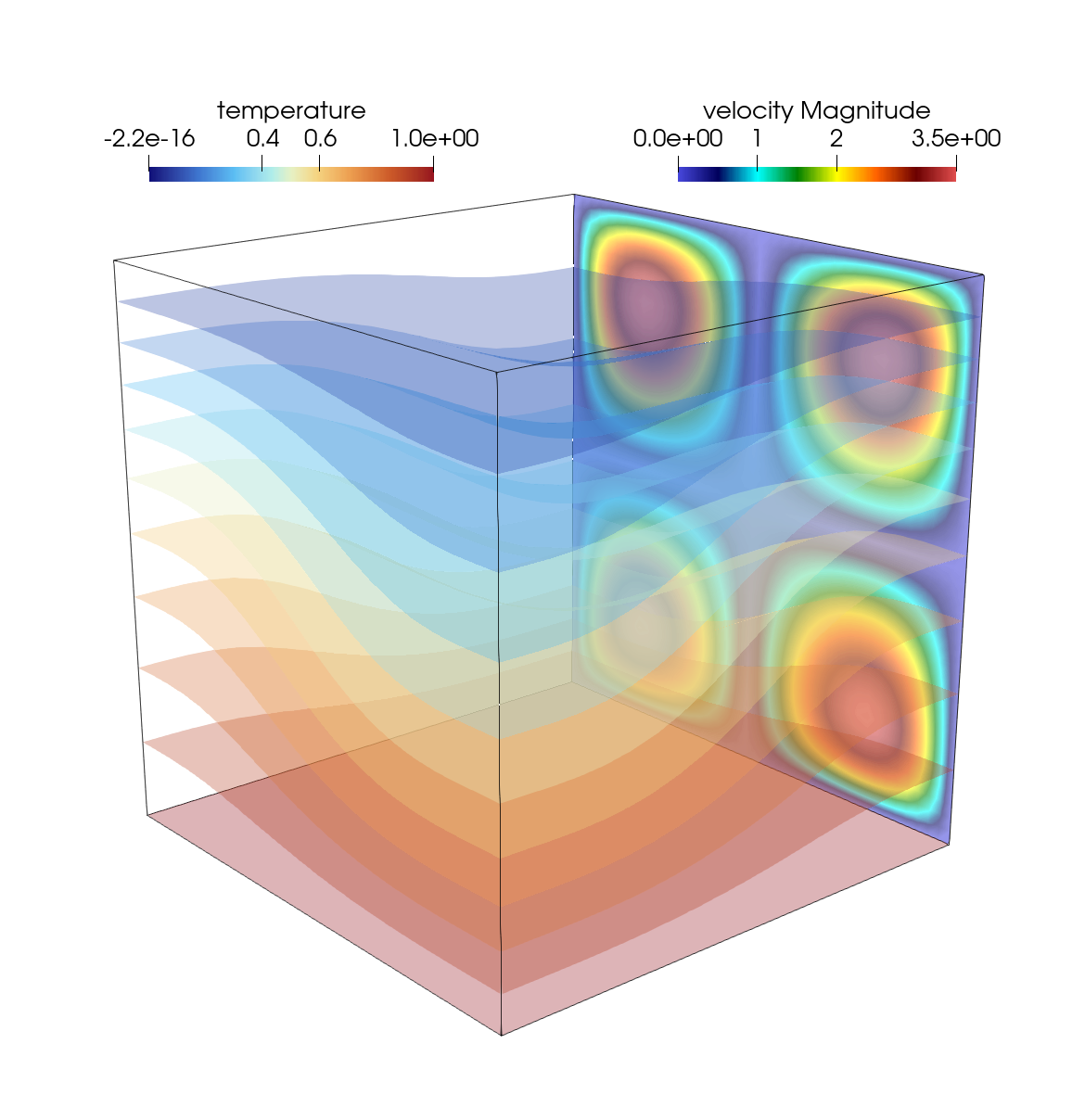}
            \caption{Velocity and temperature fields for the 3D Rayleigh--Bénard convection problem at $Ra\approx 6506$ for $Pr=1$.}
    \label{fig:problemRBC_3D_fields}
    \end{figure}

Finally, for the $3D$ context, we consider the parameter range $Ra\in[3400,6500]$ and compute the bifurcation diagram for the fixed Prandtl number $Pr=1$. The computations reveal three primary bifurcation points at approximately $Ra_c^1\approx\num{3410}$, $Ra_c^2\approx\num{3421}$ and $Ra_c^3\approx\num{5978}$. Furthermore, secondary bifurcations are detected at approximately $Ra_s^1\approx\num{4575}$ and $Ra_s^2\approx\num{4450}$, originating from the branches emerging from $Ra_c^1$ and $Ra_c^2$, respectively. The right plot in Figure  \ref{fig:problemRBC_2D_diagrams} depicts the bifurcation diagram, while Figure \ref{fig:problemRBC_3D_fields} presents, for each of the 6 different coexisting branches, a 2D slice of the velocity magnitude and the isocontours for the temperature fields. 
Similarly to the $2D$ setting, additional qualitatively similar solution branches are obtained also here due to the symmetry of the governing equations.


%% file: sections/conclusion.tex
\section{Conclusion}
In this work, we proposed a deflated arclength continuation technique for the reconstruction of bifurcation diagrams and the detection of bifurcation curves in multiparametric nonlinear PDEs. Unlike simple continuation methods, arclength continuation offers significant advantages by adaptively adjusting the step size in parameter space, taking smaller steps near bifurcation points where sharp changes occur, and larger steps in smoother regions, thus making the method more robust and reliable.

Moreover, we extended the arclength continuation framework to multiparametric settings by introducing predefined paths in the parameter space, and the integration of this with the deflation techniques led to the development of the multiparametric deflated arclength continuation method, which enables efficient and comprehensive reconstruction of bifurcation diagrams in higher-dimensional parameter spaces. This method effectively captures solution branches that emerge from different bifurcation curves. The approach was validated via three benchmark problems: the Bratu equation, exhibiting saddle-node bifurcations, the Allen--Cahn equation, featuring pitchfork bifurcations, and the Rayeigh--Benard convection problem, showing a rich bifurcation structure with complex velocity/temperature patterns. Our results demonstrate the robustness of the method, particularly in its ability to detect multiple bifurcation curves.

In addition, we introduced a zigzag path-following strategy to accurately detect bifurcation curves in multiparametric settings. The proposed strategy defines a path crossing the bifurcation curve with a small angle $\theta$, effectively and efficiently tracking the critical region, even for three-dimensional parameter spaces, showcasing its scalability and practical utility.

A promising direction for future research building on this work involves integrating reduced order modeling (ROM) techniques to further enhance computational efficiency. Indeed, having identified the bifurcation curves, one can employ ROM methods to locally reconstruct the bifurcation diagrams within each distinct region at a much lower computational cost, especially for large-scale problems. Additional challenges to be addressed in the following works include: the presence of multiple bifurcation curves, the treatment of time-dependent problems, and the analysis of more challenging bifurcating PDEs coming from continuum mechanics.

%% file: sections/acknowledgments.tex
\section*{Acknowledgments}
The authors acknowledge the support provided by INdAM-GNCS (CUP E53C25002010001) and the European Union - NextGenerationEU, in the framework of the iNEST - Interconnected Nord-Est Innovation Ecosystem (iNEST ECS00000043 - CUP G93C22000610007) consortium and its CC5 Young Researchers initiative. 

%% file: sections/Appendix.tex
\appendix
\section{Wavelet collocation method}\label{appendix1}
We briefly outline here the discretization approach employed in this work for solving the Bratu and Allen--Cahn equations. We exploited the Legendre wavelet-based approach (see \cite{kumar2024generalized,ray2018wavelet,mehra2018wavelets} for more details), primarily for its efficiency and accuracy in representing smooth solutions, but it is important to emphasize that the proposed continuation and bifurcation framework is independent of the specific discretization technique adopted.

The approximation space is constructed using Legendre scaling functions, which are obtained by scaling and translating the Legendre polynomials \( P_m(x) \) as
\begin{eqnarray*} \phi_{k,m}^J(x) = \begin{cases} 2^{J/2}\sqrt{\frac{2m+1}{2}}P_{m}(2^{J}x-\hat{k}) & ~~ \frac{(k-1)}{2^{J-1}} \leq x < \frac{k}{2^{J-1}}, \\ 0 & ~~ \text{otherwise,} \end{cases} \end{eqnarray*} where $J=1,2,\ldots$ is the level of resolution, $\hat{k}=2k-1 ~ \text{with} ~ k=1,2,\ldots,2^{(J-1)},$ is the translation parameter, and $m=0,1,2,\ldots,M-1,$ is the degree of Legendre polynomial. Then, the highest-order derivative in the governing PDE is expanded as
\begin{equation*}
u_{xx}(x,y) \approx 
\sum_{k_1=1}^{2^{J_1-1}} \sum_{m_1=0}^{M_1-1} 
\sum_{k_2=1}^{2^{J_2-1}} \sum_{m_2=0}^{M_2-1}
c_{k_1,m_1}^{k_2,m_2}\,
\phi_{k_1,m_1}^{J_1}(x)\,\phi_{k_2,m_2}^{J_2}(y)
= \big(\Phi_{K_1}^{J_1}(x)\big)^T C^T \Phi_{K_2}^{J_2}(y),
\end{equation*}
where $ c_{k_1,m_1}^{k_2,m_2} $ are the unknown coefficients to be determined, $(x,y)\in \Omega \subset \mathbb{R}^2$, and $K_i = 2^{J_i-1} M_i$ for $i = 1,2$. Furthermore, $\Phi_{K_1}^{J_1}(x)$ and $\Phi_{K_2}^{J_2}(y)$ represent the vector forms of the scaling functions $\phi_{k_1,m_1}^{J_1}(x)$ and $\phi_{k_2,m_2}^{J_2}(y)$, respectively. Similarly, we approximate $u_{yy}(x,y)$ with the unknown coefficients $ d_{k_1,m_1}^{k_2,m_2} $ and
by integrating the above relations for $u_{xx}, u_{yy}$ and enforcing the boundary conditions, approximations for the lower-order derivatives $u_x$, $u_y$, and $u$ are obtained. 

Finally, substituting these approximations into the parametric PDE and applying a collocation procedure at selected grid points yields the desired nonlinear algebraic system of equations. This system is then solved using the arclength continuation and deflation strategies as described in the manuscript to construct the bifurcation diagram.